\documentclass{article}

\usepackage{graphicx}
\usepackage[utf8]{inputenc}
\usepackage[margin=0.8in]{geometry} 
\usepackage{mathtools,amsmath,amsthm,amssymb,amsfonts}
\usepackage[backend=biber, style=alphabetic]{biblatex}
\usepackage{xcolor}
\usepackage{hyperref}
\usepackage{bbm}
\usepackage{enumitem}
\usepackage{blindtext}
\usepackage{hyperref}
\usepackage{caption}
\usepackage{subcaption}
\usepackage{float}
\usepackage{enumitem} 
\usepackage[final]{showlabels}
\usepackage{esdiff}

\DeclareMathAlphabet{\pazocal}{OMS}{zplm}{m}{n}
\DeclareFontFamily{T1}{calligra}{}
\DeclareFontShape{T1}{calligra}{m}{n}{<->s*[1.44]callig15}{}
\DeclareMathAlphabet{\mathcalligra}{T1}{calligra}{m}{n}
\DeclareMathOperator*{\esssup}{ess\,sup}

\makeatletter
\newcommand*{\rom}[1]{\expandafter\@slowromancap\romannumeral #1@}
\newcommand{\Exp}{\mathbb{E}}
\newcommand{\N}{\mathbb{N}}
\newcommand{\R}{\mathbb{R}}

\newcommand{\Z}{\mathbb{Z}}

\newcommand{\Prob}{\mathbb{P}}

\newcommand{\Div}{\mathrm{div}}
\newcommand{\W}{\mathcal{W}}

\newcommand{\bmu}{\boldsymbol{\mu}}
\newcommand{\brho}{\boldsymbol{\rho}}
\newcommand{\ind}{\mathbf{1}}
\newcommand{\footremember}[2]{%
    \footnote{#2}
    \newcounter{#1}
    \setcounter{#1}{\value{footnote}}%
}
\providecommand{\keywords}[1]
{
  \small	
  \textbf{Keywords: } #1
}
\providecommand{\classification}[1]
{
  \small	
  \textbf{AMS Subject classification (2020): } #1
}

\newtheorem{theorem}{Theorem}[section]
\newtheorem{proposition}{Proposition}[section]
\newtheorem{lemma}{Lemma}[section]
\theoremstyle{definition}
\newtheorem{definition}{Definition}[section]
\theoremstyle{definition}
\newtheorem{remark}{Remark}[section]

\newtheorem{example}{Example}[section]
\newtheorem{def-theorem}{Definition-Theorem}[section]

\numberwithin{equation}{section}

\title{On well-posedness of stable-driven McKean-Vlasov stochastic differential equations with Besov interaction kernel of negative regularity}
\author{Anna Bahrii \footremember{alley}{Laboratoire de Modélisation Mathématiques d’Evry (LAMME), UMR CNRS 8071, Université d’Evry Val d’Essonne Paris Saclay, 23 Boulevard de France, 91037 Evry, France. Email: \href{mailto:anna.bahrii@univ-evry.fr}{anna.bahrii@univ-evry.fr}}}
\date{}

\begin{document}

\maketitle

\begin{abstract}
	We prove well-posedness results for time-inhomogeneous stable-driven McKean-Vlasov stochastic differential equations with a convolution drift where the interaction kernel belongs to some Lebesgue-Besov space. The novelty of this work is that we manage to go below -1 in space regularity for such a kernel. This is achieved under additional smoothness conditions on the initial data and divergence conditions on the kernel. The proof heavily relies on a suitable product rule in Besov space. We prove smoothing properties of the law, which allow us to have the drift in a Lebesgue-Besov space of non-positive regularity.
\end{abstract}


\classification{60H10, 60H50, 35Q84, 30H25, 60G52}

\tableofcontents

\section{Introduction}

\subsection{Problem formulation and literature overview}
Let $T>0$ be a fixed time horizon. For $s\in[0,T]$, consider the \emph{formal} McKean-Vlasov stochastic differential equation (SDE)
\begin{align}
	\label{main-sde}
	\begin{split}
		dX_s^{\mu}=b\ast\bmu_s^{\mu}(s,X_s^{\mu})ds+d\W_s,\quad X_0^{\mu}\sim\mu\in\mathcal{P}(\R^d),
	\end{split}
\end{align}

where $\mathcal{P}(\R^d)$ denotes the set of probability measures on $\R^d$, $(\W_s)_{s\geq 0}$ is a symmetric $d$-dimensional $\alpha$-stable process with $\alpha\in(1,2]$, $d\geq1$, $\boldsymbol{\mu}^{\mu}_s=\text{Law}(X_s^{\mu})$, and $\mu\in\mathcal{P}(\R^d)$ is the initial distribution of the process. We are interested in the well-posedness of equation \eqref{main-sde} when the interaction kernel $b$ belongs to a Lebesgue-Besov space with negative regularity in space and the initial probability measure $\mu$ belongs to a Besov space with non-negative regularity. Precisely, we assume that
\begin{align*}
	b\in L^r((0,T],B^{\beta}_{p,q}(\R^d,\R^d))=:L^r(B^{\beta}_{p,q}), \quad \beta<0,\ p,q,r\in[1,+\infty],
\end{align*}
and
\begin{align*}
	\mu\in \mathcal{P}(\R^d)\cap B^{\beta_0}_{p_0,q_0}\cap B^{\bar{\theta}}_{1,\infty}, \quad \beta_0\geq0,\ p_0,q_0\in[1,+\infty], \quad \bar{\theta}>0.
\end{align*}

We refer to Section \ref{section-tools} for the definition and properties of Besov spaces. Moreover, we are interested in considering $b$ and $\mu$ such that the convolution drift $b\ast\bmu_s^{\mu}$ belongs to some Lebesgue-Besov space with non-positive regularity in space denoted by $-\theta\in(-1/2,0]$. When $-\theta\in(-1/2,0)$, the considered equation \eqref{main-sde} is \emph{formal} as the drift is only assumed to be a distribution, and thus the integrated drift is not well-defined in a classical sense. This leads to a modified notion of martingale and weak solutions, similarly to the one considered in \cite{deraynal2022-linear-besov}, see Definition \ref{def-martingale-problem} and Definition \ref{def-weak-solution}.

\vspace{1mm}

For ``classic" SDEs, i.e., when the drift does not depend on the law, well-posedness with the driving noise being Brownian motion (i.e., $\alpha=2$) and initial measure being only a probability measure, has been thoroughly studied in the literature. In the seminal works \cite{zvonkin} ($d=1$), \cite{Veretennikov1981} ($d>1$), the authors proved strong well-posedness of such an SDE when the drift is bounded and Borel measurable. In \cite{KrylovRockner2005LpLq}, the authors established a strong well-posedness result for the drift belonging to $L^q_t-L^p_x$ space, where it holds $d/p+2/q<1$. For an $\alpha$-stable driven SDE for $\alpha\in(1,2)$, we mention \cite{XieZhang2020LpLq}, where the Krylov-Röckner condition is extended to $d/p+\alpha/q<\alpha-1$ for ensuring strong well-posedness under some additional conditions. More precisely, $b$ has to be in $L^q-H^\beta_{p}(\R^d)$ with $\beta>1-\alpha/2$, where $H^\beta_{p}(\R^d)=(I-\Delta/2)^{\beta/2}(L^p(\R^d))$ is the Bessel potential space. In \cite{fitoussi2024weakwellposednessweakdiscretization}, the authors prove weak well-posedness under the condition $d/p+\alpha/q<\alpha-1$. We also mention \cite{rockner2025sdescriticaltimedependent} where a strong well-posedness result is proven in the critical Brownian case $d/p+2/q=1$, and \cite{krylov21pd} where a strong well-posedness is proven for a time-homogeneous drift condition to $d/p=1$.

For a Brownian motion driven SDE with a distributional drift (still measure-independent), well-posedness was studied in \cite{flandoli2015multidimensionalsde}. Therein, the authors established the existence and uniqueness of so-called virtual solutions when the drift belongs to a suitable fractional Sobolev space. In \cite{Priola2012PathwiseUniq}, the author established pathwise uniqueness for an SDE with a Hölder drift of regularity $\beta>1-\alpha/2$ with $\alpha\in(1,2]$. In \cite{ABM2020StableSDE}, the authors extended the previously mentioned result by proving strong well-posedness of a $1$-dimensional SDE driven by a symmetric $\alpha$-stable noise when the time-homogeneous drift belongs to the Besov-Hölder space $B^{\beta}_{\infty,\infty}$ with $\beta>(1-\alpha)/2$. Finally, in \cite{deraynal2022-linear-besov}, the authors proved well-posedness of an SDE with Lebesgue-Besov drift in $L^r-B^\beta_{p,q}$ for $\alpha>(1+[d/p])/(1+[d/r])$. More precisely, the martingale problem is well-posed when $\beta>(1-\alpha+[\alpha/p]+[\alpha/r])/2$, and a rigorous meaning to the formal dynamics can be given when $\beta>(1-\alpha+[2\alpha/p]+[2\alpha/r])/2$, and thus weak well-posedness is obtained. Moreover, pathwise uniqueness holds in $d=1$. For $\alpha=2$, $p,r=+\infty$, the reachable drift regularity in this work is $\beta>-1/2$. We also refer to \cite{DelarueDiel2016RoughPaths} and \cite{KrempPerkowski2022}, where for the same type of SDE, the authors reached the threshold $\beta>2(1-\alpha)/3$ by imposing additional structure on the drift (using rough paths and paracontrolled calculus techniques, respectively). For $\alpha=2$, this condition reads as $\beta>-2/3$.

\vspace{1mm}

In the framework of McKean-Vlasov SDE with a convolution-type drift and $\alpha$-stable driving noise, $\alpha\in(1,2]$, \cite{deraynal2022multidimensional} established weak well-posedness for a distributional interaction kernel in $L^r-B^{\beta}_{p,q}$, where
\begin{align}
	\label{C0}
	\beta>1-\alpha+\frac{\alpha}{r}+\frac{d}{p},
\end{align}
which allows us to reach $\beta>-1$ when $\alpha=2$, $p,r=+\infty$. Strong well-posedness holds under the slightly reinforced condition
\begin{align}
	\label{C0S}
	\beta>2-\frac{3}{2}\alpha+\frac{\alpha}{r}+\frac{d}{p}.
\end{align}

Note that the aforementioned results do not require any regularity of the initial probability measure.

\vspace{1mm}

In \cite{deraynal2023multidimensionalstabledrivenmckeanvlasov}, the authors established that additional initial regularity of $\mu$ improves \eqref{C0} and \eqref{C0S} for $\beta\in(-1,0]$. Namely, assuming that $\mu\in\mathcal{P}(\R^d)\cap B^{\beta_0}_{p_0,q_0}$ with $\beta_0+\frac{d}{p'_0}\geq0$, where $\frac{1}{p_0'}=1-\frac{1}{p_0}$, equation \eqref{main-sde} is weakly well-posed if
\begin{align}
	\label{C1}
	\beta>1-\alpha+\frac{\alpha}{r}+\Big(-\beta+\frac{d}{p}-\zeta_0\Big)_+,
\end{align}
where
\begin{align}
	\label{def-zeta0}
	\zeta_0:=\Big(\beta_0+\frac{d}{p'_0}\Big)\Big(1\wedge\frac{p'_0}{p}\Big)\geq0.
\end{align}

It is also strongly well-posed when
\begin{align}
	\label{C1S}
	\beta>\Bigg(2-\frac{3}{2}\alpha+\frac{\alpha}{r}+\Big(-\beta+\frac{d}{p}-\zeta_0\Big)\Bigg)\vee\Bigg(1-\alpha+\frac{\alpha}{r}+\Big(-\beta+\frac{d}{p}-\zeta_0\Big)_+\Bigg).
\end{align}

The positive part in \eqref{C1}, \eqref{C1S} indicates that we can consider the initial law arbitrarily smooth, but we can profit from it only up to some threshold that prevents us from reaching $\beta\leq-1$. Conditions \eqref{C1}, \eqref{C1S} are indeed an improvement w.r.t. \eqref{C0}, \eqref{C0S}, because having $\beta_0+\frac{d}{p'_0}>0$ allows us to consider a wider range of values for $p,r$ for reaching the same values of $\beta$ as in \eqref{C0}, \eqref{C0S}.

\vspace{1mm}

More importantly, exploiting the regularity of the initial measure, the threshold $\beta=-1$ has been reached in \cite{deraynal2023multidimensionalstabledrivenmckeanvlasov}. More precisely, for $\beta=-1$ the equation \eqref{main-sde} is weakly well-posed under conditions $\mu\in\mathcal{P}(\R^d)\cap B^{\beta_0}_{p_0,q_0}$ with $\beta_0+\frac{d}{p_0'}>0$ and
\begin{align}
	\label{C2}
	0>1-\alpha+\frac{\alpha}{r}+\Big(1+\frac{d}{p}-\zeta_0\Big)_+,\quad \Div(b)\in L^r(B^{\beta}_{p,q}).
\end{align}

It is additionally strongly well-posed if
\begin{align}
	\label{C2S}
	\beta=-1>\Bigg(2-\frac{3}{2}\alpha+\frac{\alpha}{r}+\Big(\frac{d}{p}-\zeta_0\Big)\Bigg)\vee\Bigg(-\alpha+\frac{\alpha}{r}+\Big(1+\frac{d}{p}-\zeta_0\Big)_+\Bigg).
\end{align}

Conditions \eqref{C2}, \eqref{C2S} hold for initial data that possesses at least some positive regularity: this can be seen for $\alpha=2$, $p,r=+\infty$. We also observe an additional condition on the divergence of the drift which may seem limiting. However, it is rather natural in the related PDE models, see \cite[Section 5]{deraynal2023multidimensionalstabledrivenmckeanvlasov}.

\vspace{1mm}

The proof in \cite{deraynal2023multidimensionalstabledrivenmckeanvlasov} as well as in \cite{deraynal2022multidimensional} consists in introducing the SDE with a mollified drift for which a weak and a strong solution exists whose time marginal distributions admit a smooth density. The a priori estimates, uniform in the mollification parameter $\varepsilon$, for the mollified density established via the Grönwall-like inequality then allow us to pass to the non-linear martingale problem. One of the main tools therein for proving estimates on the mollified density is the product rule in Besov space (see \eqref{besov-prop-pr2} introduced in Proposition \ref{prop-product-rule-usual}) which is used to make appear the same norms of the mollified densities through the Duhamel formula everywhere leading to a quadratic inequality on the suitable norm of the density. Whereas in \cite{deraynal2022multidimensional}, the authors followed the so-called ``dequadrification'' approach that allows reducing the order of the polynomial and the uniform estimates are concluded via the circular argument. 

\vspace{1mm}

The goal of the current work is to exploit the theory of \emph{linear} (measure-independent) SDEs with Lebesgue-Besov drift with negative regularity index that appeared in \cite{deraynal2022-linear-besov} in order to prove well-posedness of \eqref{main-sde} (in a sense to be specified) for an interaction kernel $b$ of negative regularity $\beta<-1$. Such a result involves conditions that relate the stable constant of the driving noise $\alpha$, the integrability parameters $r,p,q$, and the regularity index $\beta$ of the kernel, the integrability parameters $p_0,q_0$, and the regularity index $\beta_0$ of the initial measure, and the dimension $d$. These parameters are related as well through the additional parameter $\theta\in[0,\frac{1}{2})$ used to denote the regularity index of the Besov norm of the convolution drift in \eqref{main-sde}. More precisely, we derive a set of conditions such that if for given $\alpha,r,p,q,\beta,p_0,q_0,\beta_0$ in proper ranges there exists a feasible $\theta$ satisfying these conditions, then we have either well-posedness of the corresponding martingale problem or as well weak well-posedness of the McKean-Vlasov SDE \eqref{main-sde}.

\vspace{1mm}

The overall proof strategy here is the one of \cite{deraynal2023multidimensionalstabledrivenmckeanvlasov} (see Section \ref{subsection-strategy} for details). However, the main ingredient that makes the essential difference is a suitable product rule in Besov space with negative regularity (see \eqref{ineq-product-rule}, Proposition \ref{prop-product-rule-usual}) that is used to prove non-explosive estimates on the mollified density in the appropriate space. As we will see, the product rule for Besov space with negative regularity (see \eqref{besov-prop-pr2} defined in Proposition \ref{prop-product-rule-usual}) used in \cite{deraynal2023multidimensionalstabledrivenmckeanvlasov} is not enough to obtain the claimed result.

\paragraph{Organization of the paper} In Section \ref{subsection-main-results}, we introduce the appropriate notions of solutions for our framework and we state our main results. The strategy of the proof is briefly described in Section \ref{subsection-strategy}. In Section \ref{section-tools}, we state the utilities we need to address our problem. This includes the definition and properties of Besov spaces, as well as various product rules. Section \ref{section-pde} is devoted to preliminary results about the related non-linear Fokker-Planck equation. We prove well-posedness of the martingale problem, existence and uniqueness of a weak solution to \eqref{main-sde} in short time in Section \ref{section-sde}. In Section \ref{section-sde-lt}, we extend this result to an arbitrary finite time horizon under a suitable ``smallness'' condition on the initial distribution. Section \ref{section-sde-lplq} is devoted to proving a particular case of well-posedness of \eqref{main-sde} when it can be understood in a classical sense. This particular case, however, requires as well additional regularity of the initial data.

\subsection{Main results and examples}
\label{subsection-main-results}

\textbf{Assumption on the driving noise.}
Let $\nu$ be the Lévy measure of $\mathcal{W}$ when $\alpha\in(1,2)$. It is known that it admits a decomposition given by $\nu(dz)=w(d\zeta)/\rho^{1+\alpha}1_{\rho>0}$, where $w$ is a measure on the $d$-dimensional sphere $\mathbb{S}^{d-1}$, see e.g. \cite{Sato1999LvyPA}. We assume that the symmetric spherical measure $w$ satisfies the \emph{uniform non-degeneracy} condition
\begin{align*}
	\kappa^{-1}|\lambda|^\alpha\leq\int_{\mathbb{S}^{d-1}}|\zeta\cdot\lambda|^\alpha w(d\zeta)\leq\kappa|\lambda|^\alpha,
\end{align*}
for all $\lambda\in\R^d$ and for some $\kappa\geq1$. This condition in particular guarantees that the $\alpha$-stable process $\mathcal{W}$ admits a smooth density $p^\alpha_t$ for any $t>0$, see \cite{Kolokoltsov2000SymmetricSL} and \cite{WatanabeSingularMeasures} for more information about singular spectral measures.

\vspace{2mm}

\textbf{Notion of solution.}
We carefully introduce the notion of a solution to the martingale problem and a weak solution used in the sequel. We denote the canonical space
\begin{align*}
	\Omega_\alpha:=
	\begin{cases}
		\displaystyle\mathcal{C}([0,S],\R^d),\quad \alpha=2, \\
		\displaystyle\mathcal{D}([0,S],\R^d),\quad \alpha\in(1,2),
	\end{cases}
\end{align*}
where $\mathcal{D}$ denotes a Skorokhod space of càdlàg functions and $S\in(0,T]$.

\begin{definition}[Martingale problem]
\label{def-martingale-problem}
We say that a measure $\Prob^\alpha\in\mathcal{P}(\Omega_\alpha)$ where $\Omega_\alpha$ is equipped with its canonical filtration, is a \emph{martingale solution} to the \emph{non-linear martingale problem} associated with the McKean-Vlasov equation \eqref{main-sde} if for the canonical process $(X_s)_{s\in[0,S]}$ of $\Prob^\alpha$,
\begin{enumerate}
		\item $\Prob^\alpha\circ X(0)^{-1}=\mu$,
		\item For a.a. $s\in(0,S]$, $\Prob^\alpha\circ X(s)^{-1}$ is absolutely continuous w.r.t. Lebesgue measure and for some $\gamma>0$ its density satisfies 
		\begin{align*}
			\esssup_{s\in(0,S]}s^{\gamma}\Big\|\frac{d\Prob^\alpha\circ X(s)^{-1}}{d\text{Leb}}\Big\|_{B^{-\beta-\theta}_{p',1}}<+\infty.
		\end{align*}
		\item For all $f\in\mathcal{C}([0,S],\mathcal{S}(\R^d,\R))$, $g\in\mathcal{C}^1(\R^d,\R)$ s.t. $Dg\in B^{-1-\theta+\alpha-\frac{\alpha}{r_\theta}}_{\infty,\infty}$,
		      \begin{align*}
			      \Big(u(s,X_s)-u(0,X_0)-\int_0^sf(r,X_r)dr\Big)_{0\leq s\leq S}
		      \end{align*}
		      is a (square integrable if $\alpha=2$) $\Prob^\alpha$-martingale, where $u\in\mathcal{C}^{0}([0,S],\mathcal{C}^1(\R^d))$ with $Du\in\mathcal{C}^0_b([0,S],B^{-\theta-\varepsilon}_{\infty,\infty})$, $\varepsilon>0$, satisfies
		      \begin{align}
			      \label{mart-pr-duhamel-for-u}
				  \begin{split}
			      u(s,x)=p^\alpha_{S-s}\ast g(x)-\int_s^S\Big(\big(f-(b\ast(\Prob\circ X(v)^{-1}))\cdot Du\big)\ast p^\alpha_{v-s}\Big)(v,x)dv, \\
				  (s,x)\in[0,S]\times\R^d,
				  \end{split}
		      \end{align}
			  where $(b\ast(\Prob\circ X(v)^{-1}))(v,x)=\int_{\R^d}b(v,x-z)(\Prob\circ X(v)^{-1})(z)dz$. Here, parameters $\theta,\alpha,\beta,r_\theta,d,p$ satisfy some constraints that are presented later on and are omitted here for simplicity.
	\end{enumerate}
\end{definition}

\begin{remark}
	If there exists $u\in\mathcal{C}^{0}([0,S],\mathcal{C}^1(\R^d))$ with $Du\in\mathcal{C}^0_b([0,S],B^{-\theta-\varepsilon}_{\infty,\infty})$, $\varepsilon>0$, satisfying the Duhamel representation \eqref{mart-pr-duhamel-for-u}, then we say that $u:[0,S)\times\R^d\to\R$ is a \emph{mild} solution to the formal Cauchy problem
	\begin{align*}
		(\partial_s+(b\ast(\Prob\circ X(v)^{-1}))\cdot D+\mathcal{L}^\alpha)u(s,x)=f(s,x)\text{ on }[0,S)\times\R^d,\quad u(S,\cdot)=g(\cdot)\text{ on }\R^d.
	\end{align*}
	See Section 1.3 of \cite{deraynal2022-linear-besov} for more details.
\end{remark}

Since we are dealing with a distributional drift, the existence of a weak solution to \eqref{main-sde} is not straightforward from the existence of a martingale solution in the sense of Definition \ref{def-martingale-problem}. For reconstructing the dynamics of \eqref{main-sde}, we use results obtained in \cite{deraynal2022-linear-besov} under \emph{good relation for dynamics} \eqref{good-relation-enhaced}. Therein, the authors define the integrated drift of the SDE as a stochastic non-linear Young integral in the following sense.

\begin{definition}[Stochastic non-linear Young integral]
	\label{def-young-integral}
	Let $\tau>0$, $(\Omega,\mathcal{F},(\mathcal{F}_t)_{0\leq t\leq\tau},\Prob)$ be a filtered probability space and let $(\psi_t)_{0\leq t\leq\tau}$ be a progressively measurable process. Let $(A(s,t))_{0\leq s\leq t\leq\tau}$ be a continuous and progressively measurable map, i.e., for any $0\leq s\leq t$,
	\begin{align*}
		\Omega\times\{s'\in[0,s],\ t'\in[0,t],\ s'\leq t'\}\ni(\omega,s',t')\mapsto A(s',t')
	\end{align*}
	is $\mathcal{F}_t\otimes\mathcal{B}(\{s'\in[0,s],\ t'\in[0,t],\ s'\leq t'\})$-measurable and
	\begin{align*}
		\{s'\in[0,\tau],\ t'\in[0,\tau],\ s'\leq t'\}\ni(s,t)\mapsto A(s,t)
	\end{align*}
	is continuous. For $l\geq1$, we say that
	\begin{align*}
		\int_0^\tau\psi_tA(t,t+dt):=\lim_{\substack{\Delta\text{ partition of } [0,\tau], \\
				|\Delta|\to0}}\sum_{t_i\in\Delta}\psi_{t_i}A(t_i,t_{i+1})\quad \text{ in } L^l(\Omega,\Prob)
	\end{align*}
	is a \emph{$L^l$-stochastic Young integral} of $\psi$ w.r.t. the pseudo increment $A$ when it exists.
\end{definition}

With this definition in hand, we can properly define a weak solution to the non-linear McKean-Vlasov SDE \eqref{main-sde}.

\begin{definition}[Weak solution]
	\label{def-weak-solution}
	We say that a pair $(X,\mathcal{W})$ of adapted processes on \\$(\Omega,\mathcal{F},(\mathcal{F}_t)_{t\geq0},\Prob)$ is a \emph{weak solution} to the formal McKean-Vlasov equation \eqref{main-sde} on the time interval $[0,S]$ with a random initial value $\xi$ such that $\text{Law}(\xi)=\mu\in\mathcal{P}(\R^d)$, if
	\begin{itemize}
		\item $\mathcal{W}$ is an $(\mathcal{F}_t)_{t\geq0}$ $\alpha$-stable process
		\item For any $s\in[0,S]$,
		      \begin{align*}
			      X_s=\xi+\int_0^s\mathcal{G}(v,X_v,dv)+\mathcal{W}_s, \quad \Prob-\text{a.s.},
		      \end{align*}
		      where the integral is understood as a $L^1$-stochastic non-linear Young integral in the sense of Def. \ref{def-young-integral} with $\psi\equiv1$, $A(v,v+dv)=\mathcal{G}(v,X_v,dv)$, $A(v_i,v_{i+1})=\mathcal{G}(v_i,X_{v_i},v_{i+1}-v_i)$ therein, and $\mathcal{G}$ satisfies the properties below,
		\item For any $s\in[0,S]$,
		      \begin{align*}
			      \Exp[|\int_0^s\mathcal{G}(v,X_v,dv)|]<+\infty,
		      \end{align*}
		\item For any $0\leq r\leq s\leq S$, $x\in\R^d$,
		      \begin{align*}
			      \mathcal{G}(r,x,s-r)=\int_r^s\int_{\R^d}(b\ast\bmu_v^{\mu})(v,y)p^\alpha(v-r,y-x)dydv,
		      \end{align*}
		      where $p^\alpha$ is the density of $\mathcal{W}$.
	\end{itemize}

	Moreover, we say that \emph{weak uniqueness} holds for \eqref{main-sde} if for any two weak solutions $(X,\mathcal{W})$ on $(\Omega,\mathcal{F},(\mathcal{F})_{t\geq0},\Prob)$ and $(\bar{X},\bar{\mathcal{W}})$ on $(\bar{\Omega},\bar{\mathcal{F}},(\bar{\mathcal{F}})_{t\geq0},\bar{\Prob})$ with $X_0=\bar{X}_0\sim\mu\in\mathcal{P}(\R^d)$, it holds $\text{Law}((X_t)_{t\geq0})=\text{Law}((\bar{X}_t)_{t\geq0})$.
\end{definition}

\begin{definition}[Strong solution]
	We say that there exists a \emph{strong solution} to \eqref{main-sde} if for a given noise $(\mathcal{W}_t)_{t\geq0}$ on $(\Omega,\mathcal{F},(\mathcal{F}_t)_{t\geq0},\Prob)$ where $(\mathcal{F}_t)_{t\geq0}$ is the complete filtration generated by $(\mathcal{W}_t)_{t\geq0}$, there exists a càdlàg process $(X_t)_{t\geq0}$ that satisfies equation \eqref{main-sde} a.s. and is adapted to $(\mathcal{F}_t)_{t\geq0}$. We say that a strong solution is \emph{unique} if for $(X_t)_{t\geq0}$, $(\bar{X}_t)_{t\geq0}$ two strong solutions to \eqref{main-sde} such that $X_0=\bar{X}_0\sim\mu\in\mathcal{P}(\R^d)$, then $\Prob(X_t=\bar{X}_t,\ \forall t\geq0)=1$.
\end{definition}

\vspace{2mm}

\textbf{Main results.}
We now state the main results of this paper. For readability, we consider short- and long-time regimes separately. We also recall that the parameter $-\theta\leq0$ is used to denote the spatial regularity of the convolution drift: in the following theorems, it serves as an additional parameter appearing in the constraints binding the input data of the SDE \eqref{main-sde}, and at the same time, it has to satisfy some specific constraints.
\begin{theorem}[Short-time well-posedness]
	\label{thm-main-sde}
	Let
	\begin{align}
		\label{mu-range}
		\mu\in\mathcal{P}(\R^d)\cap B^{\beta_0}_{p_0,q_0}\cap B^{\bar{\theta}}_{1,\infty},
	\end{align}
	where
	\begin{align*}
		\beta_0\geq0,\ p_0,q_0\in[1,+\infty]
	\end{align*}
	are such that $\zeta_0>0$, where $\zeta_0$ is given by \eqref{def-zeta0}, and $\bar{\theta}>0$ is such that
	\begin{align*}
		\bar{\theta}>\theta,
	\end{align*}
	for some $\theta\in[0,\frac{1}{2})$. Let $\alpha\in(1,2]$, $\beta\in(-2,-1)$, $p,r\in(1,+\infty]$, $q\in[1,+\infty]$. Finally, let
	\begin{align*}
		b\in L^r((0,T],B^\beta_{p,q}),\quad \Div(b)\in L^r((0,T],B^\beta_{p,q}),
	\end{align*}
	where
	\begin{equation}
		\label{C3}
		\begin{gathered}
			\beta>-\alpha+\frac{\alpha}{r}+\Big(-\theta-\beta+\frac{d}{p}-\zeta_0\Big)_+\quad \text{ and } \quad \beta<-1-2\theta,\\
			\beta\neq-\theta-\bar{\beta}_0+\frac{d}{p},
		\end{gathered}
	\end{equation}
	where $\bar{\beta}_0=\beta_0(1\wedge\frac{p_0'}{p})$.

	\begin{enumerate}[label=(\arabic*)]
		\item{\emph{\textbf{(Martingale solution)}}} If the parameters $(\theta,r)$ satisfy conditions
		      \begin{equation}
			      \label{good-relation}
			      \begin{gathered}
				      r>\frac{\alpha}{\alpha-1},\quad \theta\in\Bigg[0,\frac{1}{2}\Big(\alpha-\frac{\alpha}{r}-1\Big)\Bigg),
			      \end{gathered}
		      \end{equation}
		      then there exists time horizon $\mathcal{T}\in(0,T]$ small enough such that the martingale problem associated with \eqref{main-sde} is well-posed up to any time $S\in(0,\mathcal{T})$ in the sense of Definition \ref{def-martingale-problem}. Moreover, for $r_\theta\in(\alpha,r)$,
			  \begin{align*}
				b\ast\bmu^{\mu}\in L^{r_\theta}(B^{-\theta}_{\infty,\infty}).
			  \end{align*}

		\item{\emph{\textbf{(Weak WP)}}}
			If additionally to \eqref{good-relation} the parameter $r$ satisfies condition
				      \begin{align}
					      \label{good-relation-enhaced}
					      \begin{split}
						      r>\frac{\alpha}{\alpha-1}\vee2\alpha,
					      \end{split}
				      \end{align}
				      then the formal equation \eqref{main-sde} admits a unique weak solution in the sense of Definition \ref{def-weak-solution} up to any time $S\in(0,\mathcal{T})$.
	  
		\item{\emph{\textbf{(Regularity of the law)}}} If either \eqref{good-relation} or \eqref{good-relation}-\eqref{good-relation-enhaced} holds, then the time marginals $(\bmu_s^{\mu})_{s\in[0,S]}$ of the solution for $S\in(0,\mathcal{T})$ have a density $\brho_{\mu}(s,\cdot)$ for a.a. time $s\in(0,S]$ and
		      \begin{align*}
			      \sup_{s\in(0,S]} s ^\gamma\|\brho_{\mu}(s,\cdot)\|_{B^{-\beta-\theta}_{p',1}}<+\infty,
		      \end{align*}
		      where
		      \begin{align*}
			      \gamma>\frac{1}{\alpha}\Big(-\theta-\beta+\frac{d}{p}-\zeta_0\Big)_+.
		      \end{align*}
	\end{enumerate}
\end{theorem}

\begin{theorem}[Long-time well-posedness]
	\label{thm-main-sde-lt}
	Let $\mu$ and $b$ be as before for $\alpha\in(1,2]$, $\beta\in(-2,-1)$, $r\in(1,+\infty]$, $p\in[1,+\infty)$, $q\in[1,+\infty]$, $\theta\in[0,\frac{1}{2})$. Assume that condition \eqref{C3} of the previous theorem holds. Assume additionally that for $\varepsilon\ll1$,
	\begin{align}
		\label{C3LT}
		1<\alpha\Big(1-\frac{1}{r}\Big)\leq\min\Big(1+\frac{d}{p},2-\varepsilon\Big).
	\end{align}
	\begin{enumerate}[label=(\arabic*)]
		\item{\emph{\textbf{(Martingale solution)}}} If the parameters $(\theta,r)$ satisfy conditions \eqref{good-relation}, then there exist a constant $\mathcal{C}_0>0$ and time $\mathcal{T}_1\in(0,T]$ defined by
		      \begin{align*}
			      \mathcal{T}_1:=\mathcal{T}\ind_{\{\|\mu\|_{B^{\bar{\beta}_0}_{\bar{p}_0,\bar{q}_0}}\geq \mathcal{C}_0\}}+T\ind_{\{\|\mu\|_{B^{\bar{\beta}_0}_{\bar{p}_0,\bar{q}_0}}<\mathcal{C}_0\}},
		      \end{align*}
		      such that the martingale problem associated with \eqref{main-sde} is well-posed up to any time $S\in(0,\mathcal{T}_1)$ in the sense of Definition \ref{def-martingale-problem}. Moreover, for $r_\theta\in(\alpha,r)$,
			  \begin{align*}
				b\ast\bmu^{\mu}\in L^{r_\theta}(B^{-\theta}_{\infty,\infty}).
			  \end{align*}

		\item{\emph{\textbf{(Weak WP)}}} If the parameter $r$ additionally satisfies condition \eqref{good-relation-enhaced}, then the formal equation \eqref{main-sde} admits a unique weak solution in the sense of Definition \ref{def-weak-solution} up to any time $S\in(0,\mathcal{T}_1)$.
		\item{\emph{\textbf{(Regularity of the law)}}} If either \eqref{good-relation} or \eqref{good-relation-enhaced} holds, then the time marginals $(\bmu_s^{\mu})_{s\in[0,S]}$ of the solution for $S\in(0,\mathcal{T}_1)$ have a density $\brho_{\mu}(s,\cdot)$ for a.a. time $s\in(0,S]$ and
		      \begin{align*}
			      \sup_{s\in(0,S]}\big( s \wedge1\big)^{\gamma_1}\big( s \vee1\big)^{\gamma_2}\|\brho_{\mu}(s,\cdot)\|_{B^{-\beta-\theta}_{p',1}}<+\infty,
		      \end{align*}
		      where
		      \begin{align*}
			      \gamma_1 & >\frac{1}{\alpha}\Big(-\theta-\beta+\frac{d}{p}-\zeta_0\Big)_+,\quad\gamma_2:=\frac{1}{r'}-\frac{1}{\alpha}.
		      \end{align*}
	\end{enumerate}
\end{theorem}

We will see that the upper bound on $\beta$ in \eqref{C3} comes from the product rule \eqref{ineq-product-rule} defined in Proposition \ref{prop-product-rule-usual}. Since we are interested in lowering the negative threshold on $\beta$, the upper bound does not appear to be too restrictive. 

\vspace{1mm}

We also emphasize that the additional regularity assumption on $\mu\in B^{\bar{\theta}}_{1,\infty}$ for $\bar{\theta}>\theta$ is needed to prove that the solution of \eqref{main-sde} admits a density, see Remark \ref{remark-mu-reg}.

\begin{remark}[About a choice of $\theta$]
	Let us discuss the condition on the parameter $\theta$ appearing in \eqref{good-relation}.
	First thing to note is that the theory of \emph{linear} (measure independent) SDE with a distributional drift $\mathcal{B}\in L^{r-\varepsilon}(B^{-\theta}_{\infty,\infty})$, \cite{deraynal2022-linear-besov}, tells us that for $r>\alpha$ the martingale problem associated with
	\begin{align*}
		dX_s=\mathcal{B}(s,X_s)ds+d\mathcal{W}_s,\quad X_t=\xi\sim\mu\in\mathcal{P}(\R^d),
	\end{align*}
	is well-posed in the sense of Definition \ref{def-martingale-problem} if
	\begin{align}
		\label{cond-theta-general}
		\theta\in\Big(0,\frac{1}{2}\big(\alpha-1\big)\Big).
	\end{align}
	If moreover $r>2\alpha$ (as in \eqref{good-relation-enhaced}), then there exists a unique weak solution in the sense of Definition \ref{def-weak-solution} to the SDE. In particular,
	\begin{align*}
		\alpha=2\implies\theta<\frac{1}{2}.
	\end{align*}
	We refer to Section \ref{subsection-strategy} for a more detailed discussion. Also from \eqref{C3} we see that,
	\begin{align*}
		-\alpha+\frac{\alpha}{r}\leq-\alpha+\frac{\alpha}{r}+\Big(-\theta-\beta+\frac{d}{p}-\zeta_0\Big)_+<\beta<-1-2\theta\implies\theta<\frac{1}{2}\Big(\alpha-\frac{\alpha}{r}-1\Big).
	\end{align*}
	However, in order for such $\theta$ to exist so that \eqref{cond-theta-general} is also satisfied, we need an additional condition on $r$. More precisely, it should hold that
	\begin{align*}
		r>\frac{\alpha}{\alpha-1},
	\end{align*} 
	so that $\alpha-\frac{\alpha}{r}-1>0$. Additionally, thanks to the embeddings between Lebesgue and Besov space of regularity zero (see Section \ref{section-tools}, \eqref{besov-prop-e1}), we also include $\theta=0$, which allows us to handle a drift that is not a distribution but a function with zero regularity. Therefore, this explains the possible values for $\theta$ in \eqref{good-relation}.
\end{remark}

The following theorem is a particular case of Theorem \ref{thm-main-sde} when $\theta=0$. In this case, the McKean-Vlasov equation \eqref{main-sde} is understood in a usual sense, and the convolution drift is a well-defined function. This implies that if a weak solution exists, it is understood in a classical sense.
\begin{theorem}[Classical well-posedness in short time]
	\label{thm-main-sde-lplq}
	Let $\mu\in\mathcal{P}(\R^d)\cap B^{\beta_0}_{p_0,q_0}$ and \\ $b\in L^r((0,T],B^{\beta}_{p,q})$ for $\alpha\in(1,2]$, $\beta\in(-2,-1)$, $r\in(1,+\infty]$, $p\in[1,+\infty)$, $q\in[1,+\infty]$.
	\begin{enumerate}[label=(\arabic*)]
		\item{\emph{\textbf{(Weak WP)}}} If
		      \begin{align}
			      \label{C4}
			      \beta>-\alpha+\frac{\alpha}{r}+\Big(-\beta+\frac{d}{p}-\zeta_0\Big)_+, \quad\Div(b)\in L^r((0,T],B^{\beta}_{p,q}),
		      \end{align}
		      then there exists a small time horizon $0<\mathcal{T}_0\leq T$ such that for any $S\leq\mathcal{T}_0$, the McKean-Vlasov SDE \eqref{main-sde} admits a unique classical weak solution. Moreover, for $r_0\in\big(\frac{\alpha}{\alpha-1},\frac{r}{1+r\gamma^*}\big)$,
			  \begin{align*}
				b\ast\bmu^{\mu}\in L^{r_0}(L^\infty),
			  \end{align*}
			  where $\gamma^*>0$ is defined by
		      \begin{align*}
			      \gamma^*=\frac{1}{\alpha}\Big(-\beta+\frac{d}{p}-\zeta_0+\Big(\frac{1+\eta}{2\eta}\Big)\Gamma\Big)>0,
		      \end{align*}
			  for $\Gamma\in(0,1)$ given by
		      \begin{align*}
			      \Gamma=\eta\Big(\alpha-1+\beta-\frac{\alpha}{r}-\frac{d}{p}+\zeta_0\Big),\quad \eta\in(0,1).
		      \end{align*}
		\item{\emph{\textbf{(Strong WP)}}} If under \eqref{C4}
		      \begin{align}
			      \label{C4S}
			      \beta>\Bigg(1-\frac{3}{2}\alpha+\frac{\alpha}{r}+\Big(-\beta+\frac{d}{p}-\zeta_0\Big)\Bigg)\vee\Bigg(-\alpha+\frac{\alpha}{r}+\Big(-\beta+\frac{d}{p}-\zeta_0\Big)_+\Bigg),
		      \end{align}
		      then the unique solution is strong.
		\item{\emph{\textbf{(Regularity of the law)}}} If \eqref{C4} holds, then the time marginals $(\bmu_s^{\mu})_{s\in[0,S]}$ of the solution for $S\in(0,\mathcal{T}_1)$ have a density $\brho_{\mu}(s,\cdot)$ for a.a. time $s\in(0,S]$ and
		      \begin{align*}
			      \sup_{s\in(0,S]} s ^{\gamma^*}\|\brho_{\mu}(s,\cdot)\|_{B^{-\beta+\Gamma}_{p',1}}<+\infty.
		      \end{align*}
	\end{enumerate}
\end{theorem}

\begin{example}[Thresholds in short time]
	\label{example-1}
	Theorem \ref{thm-main-sde} involves many parameters, and it may appear intricate to balance them. Let us apply this theorem to the simplest set of parameters when the driving noise is Brownian motion. Namely, we consider for $d\geq1$,
	\begin{align*}
		\alpha=2, \ p=+\infty, \ r=+\infty, \ q\in[1,+\infty].
	\end{align*}
	For these parameters, we can have at most $\beta>-2$ from \eqref{C3}. Then recalling that $\zeta_0:=\Big(\beta_0+\frac{d}{p'_0}\Big)\Big(1\wedge\frac{p'_0}{p}\Big)$, using convention $\frac{p_0'}{p}=1$ when $p_0',p=+\infty$ and taking for $\varepsilon>0$
	\begin{align*}
		\beta=-2+\varepsilon,\quad p_0=1,\ q_0\in[1,+\infty],\ \zeta_0=\beta_0>\frac{3}{2}\Big(1-\varepsilon\Big),
	\end{align*}
	there exists $\theta=\frac{1}{2}-\varepsilon'$, $\varepsilon<2\varepsilon'$, such that \eqref{C3} and \eqref{good-relation} are verified. The initial law $\mu$ should also be at least in $B^{\bar{\theta}}_{1,\infty}$. This is satisfied for the current choice of the parameters $\beta_0,p_0$ thanks to the embedding of Besov spaces (see \eqref{besov-prop-e2}).

	\vspace{1mm}

	Hence, we obtain the following result: for $\alpha=2$, the McKean-Vlasov SDE \eqref{main-sde} with a drift $b$ and initial distribution $\mu$ such that
	\begin{align*}
		b,\Div(b)\in L^\infty((0,T],B^{-2+}_{\infty,q})\text{ and }\mu\in\mathcal{P}(\R^d)\cap B^{3/2-}_{1,q_0},
	\end{align*}
	admits a unique weak solution on some small time interval in the sense of Definition \ref{def-weak-solution}.
\end{example}

\begin{remark}[Comparison of \eqref{C3} and \eqref{C4}]
	\label{remark-compare-c3-c2}
	We note that the lower bound on $\beta$ in \eqref{C3} is equivalent to the lower bound on $\beta$ in \eqref{C4} when taking $\theta=0$. However, when $\theta\in(0,\frac{1}{2})$, it is the parameter that allows us to demand less regularity on the initial measure $\mu$ in order to achieve the same regularity on $b$ that guarantees the well-posedness of \eqref{main-sde}.

	\vspace{1mm}

	Let us make a comparison to Example \ref{example-1}. Consider the following parameters
	\begin{align*}
		\alpha=2, \ p=+\infty, \ r=+\infty, \ q\in[1,+\infty].
	\end{align*}
	Then for $\varepsilon>0$ taking
	\begin{align*}
		\beta=-2+\varepsilon,\quad p_0=1,\ q_0\in[1,+\infty],\ \zeta_0=\beta_0>2(1-\varepsilon),
	\end{align*}
	we see that Theorem \ref{thm-main-sde-lplq} for the weak well-posedness is satisfied. Note that the gain of Example \ref{example-1} compared to the current one is the smoothness of the initial measure $\mu$ given exactly by $\theta=\frac{1}{2}-\varepsilon'$. This means that we can profit from this gain by paying the price of working in a distributional framework of the SDE \eqref{main-sde} with an adapted notion of a solution developed in \cite{deraynal2022-linear-besov}. We also note that both in Theorem \ref{thm-main-sde} and Theorem \ref{thm-main-sde-lplq}, some regularity on the initial data $\mu$ is needed. This is necessary to handle kernels of regularity $\beta\leq-1$.

\end{remark}

\begin{example}[Thresholds in long time]
	Before looking into an example satisfying Theorem \ref{thm-main-sde-lt}, let us discuss condition \eqref{C3LT} in the long time and how it changes the threshold on $\beta$. Condition \eqref{C3LT} is the main difference between short and long time regimes. In short time, the only condition we have is $\alpha(1-\frac{1}{r})>1$ (see \eqref{good-relation}), whereas in long time condition $\alpha(1-\frac{1}{r})\leq\min(1+\frac{d}{p},2-\varepsilon)$ in \eqref{C3LT} is more restrictive. For example, we cannot consider cases $\alpha=2,r=+\infty$ (as in Example \ref{example-1}). Moreover, it becomes even more restrictive when $p$ is much greater than $d$. From the lower bound $\alpha(1-\frac{1}{r})>1$, we see also that $p$ cannot be taken to infinity.

	\vspace{1mm}

	Now, let us construct an example. For $\varepsilon>0$, let
	\begin{align*}
		\alpha=2-\varepsilon,\ p=d,\ r=+\infty,\ q\in[1,+\infty]
	\end{align*}
	so that
	\begin{align*}
		\alpha\Big(1-\frac{1}{r}\Big)=2-\varepsilon,
	\end{align*}
	i.e., \eqref{C3LT} is satisfied. Let
	\begin{align*}
		\beta=-2+2\varepsilon,\quad p_0=1,\ q_0\in[1,+\infty],\ \zeta_0=\beta_0=\frac{5}{2}-\varepsilon.
	\end{align*}
	Then there exists $\theta=\frac{1}{2}-\varepsilon'$ such that \eqref{C3}, \eqref{C3LT} and \eqref{good-relation}-\eqref{good-relation-enhaced} are verified and Theorem \ref{thm-main-sde-lt} holds true if $\|\mu\|_{B^{5/2-}_{1,q_0}}<\mathcal{C}_0$ for some $\mathcal{C}_0$ small, i.e., for $\alpha=2-\varepsilon$, the McKean-Vlasov SDE \eqref{main-sde} with a drift $b$ and initial distribution $\mu$ such that
	\begin{align*}
		b,\Div(b)\in L^\infty((0,T],B^{-2+}_{d,q})\text{ and }\mu\in\mathcal{P}(\R^d)\cap B^{5/2-}_{1,q_0},\ \|\mu\|_{B^{5/2-}_{1,q_0}}<\mathcal{C}_0
	\end{align*}
	admits a unique weak solution on any finite time interval in the sense of Definition \ref{def-weak-solution}.

	\vspace{1mm}

	Comparing this result to Example \ref{example-1}, we see that in the long time, well-posedness of \eqref{main-sde} depends more on the relation between $p$ and $d$, and on the regularity of the initial measure.
\end{example}

\begin{remark}[About $d=1$ and example of a kernel in $d=2$]$ \newline$
	\textbullet\ Let us carefully note that if the coefficients of \eqref{main-sde} satisfy the conditions for weak well-posedness of Theorem \ref{thm-main-sde} (Theorem \ref{thm-main-sde-lt} in long time, respectively), then strong well-posedness in dimension one is basically automatic. More precisely, it can be deduced from a generalized Yamada-Watanabe theorem, \cite{KurtzWeakStrongSolutions}, which gives uniqueness of a strong solution once there exists a weak solution that is pathwise unique (see also \cite{ABM2020StableSDE}). In this article, we solve the non-linear Fokker-Planck PDE in some appropriate sense, so that the stochastic equation \eqref{main-sde} becomes linear, and the existing theory \cite{deraynal2022-linear-besov} on well-posedness of linear distributional SDEs is applied. Therein, pathwise uniqueness is obtained only in dimension one, and the same dimensional restriction holds for strong well-posedness. However, this consequence of Theorem \ref{thm-main-sde} does not provide an improvement in small time of the existing result \cite{deraynal2023multidimensionalstabledrivenmckeanvlasov} when $d=1$ for strong well-posedness. 
		
	\vspace{1mm}

	Indeed, this is a consequence of several things.  Firstly, our results are obtained for the kernel of regularity $\beta$ strictly in $(-2,-1)$. The upper bound comes from the second condition in \eqref{C3} and is closely linked to the Besov product rule which is used in the proof of Lemma \ref{lemma-density-estimate} and cannot be dropped. Secondly, in dimension one, divergence becomes a usual derivative. Thus, if we have a kernel $b$ so that $b,b'\in L^r(B^{-2+\varepsilon}_{p,q})$, i.e., $\beta=-2+\varepsilon$, $\varepsilon\in(0,1)$, then it implies that $b\in L^r(B^{-1+\varepsilon}_{p,q})$. Working with the weaker regularity $\beta=-2+\varepsilon$, we need to verify condition \eqref{C3}. In this case,
	\begin{align*}
		\theta<\frac{1}{2}-\frac{\varepsilon}{2},
	\end{align*}
	and the other parameters should satisfy
	\begin{align*}
		-2+\varepsilon&>-\alpha+\frac{\alpha}{r}+\Big(-\theta+(2-\varepsilon)+\frac{1}{p}-\zeta_0\Big)_+.
	\end{align*}
	We find that the inequality is true if
	\begin{align}
		\label{example-cond-zeta-c3}
		\zeta_0>\Big(4-\alpha+\frac{\alpha}{r}+\frac{1}{p}-2\varepsilon\Big)-\theta.
	\end{align}

	On the other hand, taking directly $\beta=-1+\varepsilon$ and verifying condition \eqref{C1S} for strong well-posedness, which arises from \cite[Theorem 1]{deraynal2023multidimensionalstabledrivenmckeanvlasov} and does not require any information about the divergence of the kernel, we see that
	\begin{align*}
		-1+\varepsilon&>\Bigg(1-\alpha+\frac{\alpha}{r}+\Big((1-\varepsilon)+\frac{1}{p}-\zeta_0\Big)_+\Bigg)\vee\Bigg(2-\frac{3}{2}\alpha+\frac{\alpha}{r}+\Big((1-\varepsilon)+\frac{1}{p}-\zeta_0\Big)\Bigg).
	\end{align*}
	This imposes that
	\begin{align}
		\label{example-cond-zeta-c2}
		\begin{split}
			\zeta_0&>\Bigg(\Big(4-\alpha+\frac{\alpha}{r}+\frac{1}{p}-2\varepsilon\Big)-1\Bigg)\vee\Bigg(\Big(4-\alpha+\frac{\alpha}{r}+\frac{1}{p}-2\varepsilon\Big)-\frac{\alpha}{2}\Bigg)\\
			&=\Big(4-\alpha+\frac{\alpha}{r}+\frac{1}{p}-2\varepsilon\Big)-\frac{\alpha}{2}
		\end{split}
	\end{align}
	for $\alpha\leq2$. Let us compare the r.h.s. of \eqref{example-cond-zeta-c3} and \eqref{example-cond-zeta-c2}. By the second condition in \eqref{C3},
	\begin{align*}
		\theta<\frac{1}{2}\Bigg(\alpha-\Big(\frac{\alpha}{r}+1\Big)\Bigg)\implies\theta<\frac{\alpha}{2}.
	\end{align*}

	Consequently, $-\theta>-\frac{\alpha}{2}$ which means that the r.h.s. of \eqref{example-cond-zeta-c3} is always greater than the r.h.s. of \eqref{example-cond-zeta-c2}. In other words, applying Theorem \ref{thm-main-sde} for strong well-posedness in $d=1$ to the kernel living in $L^r(B^{-1+\varepsilon}_{p,q})$ requires more regularity on the initial data encoded by the parameter $\zeta_0$ than applying Theorem 1 of \cite{deraynal2023multidimensionalstabledrivenmckeanvlasov} under condition $(\textbf{C1\textsubscript{S}})$ therein (condition \eqref{C1S} herein). This suggests that when $d=1$, Theorem \ref{thm-main-sde} does not improve existing strong well-posedness results in short time. Similar analysis can be performed with regard to Theorem \ref{thm-main-sde-lplq}. Precisely, with the parameters like above, condition \eqref{C4} reads as
	\begin{align*}
		-2+\varepsilon>-\alpha+\frac{\alpha}{r}+\Big((2-\varepsilon)+\frac{1}{p}-\zeta_0\Big)_+,
	\end{align*}
	meaning that $\zeta_0$ should satisfy
	\begin{align*}
		\zeta_0>\Big(4-\alpha+\frac{\alpha}{r}+\frac{1}{p}-2\varepsilon\Big).
	\end{align*}

	This condition is stronger than \eqref{example-cond-zeta-c2} but weaker than \eqref{example-cond-zeta-c3}. Similarly, Theorem \ref{thm-main-sde} does not improve the weak well-posedness result in $d=1$ from \cite{deraynal2023multidimensionalstabledrivenmckeanvlasov} for the same reason. Indeed, by $(\textbf{C1})$ therein (or \eqref{C1} herein), for $\beta=-1+\varepsilon$,
	\begin{align*}
		-1+\varepsilon&>1-\alpha+\frac{\alpha}{r}+\Big((1-\varepsilon)+\frac{1}{p}-\zeta_0\Big)_+,
	\end{align*}
	thus,
	\begin{align*}
		\zeta_0&>\Big(4-\alpha+\frac{\alpha}{r}+\frac{1}{p}-2\varepsilon\Big)-1.
	\end{align*}

	Again, the r.h.s of the above inequality is smaller than the r.h.s of \eqref{example-cond-zeta-c3}, meaning that \eqref{example-cond-zeta-c3} provides a stronger condition on $\zeta_0$. Importantly, the problem we have described so far appears only when $d=1$ because the divergence condition becomes very stringent.
	
	\vspace{1mm}
	
	We emphasize, however, that in the long time Theorem \ref{thm-main-sde-lt} is an improvement of \cite[Theorem 1]{deraynal2023multidimensionalstabledrivenmckeanvlasov} in any dimension $d\geq1$ since the result therein is proved only in short time. Our long-time result is also an improvement of \cite[Theorem 1]{deraynal2022multidimensional}, where the result is stated for any finite time horizon but without exploiting the regularity of initial data. This also means that if for $d=1$ equation \eqref{main-sde} is weakly well-posed in long time according to Theorem \ref{thm-main-sde-lt}, then its strong well-posedness can be deduced using the arguments stated at the beginning of the remark.
	
	\vspace{2mm}

	\textbullet\ Let us construct an example of a kernel in $d=2$ that satisfies the set of conditions in Theorem \ref{thm-main-sde}, but that is not suitable for the results in \cite{deraynal2022multidimensional,deraynal2023multidimensionalstabledrivenmckeanvlasov} in short time. For $\gamma\in(0,\alpha-1)$, denoting by $(-\Delta)^{\frac{\gamma}{2}}$ a fractional derivative and by $\delta_x$ the Dirac distribution, consider the following $2\times2$ matrix field:
	\begin{align*}
		V(x)=\begin{pmatrix}
			0 & -(-\Delta)^{\frac{\gamma}{2}}\delta_x\\
			(-\Delta)^{\frac{\gamma}{2}}\delta_x & 0
		\end{pmatrix}.
	\end{align*}

	Let $b:\R^2\to\R^2$ be the kernel given by
	\begin{align*}
		b_i(x)=\sum_{j=1}^2\partial_{x_j}V_{i,j}(x),\quad i=1,2.
	\end{align*}

	Above, the derivatives are meant to be taken in a weak sense. It is known that $\delta\in B^{-(d-\frac{d}{p})}_{p,\infty}$ for any $p\in[1,+\infty]$ (see, for example, the embedding property \eqref{besov-prop-e3} in Section \ref{section-tools}). In particular, $\delta\in B^{0}_{1,\infty}$ and $\delta\notin B^{\varepsilon}_{1,\infty}$ for any $\varepsilon>0$\footnote{This can be seen using, for example, the thermic characterization of Besov spaces introduced in \eqref{def-besovhermic} in Section \ref{section-tools}. It is enough to observe that taking $\alpha=2$, $n=1$ for $\varepsilon$ small enough, $\|\partial_v\tilde{p}_v\ast\delta\|_{L^1}=\|\partial_v\tilde{p}_v\|_{L^1}\geq\tilde{c}v^{-1}$, $\tilde{c}>0$, where the bound follows from the basic computations for the Gaussian heat kernel. And therefore, $\sup_{v\in(0,1]}v^{1-\frac{\varepsilon}{2}}\|\partial_v\tilde{p}_v\ast\delta\|_{L^1}\geq \tilde{c}\sup_{v\in(0,1]}v^{-\frac{\varepsilon}{2}}=+\infty$ and $\|\delta\|_{B^{\varepsilon}_{1,\infty}}=+\infty$.}. Therefore,
	\begin{align*}
		(-\Delta)^{\frac{\gamma}{2}}\delta\in B^{-\gamma}_{1,\infty},\quad\gamma\in(0,\alpha-1),
	\end{align*}
	and $(-\Delta)^{\frac{\gamma}{2}}\delta\notin B^{0}_{1,\infty}$. It is then clear (for example, from the thermic characterization of Besov space given by \eqref{def-besovhermic} in Section \ref{section-tools}; see also \cite[Section 5.2]{deraynal2023multidimensionalstabledrivenmckeanvlasov} for illustrative calculations) that
	\begin{align*}
		b\in B^{-1-\gamma}_{1,\infty}\quad -1-\gamma\in(-\alpha,-1).
	\end{align*}

	From the structure of the matrix $V$ it is also clear that $\Div(b)=0$ which in particular means that $\Div(b)\in B^{-1-\gamma}_{1,\infty}$. With $\beta=-1-\gamma>-\alpha$, $r=+\infty$, $p=1$ and $d=2$, condition \eqref{C3} of Theorem \ref{thm-main-sde} then reads as
	\begin{align*}
		\theta<\frac{\gamma}{2},\quad -1-\gamma&>-\alpha+\Big(-\theta+(1+\gamma)+2-\zeta_0\Big)_+.
	\end{align*}
	Then, for this inequality to be true for $\alpha\in(1,2]$ and $\gamma\in(0,\alpha-1)$, we need
	\begin{align*}
		\zeta_0\geq3-\theta+\gamma>3+\frac{\gamma}{2}\implies\Big(-\theta+(1+\gamma)+2-\zeta_0\Big)_+<\alpha-1-\gamma.
	\end{align*}

	Obviously, the chosen parameters $r,\alpha,\theta$ satisfy the condition \eqref{good-relation} needed for well-posedness of the martingale problem as well as the condition \eqref{good-relation-enhaced} needed for weak well-posedness. A similar analysis can be carried out for the long-time well-posedness given by Theorem \ref{thm-main-sde-lt} by tuning the parameters accordingly.
\end{remark}

\begin{remark}
	Last but not least, let us mention that we consider equation \eqref{main-sde} starting at initial time $t=0$. Note that in \cite{deraynal2022multidimensional,deraynal2023multidimensionalstabledrivenmckeanvlasov}, the works to which we strongly relate, the authors considered a general initial time $t\geq0$. Although a general initial time can be easily introduced in our framework, we set it to zero to lighten the notation.
\end{remark}

\subsection{Strategy of the proof}
\label{subsection-strategy}

To prove Theorem \ref{thm-main-sde} and Theorem \ref{thm-main-sde-lt}, we use the methodology of \cite{deraynal2022multidimensional} and \cite{deraynal2023multidimensionalstabledrivenmckeanvlasov}, which heavily relies on the corresponding non-linear PDE analysis. We briefly describe it here.

\vspace{1mm}

The general strategy consists of the following steps. We first introduce a mollified SDE, in which we mollify the singular kernel $b$ so that the corresponding SDE is well-posed and the time-marginals of the solution are absolutely continuous w.r.t. the Lebesgue measure (see discussion below). We prove a priori estimates for the mollified density, using the fact that it is a mild solution to the corresponding Fokker-Planck equation. This includes proving uniform in the mollification parameter estimates on the mollified density and, consequently, showing that the sequence of mollified densities forms a Cauchy sequence in an appropriate Banach space. We can then pass to the well-posedness of the limit Fokker-Planck equation and obtain a priori estimates for its solution in a suitable space; see Section \ref{section-pde}.

\vspace{1mm}

This somewhat classical approach allows us to show the existence of a solution to the martingale problem associated with the McKean-Vlasov SDE \eqref{main-sde} through the tightness argument and identification of the limit points. Note that since we consider distributional drifts, the martingale problem is understood in the sense of Definition \ref{def-martingale-problem}. Under condition \eqref{good-relation}, we then obtain uniqueness of the martingale solution thanks to \cite[Theorem 3]{deraynal2022-linear-besov}. See Section \ref{section-sde}.

\vspace{1mm}

The existence and uniqueness of a weak solution are not immediate from the well-posedness of the martingale problem in our framework. Indeed, in a classical setup when drift is a function, one recovers noise by subtracting the integrated drift from the canonical process obtained from the martingale solution, see e.g. \cite{stroock1997multidimensional}. However, when the drift is a distribution, the main challenge is to give a meaning to the integrated drift, and therefore, the noise cannot be recovered as before. In this case, it is useful to keep track of the noise when dealing with a martingale problem. This leads to considering an \emph{enlarged} martingale problem. See discussion on Page 9 of \cite{deraynal2022-linear-besov} as well as \cite{DelarueDiel2016RoughPaths}. Following the approach of \cite{deraynal2022-linear-besov}, we reconstruct the integrated drift, with a Lebesgue-Besov drift, as a stochastic non-linear Young integral, see Definition \ref{def-young-integral}. Weak well-posedness then follows directly from \cite[Theorem 6]{deraynal2022-linear-besov}.

\vspace{1mm}

Now, let us introduce some notation and equations based on the preceding discussion. For a Lebesgue-Besov kernel $b$ with negative spatial regularity and $s\in[0,T]$, introduce the following notation:
\begin{align*}
	\mathcal{B}_\nu(s,\cdot):=b(s,\cdot)\ast\nu(\cdot),
\end{align*}
where $\nu$ is any measure for which the convolution is meaningful. Let $b^\varepsilon$, $\varepsilon>0$, be a time-space mollified drift (see Proposition 3 in \cite{deraynal2022multidimensional} for a precise construction). It is smooth and bounded in time and space (see also Proposition \ref{prop-b-conv}). Denote
\begin{align*}
	\mathcal{B}^\varepsilon_\nu(s,\cdot):=b^\varepsilon(s,\cdot)\ast\nu(\cdot).
\end{align*}

Such a drift is smooth and bounded for any probability measure $\nu$. We consider the mollified McKean-Vlasov SDE for $\varepsilon>0$,
\begin{align}
	\label{main-sde-mollified}
	X^{\varepsilon,\mu}_s=\xi+\int_0^s\mathcal{B}^{\varepsilon}_{\bmu^{\varepsilon,\mu}_r}(r,X^{\varepsilon,\mu}_r)dr+\mathcal{W}_{s},\quad\bmu^{\varepsilon}_s=\text{Law}(X^{\varepsilon,\mu}_s),\ \text{Law}(\xi)=\mu.
\end{align}

It is known that for every $\varepsilon>0$, equation \eqref{main-sde-mollified} admits a unique weak (and strong) solution whose time marginal distributions admit a density $\brho^\varepsilon_{\mu}(s,\cdot)$ for any $s\in(0,T]$, for $\alpha\in(1,2]$ (see Section 1.3 in \cite{deraynal2022multidimensional}). By Lemma 3 in \cite{deraynal2022multidimensional}, for every $\varepsilon>0$, for all $s\in(0,T]$ and $y\in\R^d$, the function $\brho^{\varepsilon}_{\mu}(s,y)$ satisfies the Duhamel representation
\begin{align}
	\label{duhamel-main-mollified}
	\brho^{\varepsilon}_{\mu}(s,y)=p^\alpha_{s}\ast\mu(y)-\int_0^s\big((\mathcal{B}^{\varepsilon}_{\brho^{\varepsilon}_{\mu}}(v,\cdot)\brho^{\varepsilon}_{\mu}(v,\cdot))\ast\nabla p^\alpha_{s-v}\big)(y)dv,
\end{align}
where $p^\alpha$ is the density of the driving noise $\mathcal{W}$ with the generator denoted by $\mathcal{L}^\alpha$. Moreover, $\brho^{\varepsilon}_{\mu}$ is a mild solution to the non-linear mollified Fokker-Planck equation
\begin{align}
	\label{pde-main-mollified}
	\begin{cases}
		\partial_s\brho^{\varepsilon}_{\mu}(s,y)+\Div(\mathcal{B}^{\varepsilon}_{\brho^{\varepsilon}_{\mu}}\brho^{\varepsilon}_{\mu}(s,y))-\mathcal{L}^\alpha\brho^{\varepsilon}_{\mu}(s,y)=0, \\
		\brho^{\varepsilon}_{\mu}(0,\cdot)=\mu.
	\end{cases}
\end{align}

In Section \ref{section-pde}, we prove that for any decreasing subsequence $(\varepsilon_k)_{k\geq1}$, $\brho^{\varepsilon_k}_{\mu}(s,\cdot)$ converges to $\brho_{\mu}(s,\cdot)$ in some (appropriate for SDE analysis) space, which, in turn, allows us to take the limit in the Duhamel representation \eqref{duhamel-main-mollified} to obtain that this limit satisfies
\begin{align}
	\label{duhamel-main}
	\brho_{\mu}(s,y)=p^\alpha_{s}\ast\mu(y)-\int_0^s\big((\mathcal{B}_{\brho_{\mu}}(v,\cdot)\brho_{\mu}(v,\cdot))\ast\nabla p^\alpha_{s-v}\big)(y)dv,
\end{align}
and it is a mild solution to the non-linear Fokker-Planck equation
\begin{align}
	\label{pde-main}
	\begin{cases}
		\partial_s\brho_{\mu}(s,y)+\Div(\mathcal{B}_{\brho_{\mu}}\brho_{\mu}(s,y))-\mathcal{L}^\alpha\brho_{\mu}(s,y)=0, \\
		\brho_{\mu}(0,\cdot)=\mu.
	\end{cases}
\end{align}

With this theory in hand, we can construct a martingale solution to a martingale problem linked to \eqref{main-sde} as a limit of a martingale problem associated with \eqref{main-sde-mollified}. Namely, denoting by $\Prob^\varepsilon$ a probability measure on $C([0,S],\R^d)$ if $\alpha=2$ or on $D([0,S],\R^d)$ (space of càdlàg functions) if $\alpha\in(1,2)$, $\Prob^\varepsilon$ is a martingale solution to \eqref{main-sde-mollified} in the sense of Definition \ref{def-martingale-problem}. Provided that $\Prob^\varepsilon(s,dx)=\brho^\varepsilon_{\mu}(s,x)dx$ satisfies appropriate estimates, we deduce the existence of a corresponding limit martingale problem associated with \eqref{main-sde}. Then, suitable PDE theory for the limit Fokker-Planck equation \eqref{pde-main} and the results of \cite{deraynal2022-linear-besov} allow us to conclude that the martingale solution is unique.

\vspace{1mm}

Let us briefly describe the main result of \cite{deraynal2022-linear-besov} used here when $\theta\in(0,\frac{1}{2})$. Therein, well-posedness of a \emph{linear} SDE with a distributional Lebesgue-Besov drift has been proven. More precisely, given the \emph{formal} stochastic differential equation
\begin{align*}
	dX_s=\mathcal{B}(s,X_s)ds+d\mathcal{W}_s,\quad X_0\sim\mu\in\mathcal{P}(\R^d),
\end{align*}
the associated martingale problem (in the sense of Definition \ref{def-martingale-problem}) admits a unique solution if
\begin{align*}
	\mathcal{B}\in L^{r_\theta}((0,T],B^{-\theta}_{p_\theta,q_\theta}),
\end{align*}
where
\begin{align*}
	\alpha\in\Bigg(\frac{1-[1/p_\theta]}{1-[1/r_\theta]},2\Bigg],\quad \theta\in\Bigg(0,\frac{\alpha-1-[\alpha/p_\theta]-[\alpha/r_\theta]}{2}\Bigg).
\end{align*}
Moreover, the formal SDE \eqref{main-sde} has a unique weak solution (in the sense of Definition \ref{def-weak-solution}) that is also pathwise unique in dimension 1 if $\theta$ satisfies the reinforced condition
\begin{align*}
	\theta\in\Bigg(0,\frac{\alpha-1-[2\alpha/p_\theta]-[2\alpha/r_\theta]}{2}\Bigg).
\end{align*}

Following these results and understanding $\mathcal{B}=b\ast\boldsymbol{\mu}^{\mu}$ as a \emph{linear} distributional drift, where $\boldsymbol{\mu}^{\mu}$ is the limit of $\Prob^\varepsilon$, we can show that
\begin{align}
	\label{ineq-drift-liear-bound}
	\|b\ast\boldsymbol{\mu}^{\mu}\|_{L^{r_\theta}(B^{-\theta}_{\infty,\infty})}<+\infty,
\end{align}
where $\theta$ satisfy \eqref{good-relation}, $p_\theta=\infty$, $q_\theta=\infty$ and $r_\theta\in(\alpha,r)$ for martingale problem or $r_\theta\in(2\alpha,r)$ for weak well-posedness. The chosen lower bound on $r_\theta$ guarantees that $\theta=\frac{1}{2}-\varepsilon$ is reachable, whereas the upper bound is needed to meet particular integrability conditions of time singularity, see Lemma \ref{lemma-mollified-drift-bound}. Thus, conditions on $r_\theta$ for proving well-posedness are translated to conditions on $r$ (see \eqref{good-relation} and \eqref{good-relation-enhaced}). Herein, we also include the case $\theta=0$, which was not covered in \cite{deraynal2022-linear-besov} and was rather developed in \cite{deraynal2023multidimensionalstabledrivenmckeanvlasov}.

\vspace{1mm}

For proving Theorem \ref{thm-main-sde-lplq}, we note that it is enough to adapt the results of \cite{deraynal2023multidimensionalstabledrivenmckeanvlasov} for $\beta<-1$ instead of $\beta=-1$ under almost identical conditions in a short time. This comes from a margin in the integrability condition in the a priori estimate of the mollified density, which allows us to balance between integrability in time and negative regularity in space of the kernel. In other words, we derive the exhaustive bounds on the parameters in the classical, non-distributional setting.

\section{Auxiliary tools}
\label{section-tools}

In this section, we list the main properties of Besov spaces used to prove Theorems \ref{thm-main-sde} and \ref{thm-main-sde-lt} including various product rules; see Proposition \ref{prop-product-rule-usual}. We refer the reader to \cite{deraynal2022multidimensional} for proofs.

\begin{remark}
	Define
	\begin{align*}
		\Theta:=\{\alpha,d,r,p,q,\beta,\beta_0,p_0,q_0,\bar{\theta},\|b\|_{L^r(B^\beta_{p,q})},\|\Div(b)\|_{L^r(B^\beta_{p,q})},\|\mu\|_{B^{\beta_0}_{p_0,q_0}},\|\mu\|_{B^{\bar{\theta}}_{1,\infty}}\}
	\end{align*}
	Throughout the whole work, for simplicity, we use a generic constant denoted by $C:=C(\Theta)>0$ that may change from line to line but can always be specified explicitly. It is defined by inequalities in the Besov space, as well as by parameters appearing in \eqref{C3}. Importantly, this constant does not depend on time.
\end{remark}

\vspace{3mm}

First, let us recall the definition of Besov space. We denote by $B^\gamma_{\ell,m}:=B^\gamma_{\ell,m}(\R^d)$ the \emph{Besov space} with regularity index $\gamma\in\R$ and integrability parameters $\ell,m\in[1,+\infty]$. For $\mathcal{S}'(\R^d)$ the dual space of the Schwartz class $\mathcal{S}(\R^d)$, the \emph{thermic characterization} of Besov space is given by
\begin{align}
	\label{def-besovhermic}
	B^\gamma_{\ell,m}=\big\{f\in\mathcal{S}'(\R^d):\ \|f\|_{B^\gamma_{\ell,m}}:=\|\mathcal{F}^{-1}(\phi\mathcal{F}(f))\|_{L^\ell}+\mathcal{T}^\gamma_{\ell,m}(f)<+\infty \big\},
\end{align}
where $\mathcal{T}^\gamma_{\ell,m}(f)$ is called the \emph{thermic part} of the norm and is defined as follows:
\begin{align*}
	\mathcal{T}^\gamma_{\ell,m}(f):=
	\begin{cases}
		\displaystyle\Big(\int_0^1\frac{dv}{v}v^{(n-\frac{\gamma}{\alpha})m}\|\partial_v^n \tilde{p}^{\alpha}_v\ast f\|^m_{L^\ell}\Big)^{1/m},\quad & 1\leq m<+\infty, \\
		\\[0pt]
		\displaystyle\sup_{v\in(0,1]}(v^{(n-\frac{\gamma}{\alpha})}\|\partial_v^n \tilde{p}^{\alpha}_v\ast f\|_{L^\ell}),\quad                      & m=\infty.
	\end{cases}
\end{align*}
Here, $n\in(\frac{\gamma}{\alpha},+\infty)\cap\Z_+$, the function $\phi\in\mathcal{C}^\infty_0(\R^d)$ is s.t. $\phi(0)\neq0$, and $\tilde{p}^{\alpha}_v$ is the density function at time $v>0$ of the $d$-dimensional isotropic $\alpha$-stable process. We refer to \cite[Section 5.3]{LemariRieusset2002RecentDI} for more information about the thermic characterization of Besov spaces.

\vspace{3mm}

Now, we list some properties of Besov spaces that are commonly used in this work.

\vspace{1mm}

\noindent
\textbf{Continuous embeddings.}
\begin{enumerate}[label=(\arabic*)]
	\item(\cite[Proposition 2.1]{Sawano2018BesovSpaces}) For any $\ell\in[1,+\infty]$,
	      \begin{align}
		      \label{besov-prop-e1}
		      B^0_{\ell,1}\hookrightarrow L^\ell\hookrightarrow B^0_{\ell,\infty}.
	      \end{align}
	\item(\cite[Proposition 2.2]{Sawano2018BesovSpaces}) For any $\ell_0,\ell_1,m_0,m_1\in[1,\infty]$ s.t. $\ell_0\leq\ell_1$, $m_0\leq m_1$ and $\gamma_1-d/\ell_1\leq \gamma_0-d/\ell_0$.
	      \begin{align}
		      \label{besov-prop-e2}
		      B^{\gamma_0}_{\ell_0,m_0}\hookrightarrow B^{\gamma_1}_{\ell_1,m_1}.
	      \end{align}
	\item(\cite[Lemma 5]{deraynal2022multidimensional}) For any $\varepsilon>0$, $\ell\in[1,+\infty]$, $m\in[1,+\infty)$,
	      \begin{align}
		      \label{besov-prop-e3}
		      \mathcal{P}(\R^d)\hookrightarrow\cap_{\ell'\geq1}B^{-d/\ell'}_{\ell,\infty}, \quad \mathcal{P}(\R^d)\hookrightarrow\cap_{\ell'\geq1}B^{-d/\ell'-\varepsilon}_{\ell,m}.
	      \end{align}
	\item(\cite[Theorem 2.4]{Sawano2018BesovSpaces}) For any $\ell,m\in[1,+\infty]$, $\gamma\in\R$,
	      \begin{align}
		      \label{besov-prop-e4}
		      \mathcal{S}(\R^d)\hookrightarrow B^\gamma_{\ell,m}\hookrightarrow\mathcal{S}'(\R^d).
	      \end{align}
\end{enumerate}

\noindent
\textbf{Young's convolution inequality} (\cite[Theorem 3]{Burenkov1990BesovYoungIneq}). Let $\ell,m\in[1,\infty]$, $\gamma\in\R$. For any $\delta\in\R$, $\ell_1,\ell_2\in[1,+\infty]$ s.t. $1+\ell^{-1}=\ell_1^{-1}+\ell_2^{-1}$ and $m_1,m_2\in(0,\infty]$ s.t. $m_1^{-1}\geq(m^{-1}-m_2^{-1})\vee0$, there exists a constant $c=c(d)>0$ s.t.
\begin{align}
	\label{besov-prop-y}
	\|f\ast g\|_{B^\gamma_{\ell,m}}\leq c\|f\|_{B^{\gamma-\delta}_{\ell_1,m_1}}\|g\|_{B^{\delta}_{\ell_2,m_2}}.
\end{align}

\noindent
\textbf{Heat kernel Besov norm bound} (\cite[Lemma 11, 12, (3.19)]{deraynal2022-linear-besov}). For any $\ell,m\in[1,\infty]$ and $\gamma\in\R$ such that $\gamma\neq -d(1-\frac{1}{\ell})$, there exists a constant $c=c(\alpha,\ell,m,\gamma,d)>0$ s.t. for any multi-index $\mathbf{a}\in\N^d$ with $|\mathbf{a}|\leq1$ and $0\leq v<s<+\infty$ s.t. $s-v\leq1$,
\begin{align}
	\label{besov-prop-hk}
	\|\partial^\mathbf{a}p^\alpha_{s-v}\|_{B^\gamma_{\ell,m}}\leq\frac{c}{(s-v)^{(\frac{\gamma}{\alpha}+\frac{d}{\alpha}(1-\frac{1}{\ell}))_++\frac{|\mathbf{a}|}{\alpha}}},
\end{align}
where $\partial^\mathbf{a}=\frac{\partial^{|\mathbf{a}|}}{\partial x_1^{\alpha_1}...\partial x_d^{\alpha_d}}$, $|\mathbf{a}|=a_1+...+a_d$. For the case when $s-v>1$, we refer to Section \ref{section-sde-lt}, \eqref{besov-prop-hk-lt}.

\noindent
\textbf{Duality inequality} (\cite[Proposition 3.6]{LemariRieusset2002RecentDI}).
For $\ell,m\in[1,+\infty]$, $\gamma\in\R$ and $f\in B^\gamma_{\ell,m}$, $g\in B^{-\gamma}_{\ell',m'}$, where $\ell',m'$ are conjugates of $\ell,m$ respectfully (i.e., $1/\ell+1/\ell'=1$, $1/m+1/m'=1$),
\begin{align}
	\label{besov-prop-d}
	|\int_{\R^d}f(y)g(y)dy|\leq\|f\|_{B^{\gamma}_{\ell,m}}\|g\|_{B^{-\gamma}_{\ell',m'}}.
\end{align}

\noindent
\textbf{Lift operator} (\cite[Theorem 2.3]{Sawano2018BesovSpaces}). For $\ell,m\in[1,+\infty]$, $\gamma\in\R$, there exists $c>0$ s.t.
\begin{align}
	\label{besov-prop-l}
	\|\nabla f\|_{B^{\gamma-1}_{\ell,m}}\leq c\|f\|_{B^{\gamma}_{\ell,m}}.
\end{align}

\noindent
\textbf{Smooth approximation of the interaction kernel.}
\begin{proposition}[\cite{deraynal2022multidimensional}, Proposition 2]
	\label{prop-b-conv}
	Let $b\in L^r((0,T],B^{\beta}_{p,q})$ with $\beta\leq0$, $p,q\in[1,+\infty]$. There exists a sequence of time-space smooth bounded functions $(b^\varepsilon)_{\varepsilon>0}$ s.t.
	\begin{align*}
		\|b-b^\varepsilon\|_{L^{\bar{r}}((0,T],B^{\bar{\beta}}_{p,q})}\xrightarrow[\varepsilon\to0]{}0,\quad \forall\bar{\beta}<\beta,
	\end{align*}
	where $\bar{r}\in[1,\infty)$ s.t. $\bar{r}=r$ if $r<+\infty$. Moreover, there exists $\kappa\geq1$ s.t.
	\begin{align*}
		\sup_{\varepsilon>0}\|b^\varepsilon\|_{L^{\bar{r}}((0,T],B^{\bar{\beta}}_{p,q})}\leq\kappa\|b\|_{L^{\bar{r}}((0,T],B^{\bar{\beta}}_{p,q})}.
	\end{align*}
\end{proposition}

Note that in \cite{deraynal2022multidimensional} the result is stated for $\beta\in(-1,0]$, however, the proof easily extends to $\beta\in(-\alpha,0]$.

\vspace{1mm}

\noindent
\textbf{Product rules.} For $f,g\in\mathcal{S}'$, we define the product $f\cdot g$ as
\begin{align*}
	f\cdot g=\lim_{j\to\infty}\mathcal{F}^{-1}\big(\phi(2^{-j}\xi)\mathcal{F}(f)(\xi)\big)\times\mathcal{F}^{-1}\big(\phi(2^{-j}\xi)\mathcal{F}(g)(\xi)\big)\quad \text{ in } \mathcal{S}',
\end{align*}
where $\phi\in\mathcal{C}_0^\infty$ with $\phi(x)=\begin{cases}
		1,\ |x|\leq1, \\
		0,\ |x|\geq\frac{3}{2}.
	\end{cases}$

We refer the reader to \cite[Section 4.2]{RunstSickel1996} for more details.

\vspace{1mm}
\noindent
\begin{proposition}[Usual products of Besov spaces]~
	\label{prop-product-rule-usual}
	\begin{itemize}
		\item\emph{\textbf{(positive regularity, \cite{deraynal2023multidimensionalstabledrivenmckeanvlasov})}} For $\lambda\in\R_+$, $\ell,\ell_1,\ell_2\in[1,+\infty]$ s.t. $\frac{1}{\ell}=\frac{1}{\ell_1}+\frac{1}{\ell_2}$ and for $f\in B^\lambda_{\ell_1,\infty}$, $g\in B^{\lambda}_{\ell_2,1}$, there exists $c>0$ s.t.
		      \begin{align}
			      \label{besov-prop-pr1}
			      \|f\cdot g\|_{B^\lambda_{\ell,\infty}}\leq c\|f\|_{B^\lambda_{\ell_1,\infty}}\|g\|_{B^{\lambda}_{\ell_2,1}}.
		      \end{align}

		\item\emph{\textbf{(any regularity, \cite{Sawano2018BesovSpaces})}} For $\lambda\in\R$ and $\rho>\max(\lambda,-\lambda)$, $\ell,m\in[1,+\infty]$ and for $f\in B^\lambda_{\ell,m}$, $g\in B^\rho_{\infty,\infty}$, there exists $c>0$ s.t.
		      \begin{align}
			      \label{besov-prop-pr2}
			      \|f\cdot g\|_{B^\lambda_{\ell,m}}\leq c\|f\|_{B^\rho_{\infty,\infty}}\|g\|_{B^{\lambda}_{\ell,m}}.
		      \end{align}
		
		\item\emph{\textbf{(general regularity indices, \cite{RunstSickel1996}, p. 171)}} For $\lambda,\lambda_1,\lambda_2\in\R_+$, $\ell\in[1,+\infty)$, such that $\lambda_1<\lambda_2$, $\lambda_1\leq\lambda$ and for $f\in B^{-\lambda_1}_{\infty,\infty}$, $g\in B^{\lambda_2}_{\ell,1}$, there exists $c>0$ s.t.
			  \begin{align}
			      \label{ineq-product-rule}
			      \|f\cdot g\|_{B^{-\lambda}_{\ell,\infty}}\leq c\|f\|_{B^{-\lambda_1}_{\infty,\infty}}\|g\|_{B^{\lambda_2}_{\ell,1}}.
			  \end{align}
	\end{itemize}
\end{proposition}

\begin{remark}
	We remark that the product rule \eqref{ineq-product-rule} is stated in \cite{RunstSickel1996} in whole generality. However, here we state only a particular case for the set of parameters of our interest. Similarly, the product rule \eqref{besov-prop-pr2} can be seen as a particular case of the one in \cite{RunstSickel1996}, but is stated separately herein to keep the use case of each clear.
\end{remark}

Let us state an easy but useful for the upcoming computations lemma the proof of which follows from direct calculations.
\begin{lemma}[Beta function upper bound]
	\label{lemma-beta-function}
	For $t,r>0$, $t\leq r$, and $\gamma_1<1$, $\gamma_2<1$, there exists $C>0$ such that
	\begin{align*}
		\int_t^r(r-s)^{-\gamma_1} s ^{-\gamma_2}ds=C r^{1-\gamma_1-\gamma_2}.
	\end{align*}
\end{lemma}

\textbf{Weighted Lebesgue-Besov space.}
For notational convenience, we introduce the so-called \emph{weighted Lebesgue-Besov} spaces. Namely, let $\gamma\in\R$, $\ell,m,r\in[1,+\infty]$, $\omega\geq0$, and $S\in[0,T]$. If $r<+\infty$, we define
\begin{align}
	\label{def-leb-bes-weight}
	\begin{split}
		&L^r_\omega((0,S],B^\gamma_{\ell,m})\\
		&:=\big\{f:s\in[0,S]\mapsto f(s,\cdot)\in B^\gamma_{\ell,m}\text{ measurable and s.t. }\int_0^Ss^{r\omega}\|f(s,\cdot)\|^r_{B^\gamma_{\ell,m}}ds<+\infty\big\}.
	\end{split}
\end{align}

\noindent
If $r=+\infty$,
\begin{align*}
	\begin{split}
		&L^\infty_\omega((0,S],B^\gamma_{\ell,m})\\
		&:=\big\{f:s\in[0,S]\mapsto f(s,\cdot)\in B^\gamma_{\ell,m}\text{ measurable and s.t. }\esssup_{s\in(0,S]}s^{\omega}\|f(s,\cdot)\|_{B^\gamma_{\ell,m}}<+\infty\big\}.
	\end{split}
\end{align*}

The space $\big(L^r_\omega((0,S],B^\gamma_{\ell,m}),\|\cdot\|_{L^r_\omega((0,S],B^\gamma_{\ell,m})}\big)$ is also a Banach space (see \cite{BanachSpaces}, Chapter 1), where
\begin{align*}
	\|f\|_{L^r_\omega((0,S],B^\gamma_{\ell,m})}=\int_0^Ss^{r\omega}\|f(s,\cdot)\|^r_{B^\gamma_{\ell,m}}ds,\quad \|f\|_{L^\infty_\omega((0,S],B^\gamma_{\ell,m})}=\esssup_{s\in(0,S]}s^{\omega}\|f(s,\cdot)\|_{B^\gamma_{\ell,m}}.
\end{align*}

\section{Well-posedness and regularity of the non-linear Fokker-Planck equation in short time}
\label{section-pde}

This section is devoted to proving the following important result about the Fokker-Planck equation \eqref{pde-main} in \emph{small} time.

\begin{proposition}
	\label{prop-main-pde}
	Assume that condition \eqref{C3} is satisfied. Then, for any $\mu\in\mathcal{P}(\R^d)\cap B^{\beta_0}_{p_0,q_0}\cap B^{\bar{\theta}}_{1,\infty}$ with $\bar{\theta}>\theta$, there exists a time horizon $\mathcal{T}\in(0,T]$ with $\mathcal{T}\leq1$ small enough such that the non-linear Fokker-Planck equation \eqref{duhamel-main} admits a unique solution $\brho_{\mu}$ belonging to $L^\infty_\gamma((0,S],B^{-\beta-\theta}_{p',1})$ for any $S\in(0,\mathcal{T}]$, with $\gamma$ defined in Theorem \ref{thm-main-sde}.

	\vspace{1mm}
	\noindent
	Moreover, for all $s\in[0,S]$, $\brho_{\mu}(s,\cdot)\in\mathcal{P}(\R^d)$, and for a.e. $s\in(0,S]$, $\brho_{\mu}(s,\cdot)$ is absolutely continuous w.r.t. the Lebesgue measure and satisfies the Duhamel representation \eqref{duhamel-main}.
\end{proposition}

\begin{remark}
	We carefully notice that the condition $\mu\in B^{\bar{\theta}}_{1,\infty}$ above is needed to prove that $\brho_{\mu}(s,\cdot)$ is a probability measure for $s\in[0,S)$. We refer to Lemma \ref{lemma-density-conv} and Remark \ref{remark-mu-reg}.
\end{remark}

Proposition \ref{prop-main-pde} is the key a priori result that allows us to study the well-posedness of the SDE. Namely, in the current work, we aim at showing that the law of the McKean-Vlasov process is absolutely continuous and its density $\brho_{\mu}$ is such that
\begin{align}
	\label{prf-conv-drift-space}
	\displaystyle\mathcal{B}_{\brho_{\mu}}=\int_{\R^d}b(\cdot,\cdot-y)\brho_{\mu}(\cdot,y)dy\in L^{r_\theta}(B^{-\theta}_{\infty,\infty}),
\end{align}
for some admissible $\theta\in[0,1/2)$ and $r_\theta\in(\alpha,r)$. This information about the drift allows us to use results about linear SDEs with a distributional drift established in \cite{deraynal2022-linear-besov} as it was discussed on page \pageref{ineq-drift-liear-bound}. Let us demonstrate heuristically how this norm of the convolution drift in \eqref{prf-conv-drift-space} would be controlled with the regularity result on the density given by the proposition above. For $s\in(0,T]$, by the Young convolution inequality \eqref{besov-prop-y} and the embedding \eqref{besov-prop-e2},
\begin{align*}
	\|\displaystyle\mathcal{B}_{\brho_{\mu}}(s,\cdot)\|_{B^{-\theta}_{\infty,\infty}} & \leq C\|b(s,\cdot)\|_{B^{\beta}_{p,q}}\|\brho_{\mu}(s,\cdot)\|_{B^{-\beta-\theta}_{p',1}},
\end{align*}
and by the Hölder inequality applied in the time variable
\begin{align*}
	\|\mathcal{B}_{\brho_{\mu}}\|_{L^{r_\theta}(B^{-\theta}_{\infty,\infty})} & \leq C\|b\|_{L^r(B^{\beta}_{p,q})}\Big(\int_0^s\|\brho_{\mu}(s,\cdot)\|^{\frac{r_\theta r}{r-r_\theta}}_{B^{-\beta-\theta}_{p',1}}ds\Big)^{\frac{r-r_\theta}{r_\theta r}}                   \\
	                                                                            & \leq C\|b\|_{L^r(B^{\beta}_{p,q})}\|\brho_{\mu}\|_{L^{\infty}_\gamma(B^{-\beta-\theta}_{p',1})}\Big(\int_0^s s ^{-\gamma\frac{r_\theta r}{r-r_\theta}}ds\Big)^{\frac{r-r_\theta}{r_\theta r}}.
\end{align*}
One can see that in order to make the integral above finite, a proper tuning of the parameters $r,r_\theta,\theta,\gamma$ is needed. The parameter $\gamma$ is also chosen in a such a way that $\brho_{\mu}\in L^{\infty}_\gamma(B^{-\beta-\theta}_{p',1})$ for a given $b\in L^r(B^{\beta}_{p,q})$. For such $\gamma$, finiteness of the integral holds if $r,r_\theta,\theta$ are such that condition \eqref{good-relation} is satisfied. Additionally, in this case, we would obtain existence and uniqueness of a martingale solution to \eqref{main-sde}. Moreover, if we improve the assumption on $r,r_\theta,\theta$, i.e., if condition \eqref{good-relation-enhaced} is satisfied, then we can recover the dynamics of \eqref{main-sde} by giving the meaning to distributional drift $\int\mathcal{B}_{\brho_{\mu}}dy$ and obtain a weak well-posedness result. We refer to Section \ref{section-sde} for all the related details.

\vspace{1mm}

The strategy for proving Proposition \ref{prop-main-pde} is first to establish a control of the type
\begin{align*}
	c_1(s)(f^\varepsilon_{I}(s))^2-f^\varepsilon_{I}(s)+c_2(s)\geq0,
\end{align*}
where for $s\in I$ with $I\subseteq(0,T]$, $\varepsilon>0$, $f^\varepsilon_{I}$ is given by 
\begin{align}
	\label{def-f-I}
	f^\varepsilon_{I}(s):=\sup_{v\in(0,s]} v^\gamma\|\brho^\varepsilon_{\mu}(v,\cdot)\|_{B^{-\beta-\theta}_{p',1}},
\end{align}
and $c_1(s),c_2(s)$  are some time functions (Lemma \ref{lemma-density-estimate}). Then, by showing that $s\mapsto f^\varepsilon_{I}(s)$ is continuous (Lemma \ref{lemma-norm-density-cont}), we can deduce the uniform in small time bound on $f^\varepsilon_{I}(\cdot)$ (Lemma \ref{lemma-gronwall-control}). We will show that the mollified densities form a Cauchy sequence in $\varepsilon>0$ in the considered weighted Lebesgue-Besov space (Lemma \ref{lemma-density-conv}) which will imply the existence of a limit density (Lemma \ref{lemma-density-existence}). Finally, we will demonstrate that this limit is a unique distributional solution to the non-linear Fokker-Planck equation associated with the McKean-Vlasov SDE \eqref{main-sde} (Lemma \ref{lemma-density-uniqueness}).

\vspace{1mm}

In \cite{deraynal2023multidimensionalstabledrivenmckeanvlasov}, the authors have established well-posedness results of \eqref{main-sde} by showing that the drift $\mathcal{B}_{\brho_{\mu}}$ belongs to an appropriate $L^q_t-L^p_x$ space which, in turn, allows us to invoke the work by Krylov and Röckner \cite{KrylovRockner2005LpLq} as well as \cite{deraynal2022-linear-besov} for the regularity index $\theta=0$. For $\beta=-1$ the authors applied integration by parts in the bilinear term in the Duhamel formulation for $\brho^\varepsilon_{\mu}$ which allowed moving gradient from the heat kernel to the product $\mathcal{B}^\varepsilon_{\brho^\varepsilon_{\mu}}(s,\cdot)\brho^\varepsilon_{\mu}(s,\cdot)$ and therefore to obtain the margin for $\beta$ for going below $-1+\varepsilon$ (this threshold was proved in \cite{deraynal2022multidimensional}). This additionally imposes some conditions on the divergence of the kernel. We note that the product rule in Besov space with \emph{positive} regularity \eqref{besov-prop-pr1} on $\|\mathcal{B}^\varepsilon_{\brho^\varepsilon_{\mu}}(s,\cdot)\brho^\varepsilon_{\mu}(s,\cdot)\|_{B^{\Gamma}_{p',\infty}}$ for some $\Gamma>0$, was crucial to obtain the sharpest regularity bound in the given framework.

\vspace{1mm}

In our work, however, we want to take advantage of the ``linear'' theory of SDEs with a \emph{distributional} drift, namely of regularity $-\theta\in(-\frac{1}{2},0]$, developed in \cite{deraynal2022-linear-besov}. As one can see in Lemma \ref{lemma-density-estimate}, it is not possible to go beyond $\beta=-1$ by applying the product rule \eqref{besov-prop-pr1} or \eqref{besov-prop-pr2}. This is because we have to apply a product rule to the term $\|\mathcal{B}^\varepsilon_{\brho^\varepsilon_{\mu}}(s,\cdot)\brho^\varepsilon_{\mu}(s,\cdot)\|_{B^{-\theta}_{p',\infty}}$, and as one can see it is of negative regularity. This is one of the primary differences w.r.t. \cite{deraynal2023multidimensionalstabledrivenmckeanvlasov}. This singular product appears throughout the whole proof of Proposition \ref{prop-main-pde} in this section. Even though the proofs here often strongly resemble those in \cite{deraynal2023multidimensionalstabledrivenmckeanvlasov}, appearance of the singular product requires more care and sometimes different techniques. This leads to another product rule \eqref{ineq-product-rule}, with the help of which one overcomes the restrictive lower bound on $\beta$. The described approach allows us to recover a martingale solution in the sense of Definition \ref{def-martingale-problem} of \eqref{main-sde} for regular enough initial probability measure and with the interaction kernel having regularity associated with the Besov norm in the range $(-2,-1)$.

\vspace{1mm}

Note that when $\theta=0$, thanks to the embedding \eqref{besov-prop-e1} ($\|f^\varepsilon_{I}(s)\|_{L^{\infty}}\leq\|f^\varepsilon_{I}(s)\|_{B^0_{\infty,1}}$), we are in the framework of \cite{deraynal2023multidimensionalstabledrivenmckeanvlasov}, but the approach of this paper forces us to consider $\beta<-1$ which can be seen from the upper bound on $\beta$ given by \eqref{C3}.

\vspace{3mm}

We begin with proving control on the weighted Lebesgue-Besov norm of the mollified density.

\begin{lemma}[A priori estimates on the mollified density]
	\label{lemma-density-estimate} Let conditions \eqref{C3} and \eqref{mu-range} be satisfied. Let $\theta$ satisfy either \eqref{good-relation} or \eqref{good-relation-enhaced}. Let $\bar{\beta}_0=\beta_0\wedge\frac{\beta_0 p_0'}{p}$, $\bar{p}_0=p_0\wedge p'$, $\bar{q}_0=q_0\vee\frac{pq_0}{p_0'}$. Let $\eta,\delta>0$ be small enough. Denote
	\begin{align}
		\label{def-gamma-weight}
		\gamma:=\frac{\eta}{\alpha}+\frac{1}{\alpha}\Big(-\theta-\beta+\frac{d}{p}-\zeta_0\Big)_+.
	\end{align}
	Then for any $S\leq T$ s.t. $S\leq1$ is small enough, it holds
	\begin{align}
		\label{ineq-quadratic-density}
		\begin{split}
			 & \sup_{v\in(0,S]} v^{\gamma}\|\brho^{\varepsilon}_{\mu}(v,\cdot)\|_{B^{-\beta-\theta}_{p',1}}\leq C\|\mu\|_{B^{\bar{\beta}_0}_{\bar{p}_0,\bar{q}_0}}S^{\frac{\eta}{\alpha}}                                                                \\
			 & \quad+C\Big(\|b\|_{L^r(B^\beta_{p,q})}+\|\Div(b)\|_{L^r(B^\beta_{p,q})}\Big)\big(\sup_{v\in(0,S]} v^{\gamma}\|\brho^{\varepsilon}_{\mu}(v,\cdot)\|_{B^{-\beta-\theta}_{p',1}}\big)^2S^{1-\frac{1}{r}-\gamma+\frac{\beta-\delta}{\alpha}},
		\end{split}
	\end{align}
	where $1-\frac{1}{r}-\gamma+\frac{\beta-\delta}{\alpha}>0$.

	\begin{proof}
		Let $\beta\leq0$, $\theta\in[0,\frac{1}{2})$ satisfying \eqref{good-relation} or \eqref{good-relation-enhaced} and let \eqref{C3} hold. Fix $v\in(0,S]$. Using Duhamel formula \eqref{duhamel-main-mollified}, we write
			\begin{align}
				\label{prf-3-duhamel}
				\|\brho^{\varepsilon}_{\mu}(v,\cdot)\|_{B^{-\beta-\theta}_{p',1}}\leq\|\mu\ast p^{\alpha}_{ v}\|_{B^{-\beta-\theta}_{p',1}}+\int_0^v\|(\mathcal{B}^{\varepsilon}_{\brho^{\varepsilon}_{\mu}}(s,\cdot)\brho^{\varepsilon}_{\mu}(s,\cdot))\ast\nabla p^\alpha_{v-s}\|_{B^{-\beta-\theta}_{p',1}}ds.
			\end{align}

			\textbf{Initial condition.}
			To handle the first term of \eqref{prf-3-duhamel}, we have to distinguish different relations between $p$ and $p_0$.
			\begin{itemize}
				\item($p_0\leq p'$)
				      In this case, we can simply apply the Young inequality \eqref{besov-prop-y} on $\|\mu\ast p^{\alpha}_r\|_{B^{-\beta-\theta}_{p',1}}$,
				      \begin{align*}
					      \|\mu\ast p^{\alpha}_{ v}\|_{B^{-\beta-\theta}_{p',1}}\leq C\|\mu\|_{B^{\beta_0}_{p_0,q_0}}\|p^{\alpha}_{ v}\|_{B^{-\theta-\beta-\beta_0}_{\tilde{p},1}},
				      \end{align*}
				      where $\frac{1}{\tilde{p}}=1+\frac{1}{p'}-\frac{1}{p_0}$, so that the condition $p_0\leq p'$ gives $\frac{1}{\tilde{p}}\in[0,1]$. By the heat kernel estimate \eqref{besov-prop-hk} and from the relation between $\tilde{p},p$ and $p_0$,
				      \begin{align*}
					      \|\mu\ast p^{\alpha}_{ v}\|_{B^{-\beta-\theta}_{p',1}} & \leq C\|\mu\|_{B^{\beta_0}_{p_0,q_0}} v^{-\frac{1}{\alpha}(-\theta-\beta-\beta_0+d(1-\frac{1}{\tilde{p}}))_+}   \\
					                                                              & =C\|\mu\|_{B^{\beta_0}_{p_0,q_0}} v^{-\frac{1}{\alpha}(-\theta-\beta-\beta_0+d(\frac{1}{p_0}-\frac{1}{p'}))_+}  \\
					                                                              & =C\|\mu\|_{B^{\beta_0}_{p_0,q_0}} v^{-\frac{1}{\alpha}(-\theta-\beta-\beta_0+d(\frac{1}{p}-\frac{1}{p_0'}))_+}.
				      \end{align*}

				      Note that for \eqref{besov-prop-hk} to be applicable, it should hold
				      \begin{align}
					      \label{theta-bound1}
					      -\beta-\theta\neq\beta_0-\frac{d}{\tilde{p}'}\iff-\beta-\theta\neq\beta_0-\frac{1}{p}+\frac{1}{p_0'}.
				      \end{align}

				\item($p_0>p'$) By the embedding \eqref{besov-prop-e3}, any probability measure belongs to the Besov space $B^{0}_{1,\infty}$. Then, by interpolation between Besov spaces (see, e.g., \cite{Sawano2018BesovSpaces}), we have
				      \begin{align*}
					      \mu\in[B^{0}_{1,\infty},B^{\beta_0}_{p_0,q_0}]=B^{\frac{\beta_0 p_0'}{p}}_{p',\frac{pq_0}{p_0'}}.
				      \end{align*}
				      Thus, by the Young inequality \eqref{besov-prop-y} and the heat kernel estimate \eqref{besov-prop-hk},
				      \begin{align*}
					      \|\mu\ast p^{\alpha}_{ v}\|_{B^{-\beta-\theta}_{p',1}} & \leq C\|\mu\|_{B^{\frac{\beta_0 p_0'}{p}}_{p',\frac{pq_0}{p_0'}}}\|p^\alpha_{ v}\|_{B^{-\theta-\beta-\frac{\beta_0 p_0'}{p}}_{1,1}} \\
					                                                              & \leq C\|\mu\|_{B^{\frac{\beta_0 p_0'}{p}}_{p',\frac{pq_0}{p_0'}}} v^{-\frac{1}{\alpha}(-\theta-\beta-\frac{\beta_0 p_0'}{p})_+},
				      \end{align*}
				      where again according to \eqref{besov-prop-hk},
				      \begin{align}
					      \label{theta-bound2}
					      -\beta-\theta\neq\frac{\beta_0 p_0'}{p}.
				      \end{align}
			\end{itemize}

			Therefore, the general estimate reads as
			\begin{align}
				\label{prf-ineq-initial-data}
				\|\mu\ast p^{\alpha}_{ v}\|_{B^{-\beta-\theta}_{p',1}} & \leq C\|\mu\|_{B^{\bar{\beta}_0}_{\bar{p}_0,\bar{q}_0}} v^{-\gamma_0},
			\end{align}
			where
			\begin{align}
				\label{def-gamma0}
				\gamma_0=\frac{1}{\alpha}\Big(-\theta-\beta+\frac{d}{p}-\zeta_0\Big)_+.
			\end{align}
			Also, combining \eqref{theta-bound1} and \eqref{theta-bound2}, we have
			\begin{align*}
				\begin{split}
					-\beta-\theta+\frac{d}{p}\neq\bar{\beta}_0.
				\end{split}
			\end{align*}

			\vspace{2mm}

			\textbf{Bilinear operator.}
			Let us now handle the integrand of \eqref{prf-3-duhamel}. First, integration by parts and the Young inequality \eqref{besov-prop-y} yield
			\begin{align}
				\begin{split}
					 & \|(\mathcal{B}^{\varepsilon}_{ \boldsymbol{\rho}^{\varepsilon}_{\mu}}(s,\cdot)\boldsymbol{\rho}^{\varepsilon}_{\mu}(s,\cdot))\ast\nabla p^\alpha_{v-s}\|_{B^{-\beta-\theta}_{p',1}}                                                                                                                                                                                                                            \\
					 & \leq C\Big(\|(\Div(\mathcal{B}^{\varepsilon}_{ \boldsymbol{\rho}^{\varepsilon}_{\mu}}(s,\cdot))\boldsymbol{\rho}^{\varepsilon}_{\mu}(s,\cdot))\ast p^\alpha_{v-s}\|_{B^{-\beta-\theta}_{p',1}}+\|(\mathcal{B}^{\varepsilon}_{ \boldsymbol{\rho}^{\varepsilon}_{\mu}}(s,\cdot)\cdot\nabla\boldsymbol{\rho}^{\varepsilon}_{\mu}(s,\cdot))\ast p^\alpha_{v-s}\|_{B^{-\beta-\theta}_{p',1}}\Big)                \\
					 & \leq C\Big(\|\Div(\mathcal{B}^{\varepsilon}_{ \boldsymbol{\rho}^{\varepsilon}_{\mu}}(s,\cdot))\boldsymbol{\rho}^{\varepsilon}_{\mu}(s,\cdot)\|_{B^{-\theta}_{p',\infty}}+\|\mathcal{B}^{\varepsilon}_{ \boldsymbol{\rho}^{\varepsilon}_{\mu}}(s,\cdot)\cdot\nabla\boldsymbol{\rho}^{\varepsilon}_{\mu}(s,\cdot)\|_{B^{-\theta}_{p',\infty}}\Big)\|p^\alpha_{v-s}\|_{B^{-\beta}_{1,1}}.
					\label{prf-quad-term}
				\end{split}
			\end{align}

			For both appearing terms of the sum in \eqref{prf-quad-term}, we want to invoke the product rule \eqref{ineq-product-rule}. Let us first analyze the first term in \eqref{prf-quad-term}. Following the notation of \eqref{ineq-product-rule}, let $\delta>0$ and
		\begin{gather*}
			\lambda=\theta,\ \lambda_1=\theta,\ \lambda_2=\theta+\delta, \\
			\ell=p'                                                           \\
			F=\Div(\mathcal{B}^{\varepsilon}_{ \boldsymbol{\rho}^{\varepsilon}_{\mu}}(s,\cdot)),\ g=\boldsymbol{\rho}^{\varepsilon}_{\mu}(s,\cdot),
		\end{gather*}
		Note that $\lambda,\lambda_1,\lambda_2\geq0$, $p>1$ fulfill the conditions required for \eqref{ineq-product-rule}. Thus,
		\begin{align}
			\label{prf-quad-1}
			\begin{split}
			\|\Div(\mathcal{B}^{\varepsilon}_{ \boldsymbol{\rho}^{\varepsilon}_{\mu}}(s,\cdot))\boldsymbol{\rho}^{\varepsilon}_{\mu}(s,\cdot)\|_{B^{-\theta}_{p',\infty}}
			 & \leq C\|\Div(\mathcal{B}^{\varepsilon}_{ \boldsymbol{\rho}^{\varepsilon}_{\mu}}(s,\cdot))\|_{B^{-\theta}_{\infty,\infty}}\|\boldsymbol{\rho}^{\varepsilon}_{\mu}(s,\cdot)\|_{B^{\theta+\delta}_{p',1}} \\
			 & \leq C\|\Div(b^{\varepsilon}(s,\cdot))\|_{B^{\beta}_{p,q}}\|\boldsymbol{\rho}^{\varepsilon}_{\mu}(s,\cdot)\|^2_{B^{-\beta-\theta}_{p',1}},
			\end{split}
		\end{align}
		where the last inequality follows from the embedding \eqref{besov-prop-e2} as soon as
		\begin{align*}
			\theta+\delta\leq-\beta-\theta\iff\beta\leq-2\theta-\delta,
		\end{align*}
		and from the Young inequality \eqref{besov-prop-y} applied to $\Div(\mathcal{B}^{\varepsilon}_{ \boldsymbol{\rho}^{\varepsilon}_{\mu}}(s,\cdot))$.

		\vspace{1mm}

		For the second term in \eqref{prf-quad-term}, we apply the product rule \eqref{ineq-product-rule} with the same parameters as above, together with the Young inequality \eqref{besov-prop-y}, the embedding \eqref{besov-prop-e2} and the lifting property \eqref{besov-prop-l}:
		\begin{align}
			\label{prf-quad-2}
			\begin{split}
			\|\mathcal{B}^{\varepsilon}_{ \boldsymbol{\rho}^{\varepsilon}_{\mu}}(s,\cdot)\cdot\nabla\boldsymbol{\rho}^{\varepsilon}_{\mu}(s,\cdot)\|_{B^{-\theta}_{p',\infty}} & \leq C\|\mathcal{B}^{\varepsilon}_{ \boldsymbol{\rho}^{\varepsilon}_{\mu}}(s,\cdot)\|_{B^{-\theta}_{\infty,\infty}}\|\nabla\boldsymbol{\rho}^{\varepsilon}_{\mu}(s,\cdot)\|_{B^{\theta+\delta}_{p',1}} \\
			                                                                                                                                                                             & \leq C\|b^\varepsilon(s,\cdot)\|_{B^\beta_{p,q}}\|\boldsymbol{\rho}^{\varepsilon}_{\mu}(s,\cdot)\|^2_{B^{-\beta-\theta}_{p',1}},
			\end{split}
		\end{align}
		where the last inequality holds as soon as
		\begin{align*}
			1+\theta+\delta\leq-\beta-\theta\iff\beta\leq-1-2\theta-\delta.
		\end{align*}

		\vspace{2mm}

		\textbf{Compensating the time singularity.}
		Combining the estimates \eqref{prf-ineq-initial-data}, \eqref{prf-quad-term}, \eqref{prf-quad-1} and \eqref{prf-quad-2} and recalling that $\gamma_0$ is given by \eqref{def-gamma0}, we get the following control:
		\begin{align*}
			 & \|\boldsymbol{\rho}^{\varepsilon}_{\mu}(v,\cdot)\|_{B^{-\beta-\theta}_{p',1}}                                                                                                                                                                                                                                               \\
			 & \leq C\|\mu\|_{B^{\bar{\beta}_0}_{\bar{p}_0,\bar{q}_0}} v^{-\gamma_0}+C\int_0^v\Big(\|b^\varepsilon(s,\cdot)\|_{B^{\beta}_{p,q}}+\|\Div(b^\varepsilon(s,\cdot))\|_{B^{\beta}_{p,q}}\Big)\|\boldsymbol{\rho}^{\varepsilon}_{\mu}(s,\cdot)\|^2_{B^{-\beta-\theta}_{p',1}}(v-s)^{\frac{\beta}{\alpha}}ds             \\
			 & \leq C\|\mu\|_{B^{\bar{\beta}_0}_{\bar{p}_0,\bar{q}_0}} v^{-\gamma_0}\\
			 &\quad+C\Big(\|b\|_{L^r((0,T],B^{\beta}_{p,q})}+\|\Div(b)\|_{L^r((0,T],B^{\beta}_{p,q})}\Big)\Big(\int_0^v\|\boldsymbol{\rho}^{\varepsilon}_{\mu}(s,\cdot)\|^{2r'}_{B^{-\beta-\theta}_{p',1}}(v-s)^{\frac{\beta}{\alpha}r'}ds\Big)^{\frac{1}{r'}},
		\end{align*}
		where the last inequality follows from the Hölder inequality in time. Multiplying both sides by $ v^{\gamma}$ for some $\gamma\geq0$ to be specified later and taking the supremum over $(0,S]$, we derive
		\begin{align*}
			&\sup_{v\in(0,S]} v^{\gamma}\|\boldsymbol{\rho}^{\varepsilon}_{\mu}(v,\cdot)\|_{B^{-\beta-\theta}_{p',1}} \\
			& \leq C\|\mu\|_{B^{\bar{\beta}_0}_{\bar{p}_0,\bar{q}_0}}\sup_{v\in(0,S]} v^{\gamma-\gamma_0}                                                                                                      \\
			                                                                                                              & \quad+C\Big(\|b\|_{L^r((0,T],B^{\beta}_{p,q})}+\|\Div(b)\|_{L^r((0,T],B^{\beta}_{p,q})}\Big)\sup_{v\in(0,S]} v^{\gamma}\Big(\int_0^v\|\boldsymbol{\rho}^{\varepsilon}_{\mu}(s,\cdot)\|^{2r'}_{B^{-\beta-\theta}_{p',1}}(v-s)^{\frac{\beta}{\alpha}r'}ds\Big)^{\frac{1}{r'}}.
		\end{align*}
		Let us handle the integral. Multiplying the integrand by $ s ^{2r'\gamma} s ^{-2r'\gamma}$ and applying Lemma \ref{lemma-beta-function}, we obtain
		\begin{align}
			\label{prf-dens-est-integr}
			\begin{split}
				\int_0^v\|\boldsymbol{\rho}^{\varepsilon}_{\mu}(s,\cdot)\|^{2r'}_{B^{-\beta-\theta}_{p',1}}(v-s)^{\frac{\beta}{\alpha}r'}ds & \leq\big(\sup_{s\in(0,v]} s ^{\gamma}\|\boldsymbol{\rho}^{\varepsilon}_{\mu}(s,\cdot)\|_{B^{-\beta-\theta}_{p',1}}\big)^{2r'}\int_0^v(v-s)^{\frac{\beta}{\alpha}r'} s^{-2r'\gamma}ds \\
				                                                                                                                                    & \leq C\big(\sup_{s\in(0,v]} s ^{\gamma}\|\boldsymbol{\rho}^{\varepsilon}_{\mu}(s,\cdot)\|_{B^{-\beta-\theta}_{p',1}}\big)^{2r'} v^{1-2r'\gamma+\frac{\beta}{\alpha}r'},
			\end{split}
		\end{align}
		provided that
		\begin{align*}
			\gamma & <\frac{1}{2r'},\quad \frac{1}{r'}>-\frac{\beta}{\alpha}.
		\end{align*}

		Then,
		\begin{align*}
			&\sup_{v\in(0,S]} v^{\gamma}\|\boldsymbol{\rho}^{\varepsilon}_{\mu}(v,\cdot)\|_{B^{-\beta-\theta}_{p',1}} \\
			& \leq C\|\mu\|_{B^{\bar{\beta}_0}_{\bar{p}_0,\bar{q}_0}}S^{\gamma-\gamma_0}\\
			&\quad+C\Big(\|b\|_{L^r((0,T],B^{\beta}_{p,q})}+\|\Div(b)\|_{L^r((0,T],B^{\beta}_{p,q})}\Big)\big(\sup_{v\in(0,S]} v^{\gamma}\|\boldsymbol{\rho}^{\varepsilon}_{\mu}(v,\cdot)\|_{B^{-\beta-\theta}_{p',1}}\big)^2 S^{\frac{1}{r'}-\gamma+\frac{\beta}{\alpha}}.
		\end{align*}

		To avoid time singularities on the r.h.s., we need a $\gamma$ that satisfies
		\begin{align*}
			\gamma-\gamma_0>0
		\end{align*}
		and
		\begin{align}
			\label{prf-pos-exp1}
			\gamma<\frac{1}{r'}+\frac{\beta}{\alpha}.
		\end{align}
		Using the definition \eqref{def-gamma0} of $\gamma_0$, these conditions translate to
		\begin{align*}
			\gamma_0=\frac{1}{\alpha}\Big(-\theta-\beta+\frac{d}{p}-\zeta_0\Big)_+ & <\gamma<\frac{1}{r'}+\frac{\beta}{\alpha}.
		\end{align*}

		Take $\gamma=\gamma_0+\frac{\eta}{\alpha}$ with $\eta>0$. Then,
		\begin{align*}
			     & \ \gamma=\frac{\eta}{\alpha}+\frac{1}{\alpha}\Big(-\theta-\beta+\frac{d}{p}-\zeta_0\Big)_+<\frac{1}{r'}+\frac{\beta}{\alpha} \\
			\iff & \beta>-\alpha\Big(1-\frac{1}{r}\Big)+\Big(-\theta-\beta+\frac{d}{p}-\zeta_0\Big)_++\eta                                      \\
			\implies & \beta>-\alpha\Big(1-\frac{1}{r}\Big)+\Big(-\theta-\beta+\frac{d}{p}-\zeta_0\Big)_+.
		\end{align*}

		This reads exactly as condition \eqref{C3}, and we have
		\begin{align*}
			&\sup_{v\in(0,S]} v^{\gamma}\|\boldsymbol{\rho}^{\varepsilon}_{\mu}(v,\cdot)\|_{B^{-\beta-\theta}_{p',1}} \\
			& \leq C\|\mu\|_{B^{\bar{\beta}_0}_{\bar{p}_0,\bar{q}_0}} S^{\frac{\eta}{\alpha}}\\
			&\quad+C\Big(\|b\|_{L^r((0,T],B^{\beta}_{p,q})}+\|\Div(b)\|_{L^r((0,T],B^{\beta}_{p,q})}\Big)\big(\sup_{v\in(0,S]} v^{\gamma}\|\boldsymbol{\rho}^{\varepsilon}_{\mu}(v,\cdot)\|_{B^{-\beta-\theta}_{p',1}}\big)^2 S^{\frac{1}{r'}-\gamma+\frac{\beta}{\alpha}},
		\end{align*}
		where $\frac{1}{r'}-\gamma+\frac{\beta}{\alpha}>0$ from \eqref{prf-pos-exp1}.
	\end{proof}
\end{lemma}

\begin{remark}[About the threshold on $\beta$]
	Let us note that the restrictive lower bound on $\beta$ appears when handling the time-weighted Lebesgue norm of $\|\boldsymbol{\rho}^{\varepsilon}_{\mu}(v,\cdot)\|_{B^{-\beta-\theta}_{p',1}}$. More precisely, in order to have an integrable singularity in \eqref{prf-dens-est-integr}, we must have
	\begin{align*}
		\frac{1}{r'}>-\frac{\beta}{\alpha}\iff\beta>-\alpha\Big(1-\frac{1}{r}\Big).
	\end{align*}
	Taking $\alpha=2$, $r=+\infty$, one sees that necessarily $\beta>-2$ despite the regularity of $\mu$. This is thus the limitation of the current approach.
\end{remark}

In the next lemma, we obtain an estimate on $f^\varepsilon_{I}$ defined in \eqref{def-f-I}, for each $\varepsilon>0$. The resulting estimate blows up as $\varepsilon\to0$, however, it provides an a priori control that allows us to get a uniform control in $\varepsilon$ in Lemma \ref{lemma-gronwall-control}.

\begin{lemma}[Continuity of the weighted Lebesgue norm of the mollified density]
	\label{lemma-norm-density-cont}
	For any $\varepsilon>0$, any $I\subseteq(0,T]$, there exist $C_\varepsilon:=C(\varepsilon)$, $\tilde{\gamma}>0$ s.t. for any $s\in(t,1\wedge T]$,
	\begin{align}
		\label{ineq-cont-norm}
		f^\varepsilon_{I}(s)\leq C_\varepsilon s ^{\tilde{\gamma}},
	\end{align}
	where $f^\varepsilon_{I}$ is given by \eqref{def-f-I}. In particular, for any $\varepsilon>0$ and $t\in(0,T]$, $s\mapsto f^\varepsilon_{I}(s)$ is continuous.

	\begin{proof}

		\textbf{Plan of proof.} When showing an estimate of the type \eqref{ineq-cont-norm}, it is not obvious how to obtain an inequality which is not quadratic as in Lemma \ref{lemma-density-estimate} or how to obtain some boundedness of the norm in time. The approach used here consists of building an appropriate iteration process for the regularity index of $\|\brho^{\varepsilon}_{\mu}(s,\cdot)\|_{B^{-\beta-\theta}_{p',1}}$. More precisely, we start from some regularity $\xi\leq-\beta-\theta$ that allows us to show a good estimate on $\|\brho^{\varepsilon}_{\mu}(s,\cdot)\|_{B^{\xi}_{p',1}}$ by exploiting only the fact that $\brho^{\varepsilon}_{\mu}(s,\cdot)$ is a probability density, i.e., in lives in $L^1$. Then, we estimate $\|\brho^{\varepsilon}_{\mu}(s,\cdot)\|_{B^{\xi+j\xi_0}_{p',1}}$ for $j\in\N$ and $\xi_0\geq0$ by using the estimate from the iteration $j-1$. Finally, defining $n\in\N$ so that $\xi+n\xi_0=-\beta-\theta$, we can derive the desired estimate of the type \eqref{ineq-cont-norm} for $\|\brho^{\varepsilon}_{\mu}(s,\cdot)\|_{B^{-\beta-\theta}_{p',1}}$ by exploiting previous iterations' estimates.

		\vspace{2mm}

		First, note that by construction, $b^\varepsilon$ is smooth and bounded, however, it is explosive when $\varepsilon\to0$. Namely, considering a family of smooth compactly supported mollifiers $(\eta_\varepsilon)_{\varepsilon>0}$ $b^\varepsilon$ is given by 
		\begin{align*}
			b^\varepsilon(s,x)=\big(p^\alpha(\varepsilon,\cdot)\ast_x (b\ast_s \eta_\varepsilon)\big)(s,x),\quad(s,x)\in(0,T]\times\R^d,\ \varepsilon>0,
		\end{align*}
		see \cite{deraynal2022multidimensional} for more details. And therefore, for any $\rho\in\R$ such that $\rho-\beta\neq-\frac{d}{p}$,
		\begin{gather*}
			\|b^\varepsilon(s,\cdot)\|_{B^{\rho}_{\infty,\infty}}=\|\big(p^\alpha(\varepsilon,\cdot)\ast_x (b\ast_s \eta_\varepsilon)\big)(s,\cdot)\|_{B^{\rho}_{\infty,\infty}}\\
			\stackrel{\eqref{besov-prop-y}}{\leq} C\|p^\alpha(\varepsilon,\cdot)\|_{B^{\rho-\beta}_{p',1}}\|(b\ast_s \eta_\varepsilon)(s,\cdot)\|_{B^\beta_{p,q}}\stackrel{\eqref{besov-prop-hk}}{\leq} C\varepsilon^{-(\frac{\rho-\beta}{\alpha}+\frac{d}{\alpha p})_+}\big(\|b(\cdot,\cdot)\|_{B^\beta_{p,q}}\ast_s \eta_\varepsilon\big)(s).
		\end{gather*}

		Take the supremum in time:
		\begin{align*}
			\|b^\varepsilon\|_{L^\infty(B^{\rho}_{\infty,\infty})}&\leq C\varepsilon^{-(\frac{\rho-\beta}{\alpha}+\frac{d}{\alpha p})_+}\big\|\|b(\cdot,\cdot)\|_{B^\beta_{p,q}}\ast_s \eta_\varepsilon\big\|_{L^\infty}\\
			&\leq C\varepsilon^{-(\frac{\rho-\beta}{\alpha}+\frac{d}{\alpha p})_+}\|b\|_{L^r(B^\beta_{p,q})}\|\eta_\varepsilon\|_{L^{r'}}.
		\end{align*}
		Thus, there exists a positive but explosive constant $C_\varepsilon$ such that
		\begin{align}
			\label{ineq-B-rho-reg}
			\begin{split}
				\|\mathcal{B}^{\varepsilon}_{\brho^{\varepsilon}_{\mu}}\|_{L^\infty(B^{\rho}_{\infty,\infty})} & \stackrel{\eqref{besov-prop-y}}{\leq} C\|b^\varepsilon\|_{L^\infty(B^{\rho}_{\infty,\infty})}\|\brho^\varepsilon_{\mu}\|_{L^\infty(B^0_{1,\infty})} \\
																										& \stackrel{\eqref{besov-prop-e1}}{\leq} C\|b^\varepsilon\|_{L^\infty(B^{\rho}_{\infty,\infty})}\|\brho^\varepsilon_{\mu}\|_{L^\infty(L^1)}                              \\
																										& \leq C_\varepsilon.
			\end{split}
		\end{align}

		Let $\xi\in\R$. We utilize the Duhamel expansion as in \eqref{prf-3-duhamel}:
		\begin{align*}
			\|\brho^{\varepsilon}_{\mu}(v,\cdot)\|_{B^{\xi}_{p',1}}&\leq\|\mu\ast p^{\alpha}_{v}\|_{B^{\xi}_{p',1}}+C\int_0^v\|(\mathcal{B}^{\varepsilon}_{ \boldsymbol{\rho}^{\varepsilon}_{\mu}}(s,\cdot)\boldsymbol{\rho}^{\varepsilon}_{\mu}(s,\cdot))\ast\nabla p^\alpha_{v-s}\|_{B^{\xi}_{p',1}}ds.
		\end{align*}

		The initial data term is handled similarly to \eqref{prf-ineq-initial-data} up to the change of the regularity index:
		\begin{align*}
			\|\mu\ast p^{\alpha}_{v}\|_{B^{\xi}_{p',1}}\leq C\|\mu\|_{B^{\bar{\beta}_0}_{\bar{p}_0,\bar{q}_0}}v^{-\frac{1}{\alpha}(\xi+\frac{d}{p}-\zeta_0)_+}.
		\end{align*}

		For the bilinear term, we use the Young inequality \eqref{besov-prop-y} putting the regularity $\xi$ onto the heat kernel:
		\begin{align*}
			&\int_0^v\|(\mathcal{B}^{\varepsilon}_{ \boldsymbol{\rho}^{\varepsilon}_{\mu}}(s,\cdot)\boldsymbol{\rho}^{\varepsilon}_{\mu}(s,\cdot))\ast\nabla p^\alpha_{v-s}\|_{B^{\xi}_{p',1}}ds\\
			&\leq C\int_0^v\|\mathcal{B}^{\varepsilon}_{\boldsymbol{\rho}^{\varepsilon}_{\mu}}(s,\cdot)\boldsymbol{\rho}^{\varepsilon}_{\mu}(s,\cdot)\|_{B^{0}_{1,\infty}}\|\nabla p^\alpha_{v-s}\|_{B^{\xi}_{p',1}}ds\\
			&\stackrel{\eqref{besov-prop-pr2}, \eqref{besov-prop-l}, \eqref{besov-prop-hk}}{\leq}C\int_0^v\|\mathcal{B}^{\varepsilon}_{\boldsymbol{\rho}^{\varepsilon}_{\mu}}(s,\cdot)\|_{B^{0}_{\infty,\infty}}\|\boldsymbol{\rho}^{\varepsilon}_{\mu}(s,\cdot)\|_{B^{0}_{1,\infty}}(v-s)^{-\frac{1}{\alpha}(1+\xi+\frac{d}{p})_+}ds\\
			&\stackrel{\eqref{ineq-B-rho-reg}, \eqref{besov-prop-e1}}{\leq}C_{\varepsilon}\int_0^v\|\boldsymbol{\rho}^{\varepsilon}_{\mu}(s,\cdot)\|_{L^1}(v-s)^{-\frac{1}{\alpha}(1+\xi+\frac{d}{p})_+}ds\\
			&\leq C_{\varepsilon}\int_0^v(v-s)^{-\frac{1}{\alpha}(1+\xi+\frac{d}{p})_+}ds.
		\end{align*}

		In order to use Lemma \ref{lemma-beta-function} we need that $\xi$ is such that 
		\begin{align}
			\label{xi-upper-bound}
			\xi<\alpha-1-\frac{d}{p}.
		\end{align}

		Combining the above estimates, we get
		\begin{align}
			\label{ineq-rhp-reg-xi}
			\|\brho^{\varepsilon}_{\mu}(v,\cdot)\|_{B^{\xi}_{p',1}}\leq C\|\mu\|_{B^{\bar{\beta}_0}_{\bar{p}_0,\bar{q}_0}}v^{-\frac{1}{\alpha}(\xi+\frac{d}{p}-\zeta_0)_+}+C_{\varepsilon}v^{1-\frac{1}{\alpha}(1+\xi+\frac{d}{p})_+}.
		\end{align}
		
		Introduce $\xi_0\geq0$ satisfying
		\begin{align}
			\label{xi0-upper-bound}
			\xi_0<\alpha-1.
		\end{align}

		We bootstrap the norm we bounded previously by regularity $\xi_0$ using the same techniques as before:
		\begin{align*}
			\|\brho^{\varepsilon}_{\mu}(v,\cdot)\|_{B^{\xi+\xi_0}_{p',1}}&\leq\|\mu\ast p^{\alpha}_{v}\|_{B^{\xi+\xi_0}_{p',1}}+C\int_0^v\|(\mathcal{B}^{\varepsilon}_{ \boldsymbol{\rho}^{\varepsilon}_{\mu}}(s,\cdot)\boldsymbol{\rho}^{\varepsilon}_{\mu}(s,\cdot))\ast\nabla p^\alpha_{v-s}\|_{B^{\xi+\xi_0}_{p',1}}ds\\
			&\stackrel{\eqref{prf-ineq-initial-data},\eqref{besov-prop-y}}{\leq} C\|\mu\|_{B^{\bar{\beta}_0}_{\bar{p}_0,\bar{q}_0}}v^{-\frac{1}{\alpha}(\xi+\xi_0+\frac{d}{p}-\zeta_0)_+}+C\int_0^v\|\mathcal{B}^{\varepsilon}_{\boldsymbol{\rho}^{\varepsilon}_{\mu}}(s,\cdot)\boldsymbol{\rho}^{\varepsilon}_{\mu}(s,\cdot)\|_{B^{\xi}_{p',\infty}}\|\nabla p^\alpha_{v-s}\|_{B^{\xi_0}_{1,1}}ds\\
			&\stackrel{\eqref{ineq-B-rho-reg},\eqref{besov-prop-pr2}, \eqref{besov-prop-hk}}{\leq} C\|\mu\|_{B^{\bar{\beta}_0}_{\bar{p}_0,\bar{q}_0}}v^{-\frac{1}{\alpha}(\xi+\xi_0+\frac{d}{p}-\zeta_0)_+}+C_{\varepsilon}\int_0^v\|\boldsymbol{\rho}^{\varepsilon}_{\mu}(s,\cdot)\|_{B^{\xi}_{p',1}}(v-s)^{-\frac{1}{\alpha}(1+\xi_0)}ds.
		\end{align*}

		Substituting the control \eqref{ineq-rhp-reg-xi} on $\|\boldsymbol{\rho}^{\varepsilon}_{\mu}(s,\cdot)\|_{B^{\xi}_{p',1}}$ in the above integral, we get
		\begin{align*}
			&\int_0^v\|\boldsymbol{\rho}^{\varepsilon}_{\mu}(s,\cdot)\|_{B^{\xi}_{p',1}}(v-s)^{-\frac{1}{\alpha}(1+\xi_0)}ds\\
			&\leq C\|\mu\|_{B^{\bar{\beta}_0}_{\bar{p}_0,\bar{q}_0}}\int_0^vs^{-\frac{1}{\alpha}(\xi+\frac{d}{p}-\zeta_0)_+}(v-s)^{-\frac{1}{\alpha}(1+\xi_0)}ds+C_{\varepsilon}\int_0^vs^{1-\frac{1}{\alpha}(1+\xi+\frac{d}{p})_+}(v-s)^{-\frac{1}{\alpha}(1+\xi_0)}ds\\
			&\leq C\|\mu\|_{B^{\bar{\beta}_0}_{\bar{p}_0,\bar{q}_0}}v^{1-\frac{1}{\alpha}(\xi+\frac{d}{p}-\zeta_0)_+-\frac{1}{\alpha}(1+\xi_0)}+C_{\varepsilon}v^{2-\frac{1}{\alpha}(1+\xi+\frac{d}{p})_+-\frac{1}{\alpha}(1+\xi_0)}.
		\end{align*}

		Above, we applied Lemma \ref{lemma-beta-function} since
		\begin{align*}
			\frac{1}{\alpha}\Big(\xi+\frac{d}{p}-\zeta_0\Big)<1\iff\xi<\alpha-\frac{d}{p}+\zeta_0
		\end{align*}
		and
		\begin{align*}
			-1+\frac{1}{\alpha}\Big(1+\xi+\frac{d}{p}\Big)<1\iff\xi<2\alpha-1-\frac{d}{p},
		\end{align*}
		which follows from \eqref{xi-upper-bound}, i.e., from $\xi<\alpha-1-\frac{d}{p}$. Therefore,
		\begin{align*}
			&\|\brho^{\varepsilon}_{\mu}(v,\cdot)\|_{B^{\xi+\xi_0}_{p',1}}\\
			&\leq C\|\mu\|_{B^{\bar{\beta}_0}_{\bar{p}_0,\bar{q}_0}}v^{-\frac{1}{\alpha}(\xi+\xi_0+\frac{d}{p}-\zeta_0)_+}+C_{\varepsilon}\|\mu\|_{B^{\bar{\beta}_0}_{\bar{p}_0,\bar{q}_0}}v^{1-\frac{1}{\alpha}(\xi+\frac{d}{p}-\zeta_0)_+-\frac{1}{\alpha}(1+\xi_0)}+C_{\varepsilon}^2v^{2-\frac{1}{\alpha}(1+\xi+\frac{d}{p})_+-\frac{1}{\alpha}(1+\xi_0)}\\
			&=C\|\mu\|_{B^{\bar{\beta}_0}_{\bar{p}_0,\bar{q}_0}}\sum_{k=0}^{1}C_{\varepsilon}^kv^{k-\frac{1}{\alpha}\big(\xi+(1-k)\xi_0+\frac{d}{p}-\zeta_0\big)_+-\frac{k}{\alpha}(1+\xi_0)}+C_{\varepsilon}^{2}v^{2-\frac{1}{\alpha}(1+\xi+\frac{d}{p})_+-\frac{1}{\alpha}(1+\xi_0)}\\
			&=\Bigg(C\|\mu\|_{B^{\bar{\beta}_0}_{\bar{p}_0,\bar{q}_0}}\sum_{k=0}^{j}C_{\varepsilon}^kv^{k-\frac{1}{\alpha}\big(\xi+(j-k)\xi_0+\frac{d}{p}-\zeta_0\big)_+-\frac{k}{\alpha}(1+\xi_0)}+C_{\varepsilon}^{j+1}v^{(j+1)-\frac{1}{\alpha}(1+\xi+\frac{d}{p})_+-\frac{j}{\alpha}(1+\xi_0)}\Bigg)\ind_{\{j=1\}}.
		\end{align*}

		The last identity is a compact form of the estimate. We denote by $j\in\N$ a step of the added regularity $\xi_0$, i.e., above, $\xi+j\xi_0=\xi+\xi_0$. The above control easily extends to a general $j\in\N$. We just must verify that all time singularities are integrable. First, let $n\in\N$ be such that $\xi+n\xi_0=-\beta-\theta$. By induction, suppose that we have the estimate at step $1\leq j-1\leq n-1$. Then, at step $2\leq j\leq n$, following the previous approach, we have
		\begin{align*}
			\|\brho^{\varepsilon}_{\mu}(v,\cdot)\|_{B^{\xi+j\xi_0}_{p',1}}&\leq\|\mu\ast p^{\alpha}_{v}\|_{B^{\xi+j\xi_0}_{p',1}}+C\int_0^v\|(\mathcal{B}^{\varepsilon}_{ \boldsymbol{\rho}^{\varepsilon}_{\mu}}(s,\cdot)\boldsymbol{\rho}^{\varepsilon}_{\mu}(s,\cdot))\ast\nabla p^\alpha_{v-s}\|_{B^{\xi+j\xi_0}_{p',1}}ds\\
			&\leq C\|\mu\|_{B^{\bar{\beta}_0}_{\bar{p}_0,\bar{q}_0}}v^{-\frac{1}{\alpha}(\xi+j\xi_0+\frac{d}{p}-\zeta_0)_+}+C\int_0^v\|\mathcal{B}^{\varepsilon}_{\boldsymbol{\rho}^{\varepsilon}_{\mu}}(s,\cdot)\boldsymbol{\rho}^{\varepsilon}_{\mu}(s,\cdot)\|_{B^{\xi+(j-1)\xi_0}_{p',\infty}}\|\nabla p^\alpha_{v-s}\|_{B^{\xi_0}_{1,1}}ds\\
			&\leq C\|\mu\|_{B^{\bar{\beta}_0}_{\bar{p}_0,\bar{q}_0}}v^{-\frac{1}{\alpha}(\xi+j\xi_0+\frac{d}{p}-\zeta_0)_+}+C_{\varepsilon}\int_0^v\|\boldsymbol{\rho}^{\varepsilon}_{\mu}(s,\cdot)\|_{B^{\xi+(j-1)\xi_0}_{p',1}}(v-s)^{-\frac{1}{\alpha}(1+\xi_0)}ds.
		\end{align*}

		Utilizing the estimate from the step $j-1$, we derive
		\begin{align*}
			&\int_0^v\|\boldsymbol{\rho}^{\varepsilon}_{\mu}(s,\cdot)\|_{B^{\xi+(j-1)\xi_0}_{p',1}}(v-s)^{-\frac{1}{\alpha}(1+\xi_0)}ds\\
			&\leq C\|\mu\|_{B^{\bar{\beta}_0}_{\bar{p}_0,\bar{q}_0}}\sum_{k=0}^{j-1}C_{\varepsilon}^k\int_0^vs^{k-\frac{1}{\alpha}\big(\xi+(j-1-k)\xi_0+\frac{d}{p}-\zeta_0\big)_+-\frac{k}{\alpha}(1+\xi_0)}(v-s)^{-\frac{1}{\alpha}(1+\xi_0)}ds\\
			&\quad+C_{\varepsilon}^{j}\int_0^vs^{j-\frac{1}{\alpha}(1+\xi+\frac{d}{p})_+-\frac{j-1}{\alpha}(1+\xi_0)}(v-s)^{-\frac{1}{\alpha}(1+\xi_0)}ds
		\end{align*}

		Let us check the integrability in the first integral. If $\xi+(j-1-k)\xi_0+\frac{d}{p}-\zeta_0\leq0$, then
		\begin{align*}
			k-\frac{1}{\alpha}\Big(\xi+(j-1-k)\xi_0+\frac{d}{p}-\zeta_0\Big)_+-\frac{k}{\alpha}(1+\xi_0)>-1\iff \xi_0<\Big(1+\frac{1}{k}\Big)\alpha-1,
		\end{align*}
		which is satisfied from \eqref{xi0-upper-bound}, i.e., from $\xi_0<\alpha-1$ for $\alpha>1$. If $\xi+(j-1-k)\xi_0+\frac{d}{p}-\zeta_0>0$, then
		\begin{align*}
			k-\frac{1}{\alpha}\Big(\xi+(j-1-k)\xi_0+\frac{d}{p}-\zeta_0\Big)_+-\frac{k}{\alpha}(1+\xi_0)>-1\\
			\iff\xi+(j-1)\xi_0<k(\alpha-1)+\Big(\alpha-\frac{d}{p}+\zeta_0\Big).
		\end{align*}

		Since the r.h.s. in the above inequality is increasing in $k$ for $\alpha>1$, it is enough to check the validity of this inequality for $k=0$:
		\begin{align*}
			\xi+(j-1)\xi_0\leq\xi+n\xi_0=-\beta-\theta<\alpha-\frac{d}{p}+\zeta_0.
		\end{align*}
		
		This inequality is true thanks to \eqref{C3}:
		\begin{align*}
			\beta>-\alpha+\frac{\alpha}{r}+\Big(-\theta-\beta+\frac{d}{p}-\zeta_0\Big)_+\geq-\alpha+\frac{\alpha}{r}-\theta-\beta+\frac{d}{p}-\zeta_0\\
			>-\alpha-\theta+\frac{d}{p}-\zeta_0,
		\end{align*}
		where we used the fact that $\frac{\alpha}{r}-\beta>0$ for the last inequality. Let us study the second integral. If $1+\xi+\frac{d}{p}\leq0$, then
		\begin{align*}
			j-\frac{1}{\alpha}\Big(1+\xi+\frac{d}{p}\Big)_+-\frac{j-1}{\alpha}(1+\xi_0)>-1\iff\xi_0<\frac{j+1}{j-1}\alpha-1,
		\end{align*}
		and this follows from \eqref{xi0-upper-bound}. If $1+\xi+\frac{d}{p}>0$, then
		\begin{align*}
			j-\frac{1}{\alpha}\Big(1+\xi+\frac{d}{p}\Big)_+-\frac{j-1}{\alpha}(1+\xi_0)>-1\iff\xi+(j-1)\xi_0<j(\alpha-1)+\alpha-\frac{d}{p},
		\end{align*}
		which follows from the choice of $\xi,\xi_0$, i.e., from \eqref{xi-upper-bound} and \eqref{xi0-upper-bound}. This means that Lemma \ref{lemma-beta-function} is applicable:
		\begin{align*}
			&\int_0^v\|\boldsymbol{\rho}^{\varepsilon}_{\mu}(s,\cdot)\|_{B^{\xi+(j-1)\xi_0}_{p',1}}(v-s)^{-\frac{1}{\alpha}(1+\xi_0)}ds\\
			&\leq C\|\mu\|_{B^{\bar{\beta}_0}_{\bar{p}_0,\bar{q}_0}}\sum_{k=0}^{j-1}C_{\varepsilon}^{k}v^{(k+1)-\frac{1}{\alpha}\big(\xi+(j-1-k)\xi_0+\frac{d}{p}-\zeta_0\big)_+-\frac{k+1}{\alpha}(1+\xi_0)}\\
			&\quad+C_{\varepsilon}^{j}v^{(j+1)-\frac{1}{\alpha}(1+\xi+\frac{d}{p})_+-\frac{j}{\alpha}(1+\xi_0)}\\
			&\leq C\|\mu\|_{B^{\bar{\beta}_0}_{\bar{p}_0,\bar{q}_0}}\sum_{k=1}^{j}C_{\varepsilon}^{k-1}v^{k-\frac{1}{\alpha}\big(\xi+(j-k)\xi_0+\frac{d}{p}-\zeta_0\big)_+-\frac{k}{\alpha}(1+\xi_0)}\\
			&\quad+C_{\varepsilon}^{j}v^{(j+1)-\frac{1}{\alpha}(1+\xi+\frac{d}{p})_+-\frac{j}{\alpha}(1+\xi_0)}.
		\end{align*}

		Therefore,
		\begin{align*}
			&\|\brho^{\varepsilon}_{\mu}(v,\cdot)\|_{B^{\xi+j\xi_0}_{p',1}}\leq C\|\mu\|_{B^{\bar{\beta}_0}_{\bar{p}_0,\bar{q}_0}}v^{-\frac{1}{\alpha}(\xi+j\xi_0+\frac{d}{p}-\zeta_0)_+}\\
			&\quad+C\|\mu\|_{B^{\bar{\beta}_0}_{\bar{p}_0,\bar{q}_0}}\sum_{k=1}^{j}C_{\varepsilon}^kv^{k-\frac{1}{\alpha}\big(\xi+(j-k)\xi_0+\frac{d}{p}-\zeta_0\big)_+-\frac{k}{\alpha}(1+\xi_0)}+C_{\varepsilon}^{j+1}v^{(j+1)-\frac{1}{\alpha}(1+\xi+\frac{d}{p})_+-\frac{j}{\alpha}(1+\xi_0)}\\
			&\leq C\|\mu\|_{B^{\bar{\beta}_0}_{\bar{p}_0,\bar{q}_0}}\sum_{k=0}^{j}C_{\varepsilon}^kv^{k-\frac{1}{\alpha}\big(\xi+(j-k)\xi_0+\frac{d}{p}-\zeta_0\big)_+-\frac{k}{\alpha}(1+\xi_0)}+C_{\varepsilon}^{j+1}v^{(j+1)-\frac{1}{\alpha}(1+\xi+\frac{d}{p})_+-\frac{j}{\alpha}(1+\xi_0)}.
		\end{align*}

		Finally, for $j=n$,
		\begin{align*}
			&\|\brho^{\varepsilon}_{\mu}(v,\cdot)\|_{B^{\xi+n\xi_0}_{p',1}}=\|\brho^{\varepsilon}_{\mu}(v,\cdot)\|_{B^{-\beta-\theta}_{p',1}}\\
			&\leq C\|\mu\|_{B^{\bar{\beta}_0}_{\bar{p}_0,\bar{q}_0}}\sum_{k=0}^{n}C_{\varepsilon}^kv^{k-\frac{1}{\alpha}\big(\xi+(n-k)\xi_0+\frac{d}{p}-\zeta_0\big)_+-\frac{k}{\alpha}(1+\xi_0)}+C_{\varepsilon}^{n+1}v^{(n+1)-\frac{1}{\alpha}(1+\xi+\frac{d}{p})_+-\frac{n}{\alpha}(1+\xi_0)}.
		\end{align*}

		Multiply the expression by the time weight:
		\begin{align*}
			&v^{\gamma}\|\brho^{\varepsilon}_{\mu}(v,\cdot)\|_{B^{-\beta-\theta}_{p',1}}\\
			&\leq C\|\mu\|_{B^{\bar{\beta}_0}_{\bar{p}_0,\bar{q}_0}}\sum_{k=0}^{n}C_{\varepsilon}^kv^{\gamma+k-\frac{1}{\alpha}\big(\xi+(n-k)\xi_0+\frac{d}{p}-\zeta_0\big)_+-\frac{k}{\alpha}(1+\xi_0)}+C_{\varepsilon}^{n+1}v^{\gamma+(n+1)-\frac{1}{\alpha}(1+\xi+\frac{d}{p})_+-\frac{n}{\alpha}(1+\xi_0)}.
		\end{align*}

		The last step is to check that all the exponents are positive so that in small time we do not have an explosion.

		\vspace{1mm}

		\textbullet\ Let us verify that $\gamma+k-\frac{1}{\alpha}\big(\xi+(n-k)\xi_0+\frac{d}{p}-\zeta_0\big)_+-\frac{k}{\alpha}(1+\xi_0)>0$. If $\xi+(n-k)\xi_0+\frac{d}{p}-\zeta_0\leq0$, then
		\begin{align*}
			\gamma+k-\frac{1}{\alpha}\Big(\xi+(n-k)\xi_0+\frac{d}{p}-\zeta_0\Big)_+-\frac{k}{\alpha}(1+\xi_0)>0\iff\xi_0<\alpha-1+\frac{\alpha\gamma}{k},
		\end{align*}
		following from \eqref{xi0-upper-bound} since $\gamma>0$. If $\xi+(n-k)\xi_0+\frac{d}{p}-\zeta_0>0$, then
		\begin{align*}
			\gamma+k-\frac{1}{\alpha}\Big(\xi+(n-k)\xi_0+\frac{d}{p}-\zeta_0\Big)_+-\frac{k}{\alpha}(1+\xi_0)>0\iff\xi+n\xi_0<k\alpha+\alpha\gamma-\frac{d}{p}+\zeta_0.
		\end{align*}
		It is enough to verify this for $k=0$ exploiting the definition of $\gamma$ and $n$:
		\begin{align*}
			\xi+n\xi_0=-\beta-\theta<\alpha\gamma-\frac{d}{p}+\zeta_0\iff\alpha\gamma>-\beta-\theta+\frac{d}{p}-\zeta_0.
		\end{align*}
		But by definition of $\gamma$ \eqref{def-gamma-weight}
		\begin{align*}
			\alpha\gamma=\eta+\Big(-\theta-\beta+\frac{d}{p}-\zeta_0\Big)_+>-\beta-\theta+\frac{d}{p}-\zeta_0.
		\end{align*}

		\vspace{1mm}
		
		\textbullet\ Now, let us check that $\gamma+(n+1)-\frac{1}{\alpha}(1+\xi+\frac{d}{p})_+-\frac{n}{\alpha}(1+\xi_0)>0$. If $1+\xi+\frac{d}{p}\leq0$, then
		\begin{align*}
			\gamma+(n+1)-\frac{1}{\alpha}\Big(1+\xi+\frac{d}{p}\Big)_+-\frac{n}{\alpha}(1+\xi_0)>0\iff\xi_0<\alpha-1+\alpha(\gamma+1),
		\end{align*}
		which follows from \eqref{xi0-upper-bound}. If $1+\xi+\frac{d}{p}>0$, then
		\begin{align*}
			\gamma+(n+1)-\frac{1}{\alpha}\Big(1+\xi+\frac{d}{p}\Big)_+-\frac{n}{\alpha}(1+\xi_0)>0\iff\xi+n\xi_0<\Big(\alpha-1-\frac{d}{p}\Big)+n(\alpha-1)+\alpha\gamma,
		\end{align*}
		which follows from \eqref{xi-upper-bound} and \eqref{xi0-upper-bound} for $\gamma>0$. And the result follows.
	\end{proof}
\end{lemma}

From the previous lemma, we are able to derive the following result.
\begin{lemma}[A priori control through a Grönwall-type inequality]
	\label{lemma-gronwall-control}
	Let conditions \eqref{C3} and \eqref{mu-range} be satisfied. Then, there exists a time $\mathcal{T}\in(0,T)$ and a constant $\tilde{C}>0$ independent of $\varepsilon$, such that for any $S\in(0,\mathcal{T})$,
	\begin{equation*}
		\sup_{s\in(0,S]} s ^\gamma\|\brho^\varepsilon_{\mu}(s,\cdot)\|_{B^{-\beta-\theta}_{p',1}}\leq \tilde{C}.
	\end{equation*}

	\begin{proof}
		Let $\mathcal{T}\in(0,T]$ be some time to be specified later and let $S\in(0,\mathcal{T})$. Recalling that by \eqref{def-f-I} $f^\varepsilon_{(0,S]}(s)=\sup_{v\in(0,s]} v^\gamma\|\brho^\varepsilon_{\mu}(v,\cdot)\|_{B^{-\beta-\theta}_{p',1}}$, rewrite \eqref{ineq-quadratic-density} as
		\begin{align}
			\label{prf-4-polynom}
			C_{b}(s)\big(f^\varepsilon_{(0,S]}(s)\big)^2-f^\varepsilon_{(0,S]}(s)+C_{\mu}(s)\geq0,
		\end{align}
		where
		\begin{align*}
			C_{\mu}(s) & =C\|\mu\|_{B^{\bar{\beta}_0}_{\bar{p}_0,\bar{q}_0}} s ^{\frac{\eta}{\alpha}}=:C_0 s ^{\frac{\eta}{\alpha}},                                                                             \\
			C_b(s)     & =C\Big(\|b\|_{L^r(B^\beta_{p,q})}+\|\Div(b)\|_{L^r(B^\beta_{p,q})}\Big) s ^{-\gamma+1-\frac{1}{r}+\frac{\beta-\delta}{\alpha}}=:C_b s ^{-\gamma+1-\frac{1}{r}+\frac{\beta-\delta}{\alpha}}.
		\end{align*}

		This polynomial admits two positive roots up to time $\mathcal{T}>S>t$
		\begin{align}
			\label{prf-4-roots}
			R[f^\varepsilon_{(0,S]}]_{+/-}(s)=\frac{1\pm\sqrt{1-4C_{\mu}(s)C_{b}(s)}}{2C_{b}(s)}
		\end{align}
		provided that
		\begin{align*}
			C_{\mu}(s)C_{b}(s)<\frac{1}{4}.
		\end{align*}

		The r.h.s. of \eqref{prf-4-roots} is continuous in $s$ and increasing (see computations below). Moreover, when $s\to t$, $R[f^\varepsilon_{(0,S]}]_{+}(s)\to 0^-$, $R[f^\varepsilon_{(0,S]}]_{-}(s)\to 0^+$. By Lemma \ref{lemma-norm-density-cont}, for fixed $\varepsilon>0$ $s\mapsto f^\varepsilon_{(0,S]}(s)$ is continuous and $f^\varepsilon_{(0,S]}(s)\xrightarrow[s\downarrow 0]{}0$ which implies that
		\begin{align*}
			C_{b}(s)\big(f^\varepsilon_{(0,S]}(s)\big)^2-f^\varepsilon_{(0,S]}(s)+C_{\mu}(s)\xrightarrow[s\downarrow 0]{}0.
		\end{align*}
		Then for $x_s\in\Big(\frac{1-\sqrt{1-4C_{\mu}(s)C_{b}(s)}}{2C_{b}(s)},\frac{1+\sqrt{1-4C_{\mu}(s)C_{b}(s)}}{2C_{b}(s)}\Big)$, we get that $C_{b}(s)x_s^2-x_s+C_{\mu}(s)<0$, which from the non-negativity of the polynomial \eqref{prf-4-polynom} implies that necessarily
		\begin{align}
			\label{prf-4-f-eps-bound}
			f^\varepsilon_{(0,S]}(s)\leq\frac{1-\sqrt{1-4C_{\mu}(s)C_{b}(s)}}{2C_{b}(s)},
		\end{align}
		which is the smallest positive root of the polynomial under the condition \eqref{prf-4-roots}. Note that
		\begin{align*}
			C_{\mu}(s)C_{b}(s)<\frac{1}{4}\iff s<(4C_0C_b)^{1/(\gamma-1+\frac{1}{r}-\frac{\beta-\delta+\eta}{\alpha})},
		\end{align*}
		where $\gamma-1+\frac{1}{r}-\frac{\beta-\delta+\eta}{\alpha}<0$ which follows from the choice of $\gamma$ given by \eqref{def-gamma-weight}. We then set $\mathcal{T}:=(4C_0C_b)^{1/(\gamma-1+\frac{1}{r}-\frac{\beta-\delta+\eta}{\alpha})}$. Using the fact that the r.h.s. of \eqref{prf-4-f-eps-bound} is increasing, we have
		\begin{align*}
			f^\varepsilon_{(0,S]}(s) & <\frac{1-\sqrt{1-4C_{\mu}(\mathcal{T})C_{b}(\mathcal{T})}}{2C_{b}(\mathcal{T})}\leq\frac{1}{2C_{b}(\mathcal{T})} \\
			                         & =\frac{1}{2C_b}(4C_0C_b)^{1-\frac{\eta}{\alpha(-\gamma+1-\frac{1}{r}+\frac{\beta-\delta+\eta}{\alpha})}}=:\tilde{C}.
		\end{align*}
		Finally, we conclude that
		\begin{align*}
			\sup_{s\in(0,S]}f^\varepsilon_{(0,S]}(s)=\sup_{s\in(0,S]} s ^\gamma\|\brho^\varepsilon_{\mu}(s,\cdot)\|_{B^{-\beta-\theta}_{p',1}}\leq\tilde{C}.
		\end{align*}

		Let us show that the r.h.s. of \eqref{prf-4-f-eps-bound} is indeed increasing. Denote for simplicity
		\begin{align*}
			z_s=\sqrt{1-4C_{\mu}(s)C_{b}(s)}\in(0,1].
		\end{align*} 
		Then
		\begin{align*}
			\big(R[f^\varepsilon_{(0,S]}]_{-}(s)\big)'=\frac{2\big(C_{\mu}'(s)C_{b}(s)+C_{\mu}(s)C_{b}'(s)\big)C_{b}(s)-z_s(1-z_s)C_{b}'(s)}{2C_b^2(s)z_s},
		\end{align*}
		and we only have to show that $z_s(1-z_s)C_{b}'(s)<2\big(C_{\mu}'(s)C_{b}(s)+C_{\mu}(s)C_{b}'(s)\big)C_{b}(s)$. Rewriting
		\begin{align*}
			1-z_s=\frac{(1-z_s)(1+z_s)}{1+z_s}=\frac{4C_{\mu}(s)C_{b}(s)}{1+z_s},
		\end{align*}
		and using
		\begin{align*}
			\frac{2z_s}{1+z_s}\leq1\quad\text{for}\quad z_s\leq1,
		\end{align*}
		we get
		\begin{align*}
			z_s(1-z_s)C_{b}'(s) & =2C_{\mu}(s)C_{b}(s)C_{b}'(s)\frac{2z_s}{1+z_s}\leq2C_{\mu}(s)C_{b}(s)C_{b}'(s) \\
			                    & <2\big(C_{\mu}(s)C_{b}'(s)+C_{\mu}'(s)C_{b}(s)\big)C_{b}(s),
		\end{align*}
		where the last inequality follows from $C_{\mu}(s),C_{b}(s),C_{\mu}'(s)>0$ on the interval $(0,S]$.
	\end{proof}
\end{lemma}

\begin{lemma}[Convergence of the mollified densities]
	\label{lemma-density-conv}
	Let conditions \eqref{C3} and \eqref{mu-range} be satisfied. Then, for all $S\in(0,\mathcal{T})$, with $\mathcal{T}$ defined in Lemma \ref{lemma-gronwall-control}, for any decreasing sequence $(\varepsilon_k)_{k\geq1}$ s.t. $\varepsilon_k\xrightarrow[k\to\infty]{}0$, the sequence $\big(\brho_{\mu}^{\varepsilon_k}\big)_{k\geq1}$ is \emph{Cauchy} in $L^\infty_\gamma((0,S],B^{-\beta-\theta}_{p',1})\cap L^\infty((0,S],L^1)$. In particular, there exists $\brho_{\mu}\in L^\infty_{\gamma}((0,S],B^{-\beta-\theta}_{p',1})$ s.t.
	\begin{align}
		\label{conv-mol-dens}
		\sup_{s\in(0,S]} s ^{\gamma}\|(\brho^{\varepsilon_k}_{\mu}-\brho_{\mu})(s,\cdot)\|_{B^{-\beta-\theta}_{p',1}}+\sup_{s\in(0,S]}\|(\brho^{\varepsilon_k}_{\mu}-\brho_{\mu})(s,\cdot)\|_{L^1}\xrightarrow[k\to+\infty]{}0.
	\end{align}

	\begin{proof}
		Fix $k,j\in\N$ s.t. $k\geq j$. For fixed $s\in(0,S]$, $y\in\R^d$, write from the Duhamel representation \eqref{duhamel-main-mollified}
		\begin{align*}
			\big(\brho^{\varepsilon_k}_{\mu}-\brho^{\varepsilon_j}_{\mu}\big)(s,y)=-\int_0^s\Big(\mathcal{B}^{\varepsilon_k}_{\brho^{\varepsilon_k}_{\mu}}(v,\cdot)\brho^{\varepsilon_k}_{\mu}(v,\cdot)-\mathcal{B}^{\varepsilon_j}_{\brho^{\varepsilon_j}_{\mu}}(v,\cdot)\brho^{\varepsilon_j}_{\mu}(v,\cdot)\Big)\ast\nabla p^\alpha_{s-v}(y)dv.
		\end{align*}

		\textbullet\ We first show that $\big(\brho_{\mu}^{\varepsilon_k}\big)_{k\geq1}$ is a Cauchy sequence in $L^\infty_\gamma((0,S],B^{-\beta-\theta}_{p',1})$. From the Young inequality \eqref{besov-prop-y} and integration by parts for $\delta\in[0,1)$ small enough, we have
		\begin{align}
			\label{prf-conv-ineq-1}
			\begin{split}
				 & \Big\|\Big(\mathcal{B}^{\varepsilon_k}_{\brho^{\varepsilon_k}_{\mu}}(v,\cdot)\brho^{\varepsilon_k}_{\mu}(v,\cdot)-\mathcal{B}^{\varepsilon_j}_{\brho^{\varepsilon_j}_{\mu}}(v,\cdot)\brho^{\varepsilon_j}_{\mu}(v,\cdot)\Big)\ast\nabla p^\alpha_{s-v}\Big\|_{B^{-\beta-\theta}_{p',1}}                                                   \\
				 &\leq C\Big(\|\Div(\mathcal{B}^{\varepsilon_k}_{\brho^{\varepsilon_k}_{\mu}}(v,\cdot))\brho^{\varepsilon_k}_{\mu}(v,\cdot)-\Div(\mathcal{B}^{\varepsilon_j}_{\brho^{\varepsilon_j}_{\mu}}(v,\cdot))\brho^{\varepsilon_j}_{\mu}(v,\cdot)\|_{B^{-\theta-\delta}_{p',\infty}}                                                    \\
				 & \quad+\|\mathcal{B}^{\varepsilon_k}_{\brho^{\varepsilon_k}_{\mu}}(v,\cdot)\cdot\nabla\brho^{\varepsilon_k}_{\mu}(v,\cdot)-\mathcal{B}^{\varepsilon_j}_{\brho^{\varepsilon_j}_{\mu}}(v,\cdot)\cdot\nabla\brho^{\varepsilon_j}_{\mu}(v,\cdot)\|_{B^{-\theta-\delta}_{p',\infty}}\Big)\| p^\alpha_{s-v}\|_{B^{-\beta+\delta}_{1,1}}.
			\end{split}
		\end{align}

		Let us handle the second term in the sum of \eqref{prf-conv-ineq-1}. By adding and subtracting $\mathcal{B}^{\varepsilon_k}_{\brho^{\varepsilon_k}_{\mu}}(v,\cdot)\cdot\nabla\brho^{\varepsilon_j}_{\mu}(v,\cdot)$, we get
		\begin{align*}
			 & \|\mathcal{B}^{\varepsilon_k}_{\brho^{\varepsilon_k}_{\mu}}(v,\cdot)\cdot\nabla\brho^{\varepsilon_k}_{\mu}(v,\cdot)-\mathcal{B}^{\varepsilon_j}_{\brho^{\varepsilon_j}_{\mu}}(v,\cdot)\cdot\nabla\brho^{\varepsilon_j}_{\mu}(v,\cdot)\|_{B^{-\theta-\delta}_{p',\infty}}                                                                                                                                                             \\
			 & \quad\leq\|\mathcal{B}^{\varepsilon_k}_{\brho^{\varepsilon_k}_{\mu}}(v,\cdot)\cdot\nabla\big(\brho^{\varepsilon_k}_{\mu}-\brho^{\varepsilon_j}_{\mu}\big)(v,\cdot)\|_{B^{-\theta-\delta}_{p',\infty}}+\|\big(\mathcal{B}^{\varepsilon_k}_{\brho^{\varepsilon_k}_{\mu}}-\mathcal{B}^{\varepsilon_j}_{\brho^{\varepsilon_j}_{\mu}}\big)(v,\cdot)\cdot\nabla\brho^{\varepsilon_j}_{\mu}(v,\cdot)\|_{B^{-\theta-\delta}_{p',\infty}}.
		\end{align*}

		Now, we want to apply the product rule given by \eqref{ineq-product-rule} on both terms with the following parameters:
		\begin{align*}
			\lambda=\theta+\delta,\quad\lambda_1=\theta+\delta',\quad\lambda_2=\theta+\delta,\quad \ell=p',
		\end{align*}
		where $\delta'\in[0,1)$. The margin $\delta'$ is needed in order to profit from the convergence of the mollified drift given by Proposition \ref{prop-b-conv}. Indeed, the convergence holds only for any $\bar{\beta}<\beta$, and this is exactly guaranteed by $\delta'$. On the other hand, $\delta$ is needed to ensure that the product of two terms has at most the regularity of the less regular term. Taking $\delta'<\delta$, we see that
		\begin{align*}
			\lambda_2>\lambda_1,\quad \lambda>\lambda_1.
		\end{align*}

		Therefore, we have
		\begin{align}
			\label{prf-conv-pr2}
			\|\big(\mathcal{B}^{\varepsilon_k}_{\brho^{\varepsilon_k}_{\mu}}-\mathcal{B}^{\varepsilon_j}_{\brho^{\varepsilon_j}_{\mu}}\big)(v,\cdot)\cdot\nabla\brho^{\varepsilon_j}_{\mu}(v,\cdot)\|_{B^{-\theta-\delta}_{p',\infty}}\leq C\|\big(\mathcal{B}^{\varepsilon_k}_{\brho^{\varepsilon_k}_{\mu}}-\mathcal{B}^{\varepsilon_j}_{\brho^{\varepsilon_j}_{\mu}}\big)(v,\cdot)\|_{B^{-\theta-\delta'}_{\infty,\infty}}\|\nabla\brho^{\varepsilon_j}_{\mu}(v,\cdot)\|_{B^{\theta+\delta}_{p',1}}
		\end{align}
		and for $\delta'=0$,
		\begin{align}
			\label{prf-conv-pr1}
			\|\mathcal{B}^{\varepsilon_k}_{\brho^{\varepsilon_k}_{\mu}}(v,\cdot)\cdot\nabla\big(\brho^{\varepsilon_k}_{\mu}-\brho^{\varepsilon_j}_{\mu}\big)(v,\cdot)\|_{B^{-\theta-\delta}_{p',\infty}}\leq C\|\mathcal{B}^{\varepsilon_k}_{\brho^{\varepsilon_k}_{\mu}}(v,\cdot)\|_{B^{-\theta}_{\infty,\infty}}\|\nabla\big(\brho^{\varepsilon_k}_{\mu}-\brho^{\varepsilon_j}_{\mu}\big)(v,\cdot)\|_{B^{\theta+\delta}_{p',1}}.
		\end{align}

		Now, in \eqref{prf-conv-pr1} we simply apply the Young inequality \eqref{besov-prop-y}, the lifting property \eqref{besov-prop-l} and the embedding \eqref{besov-prop-e1}:
		\begin{align*}
			\|\mathcal{B}^{\varepsilon_k}_{\brho^{\varepsilon_k}_{\mu}}(v,\cdot)\cdot\nabla\big(\brho^{\varepsilon_k}_{\mu}-\brho^{\varepsilon_j}_{\mu}\big)(v,\cdot)\|_{B^{-\theta-\delta}_{p',\infty}} & \leq C\|b^{\varepsilon_k}(v,\cdot)\|_{B^\beta_{p,q}}\|\brho^{\varepsilon_k}_{\mu}(v,\cdot)\|_{B^{-\beta-\theta}_{p',1}}\|\big(\brho^{\varepsilon_k}_{\mu}-\brho^{\varepsilon_j}_{\mu}\big)(v,\cdot)\|_{B^{1+\theta+\delta}_{p',1}} \\
			                                                                                                                                                                                                   & \leq C\|b^{\varepsilon_k}(v,\cdot)\|_{B^\beta_{p,q}}\|\brho^{\varepsilon_k}_{\mu}(v,\cdot)\|_{B^{-\beta-\theta}_{p',1}}\|\big(\brho^{\varepsilon_k}_{\mu}-\brho^{\varepsilon_j}_{\mu}\big)(v,\cdot)\|_{B^{-\beta-\theta}_{p',1}},
		\end{align*}
		noting that $1+\theta+\delta\leq-\beta-\theta$, i.e., $\beta\leq-1-2\theta-\delta$ from \eqref{C3}. Next, for \eqref{prf-conv-pr2}, using again the Young inequality \eqref{besov-prop-y}, the lifting property \eqref{besov-prop-l} and the embedding \eqref{besov-prop-e1},
		\begin{align*}
			 & \|\big(\mathcal{B}^{\varepsilon_k}_{\brho^{\varepsilon_k}_{\mu}}-\mathcal{B}^{\varepsilon_j}_{\brho^{\varepsilon_j}_{\mu}}\big)(v,\cdot)\cdot\nabla\brho^{\varepsilon_j}_{\mu}(v,\cdot)\|_{B^{-\theta-\delta}_{p',\infty}}                                                                                                                                                                                                                                       \\
			 & \leq C\Big(\|\big(\mathcal{B}^{\varepsilon_k}_{\brho^{\varepsilon_k}_{\mu}}-\mathcal{B}^{\varepsilon_k}_{\brho^{\varepsilon_j}_{\mu}}\big)(v,\cdot)\|_{B^{-\theta-\delta'}_{\infty,\infty}}+\|\big(\mathcal{B}^{\varepsilon_k}_{\brho^{\varepsilon_j}_{\mu}}-\mathcal{B}^{\varepsilon_j}_{\brho^{\varepsilon_j}_{\mu}}\big)(v,\cdot)\|_{B^{-\theta-\delta'}_{\infty,\infty}}\Big)\|\nabla\brho^{\varepsilon_j}_{\mu}(v,\cdot)\|_{B^{\theta+\delta}_{p',1}} \\
			 & \leq C\Big(\|b^{\varepsilon_k}(v,\cdot)\|_{B^{\beta-\delta'}_{p,q}}\|\big(\brho^{\varepsilon_k}_{\mu}-\brho^{\varepsilon_j}_{\mu}\big)(v,\cdot)\|_{B^{-\beta-\theta}_{p',1}}+\|\big(b^{\varepsilon_k}-b^{\varepsilon_j}\big)(v,\cdot)\|_{B^{\beta-\delta'}_{p,q}}\|\brho^{\varepsilon_j}_{\mu}(v,\cdot)\|_{B^{-\beta-\theta}_{p',1}}\Big)\\
			 &\quad\cdot\|\brho^{\varepsilon_j}_{\mu}(v,\cdot)\|_{B^{1+\theta+\delta}_{p',1}}                                               \\
			 & \leq C\Big(\|b^{\varepsilon_k}(v,\cdot)\|_{B^{\beta}_{p,q}}\|\big(\brho^{\varepsilon_k}_{\mu}-\brho^{\varepsilon_j}_{\mu}\big)(v,\cdot)\|_{B^{-\beta-\theta}_{p',1}}+\|\big(b^{\varepsilon_k}-b^{\varepsilon_j}\big)(v,\cdot)\|_{B^{\beta-\delta'}_{p,q}}\|\brho^{\varepsilon_j}_{\mu}(v,\cdot)\|_{B^{-\beta-\theta}_{p',1}}\Big)\\
			 &\quad\cdot\|\brho^{\varepsilon_j}_{\mu}(v,\cdot)\|_{B^{-\beta-\theta}_{p',1}}.
		\end{align*}
		Therefore, we obtain that the second term in \eqref{prf-conv-ineq-1} is bounded as
		\begin{align*}
			 & \|\mathcal{B}^{\varepsilon_k}_{\brho^{\varepsilon_k}_{\mu}}(v,\cdot)\cdot\nabla\brho^{\varepsilon_k}_{\mu}(v,\cdot)-\mathcal{B}^{\varepsilon_j}_{\brho^{\varepsilon_j}_{\mu}}(v,\cdot)\cdot\nabla\brho^{\varepsilon_j}_{\mu}(v,\cdot)\|_{B^{-\theta-\delta}_{p',\infty}}                                                                                                                                     \\
			 & \leq C\Big(\|b^{\varepsilon_k}(v,\cdot)\|_{B^{\beta}_{p,q}}\|\big(\brho^{\varepsilon_k}_{\mu}-\brho^{\varepsilon_j}_{\mu}\big)(v,\cdot)\|_{B^{-\beta-\theta}_{p',1}}\|\brho^{\varepsilon_j}_{\mu}(v,\cdot)\|_{B^{-\beta-\theta}_{p',1}}\\
			 &\quad+\|\big(b^{\varepsilon_k}-b^{\varepsilon_j}\big)(v,\cdot)\|_{B^{\beta-\delta'}_{p,q}}\|\brho^{\varepsilon_j}_{\mu}(v,\cdot)\|^2_{B^{-\beta-\theta}_{p',1}}\Big).
		\end{align*}
		Now, performing a similar analysis on the first term in \eqref{prf-conv-ineq-1}, as a result, we obtain
		\begin{align*}
			 & \Big\|\Div(\mathcal{B}^{\varepsilon_k}_{\brho^{\varepsilon_k}_{\mu}}(v,\cdot))\brho^{\varepsilon_k}_{\mu}(v,\cdot)-\Div(\mathcal{B}^{\varepsilon_j}_{\brho^{\varepsilon_j}_{\mu}}(v,\cdot))\brho^{\varepsilon_j}_{\mu}(v,\cdot)\Big\|_{B^{-\theta-\delta}_{p',\infty}}                                                                                                                                                    \\
			 & \leq C\Big(\|\Div(b^{\varepsilon_k}(v,\cdot))\|_{B^{\beta}_{p,q}}\|\big(\brho^{\varepsilon_k}_{\mu}-\brho^{\varepsilon_j}_{\mu}\big)(v,\cdot)\|_{B^{-\beta-\theta}_{p',1}}\|\brho^{\varepsilon_j}_{\mu}(v,\cdot)\|_{B^{-\beta-\theta}_{p',1}}\\
			 &\quad+\|\Div\big(b^{\varepsilon_k}-b^{\varepsilon_j}\big)(v,\cdot)\|_{B^{\beta-\delta'}_{p,q}}\|\brho^{\varepsilon_j}_{\mu}(v,\cdot)\|^2_{B^{-\beta-\theta}_{p',1}}\Big).
		\end{align*}

		Substituting the previous two bounds into \eqref{prf-conv-ineq-1}, we get
		\begin{align*}
			 & \Big\|\Big(\mathcal{B}^{\varepsilon_k}_{\brho^{\varepsilon_k}_{\mu}}(v,\cdot)\brho^{\varepsilon_k}_{\mu}(v,\cdot)-\mathcal{B}^{\varepsilon_j}_{\brho^{\varepsilon_j}_{\mu}}(v,\cdot)\brho^{\varepsilon_j}_{\mu}(v,\cdot)\Big)\ast\nabla p^\alpha_{s-v}\Big\|_{B^{-\beta-\theta}_{p',1}}                               \\
			 & \leq C\Bigg(\Big(\|b^{\varepsilon_k}(v,\cdot)\|_{B^{\beta}_{p,q}}+\|\Div(b^{\varepsilon_k}(v,\cdot))\|_{B^{\beta}_{p,q}}\Big)\|\big(\brho^{\varepsilon_k}_{\mu}-\brho^{\varepsilon_j}_{\mu}\big)(v,\cdot)\|_{B^{-\beta-\theta}_{p',1}}\|\brho^{\varepsilon_j}_{\mu}(v,\cdot)\|_{B^{-\beta-\theta}_{p',1}}   \\
			 & \quad+\Big(\|\big(b^{\varepsilon_k}-b^{\varepsilon_j}\big)(v,\cdot)\|_{B^{\beta-\delta'}_{p,q}}+\|\Div\big(b^{\varepsilon_k}-b^{\varepsilon_j}\big)(v,\cdot)\|_{B^{\beta-\delta'}_{p,q}}\Big)\|\brho^{\varepsilon_j}_{\mu}(v,\cdot)\|^2_{B^{-\beta-\theta}_{p',1}}\Bigg)\|p^\alpha_{s-v}\|_{B^{-\beta+\delta}_{1,1}} \\
			 & \leq C\Bigg(\Big(\|b^{\varepsilon_k}(v,\cdot)\|_{B^{\beta}_{p,q}}+\|\Div(b^{\varepsilon_k}(v,\cdot))\|_{B^{\beta}_{p,q}}\Big)\|\big(\brho^{\varepsilon_k}_{\mu}-\brho^{\varepsilon_j}_{\mu}\big)(v,\cdot)\|_{B^{-\beta-\theta}_{p',1}} v^{-\gamma}                                                         \\
			 & \quad+\Big(\|\big(b^{\varepsilon_k}-b^{\varepsilon_j}\big)(v,\cdot)\|_{B^{\beta-\delta'}_{p,q}}+\|\Div\big(b^{\varepsilon_k}-b^{\varepsilon_j}\big)(v,\cdot)\|_{B^{\beta-\delta'}_{p,q}}\Big) v^{-2\gamma}\Bigg)(s-v)^{\frac{\beta-\delta}{\alpha}},
		\end{align*}
		where the last inequality comes from Lemma \ref{lemma-density-estimate} and the heat kernel estimate \eqref{besov-prop-hk}. Hence, with $\bar{r}<+\infty$ defined as $\bar{r}=r\ind_{r<+\infty}+\tilde{r}\ind_{r=+\infty}$ for any finite positive $\tilde{r}$,
		\begin{align*}
			 & \|\big(\brho^{\varepsilon_k}_{\mu}-\brho^{\varepsilon_j}_{\mu}\big)(s,\cdot)\|_{B^{-\beta-\theta}_{p',1}}                                                                                                                                                                                                                           \\
			 & \leq C \int_0^s\Big\|\Big(\mathcal{B}^{\varepsilon_k}_{\brho^{\varepsilon_k}_{\mu}}(v,\cdot)\brho^{\varepsilon_k}_{\mu}(v,\cdot)-\mathcal{B}^{\varepsilon_j}_{\brho^{\varepsilon_j}_{\mu}}(v,\cdot)\brho^{\varepsilon_j}_{\mu}(v,\cdot)\Big)\ast\nabla p^\alpha_{s-v}\Big\|_{B^{-\beta-\theta}_{p',1}}dv                                \\
			 & \leq C\int_0^s\Bigg(\Big(\|b^{\varepsilon_k}(v,\cdot)\|_{B^{\beta}_{p,q}}+\|\Div(b^{\varepsilon_k}(v,\cdot))\|_{B^{\beta}_{p,q}}\Big)\|\big(\brho^{\varepsilon_k}_{\mu}-\brho^{\varepsilon_j}_{\mu}\big)(v,\cdot)\|_{B^{-\beta-\theta}_{p',1}} v^{-\gamma}                                                                        \\
			 & \quad+\Big(\|\big(b^{\varepsilon_k}-b^{\varepsilon_j}\big)(v,\cdot)\|_{B^{\beta-\delta'}_{p,q}}+\|\Div\big(b^{\varepsilon_k}-b^{\varepsilon_j}\big)(v,\cdot)\|_{B^{\beta-\delta'}_{p,q}}\Big) v^{-2\gamma}\Bigg)(s-v)^{\frac{\beta-\delta}{\alpha}}dv.
		\end{align*}

		Applying the Hölder inequality in time yields
		\begin{align*}
			& \|\big(\brho^{\varepsilon_k}_{\mu}-\brho^{\varepsilon_j}_{\mu}\big)(s,\cdot)\|_{B^{-\beta-\theta}_{p',1}}                                                                                                                                                                                                                           \\
			&\leq C\Big(\|b^{\varepsilon_k}\|_{L^{r}(B^{\beta}_{p,q})}+\|\Div(b^{\varepsilon_k})\|_{L^{r}(B^{\beta}_{p,q})}\Big)\Bigg(\int_0^s\|\big(\brho^{\varepsilon_k}_{\mu}-\brho^{\varepsilon_j}_{\mu}\big)(v,\cdot)\|_{B^{-\beta-\theta}_{p',1}}^{r'} v^{-r'\gamma}(s-v)^{r'\frac{\beta-\delta}{\alpha}}dv\Bigg)^{\frac{1}{r'}}        \\
			 & \quad+C\Big(\|(b^{\varepsilon_k}-b^{\varepsilon_j})\|_{L^{\bar{r}}(B^{\beta-\delta'}_{p,q})}+\|\Div(b^{\varepsilon_k}-b^{\varepsilon_j})\|_{L^{\bar{r}}(B^{\beta-\delta'}_{p,q})}\Big)\Bigg(\int_0^s v^{-2\bar{r}'\gamma}(s-v)^{\bar{r}'\frac{\beta-\delta}{\alpha}}dv\Bigg)^{\frac{1}{\bar{r}'}}                                    \\
			 & \leq C\Big(\|b\|_{L^{r}(B^{\beta}_{p,q})}+\|\Div(b)\|_{L^{r}(B^{\beta}_{p,q})}\Big)\Bigg(\int_0^s\Big(\sup_{r\in(0,v]} r^\gamma\|\big(\brho^{\varepsilon_k}_{\mu}-\brho^{\varepsilon_j}_{\mu}\big)(r,\cdot)\|_{B^{-\beta-\theta}_{p',1}}\Big)^{r'} v^{-2r'\gamma}\\
			 &\quad\cdot(s-v)^{r'\frac{\beta-\delta}{\alpha}}dv\Bigg)^{\frac{1}{r'}}+C\Big(\|(b^{\varepsilon_k}-b^{\varepsilon_j})\|_{L^{\bar{r}}(B^{\beta-\delta'}_{p,q})}+\|\Div(b^{\varepsilon_k}-b^{\varepsilon_j})\|_{L^{\bar{r}}(B^{\beta-\delta'}_{p,q})}\Big) s ^{\frac{1}{\bar{r}'}-2\gamma+\frac{\beta-\delta}{\alpha}},
		\end{align*}
		where the last inequality follows from Lemma \ref{lemma-beta-function}. Multiplying the above inequality by the time factor $s^\gamma$ and taking the supremum over $s\in(0,S]$, we have
		\begin{align*}
			& \sup_{s\in(0,S]} s^\gamma\|\big(\brho^{\varepsilon_k}_{\mu}-\brho^{\varepsilon_j}_{\mu}\big)(s,\cdot)\|_{B^{-\beta-\theta}_{p',1}}\\
			& \leq C\Big(\|(b^{\varepsilon_k}-b^{\varepsilon_j})\|_{L^{\bar{r}}(B^{\beta-\delta'}_{p,q})}+\|\Div(b^{\varepsilon_k}-b^{\varepsilon_j})\|_{L^{\bar{r}}(B^{\beta-\delta'}_{p,q})}\Big)S^{\frac{1}{\bar{r}'}-\gamma+\frac{\beta-\delta}{\alpha}}\\
			&\quad+C\Big(\|b\|_{L^{r}(B^{\beta}_{p,q})}+\|\Div(b)\|_{L^{r}(B^{\beta}_{p,q})}\Big)S^\gamma\\
			&\quad\times\Bigg(\int_0^S\Big(\sup_{s\in(0,v]} s^\gamma\|\big(\brho^{\varepsilon_k}_{\mu}-\brho^{\varepsilon_j}_{\mu}\big)(s,\cdot)\|_{B^{-\beta-\theta}_{p',1}}\Big)^{r'} v^{-2r'\gamma}(S-v)^{r'\frac{\beta-\delta}{\alpha}}dv\Bigg)^{\frac{1}{r'}}\\
			& \leq C\Big(\|(b^{\varepsilon_k}-b^{\varepsilon_j})\|_{L^{\bar{r}}(B^{\beta-\delta'}_{p,q})}+\|\Div(b^{\varepsilon_k}-b^{\varepsilon_j})\|_{L^{\bar{r}}(B^{\beta-\delta'}_{p,q})}\Big)S^{\frac{1}{\bar{r}'}-\gamma+\frac{\beta-\delta}{\alpha}}\\
			&\quad+C\Big(\|b\|_{L^{r}(B^{\beta}_{p,q})}+\|\Div(b)\|_{L^{r}(B^{\beta}_{p,q})}\Big)\sup_{s\in(0,S]}s^\gamma\|\big(\brho^{\varepsilon_k}_{\mu}-\brho^{\varepsilon_j}_{\mu}\big)(s,\cdot)\|_{B^{-\beta-\theta}_{p',1}}S^{\frac{1}{r'}-\gamma+\frac{\beta-\delta}{\alpha}}.
		\end{align*}

		Consequently, from the circularizing argument, we obtain
		\begin{align*}
			&\sup_{s\in(0,S]} s^\gamma\|\big(\brho^{\varepsilon_k}_{\mu}-\brho^{\varepsilon_j}_{\mu}\big)(s,\cdot)\|_{B^{-\beta-\theta}_{p',1}}\\
			&\leq \frac{C\Big(\|b^{\varepsilon_k}-b^{\varepsilon_j}\|_{L^{\bar{r}}(B^{\beta-\delta'}_{p,q})}+\|\Div(b^{\varepsilon_k}-b^{\varepsilon_j})\|_{L^{\bar{r}}(B^{\beta-\delta'}_{p,q})}\Big)S^{\frac{1}{\bar{r}'}-\gamma+\frac{\beta-\delta}{\alpha}}}{1-C\Big(\|b\|_{L^{r}(B^{\beta}_{p,q})}+\|\Div(b)\|_{L^{r}(B^{\beta}_{p,q})}\Big)S^{\frac{1}{\bar{r}'}-\gamma+\frac{\beta-\delta}{\alpha}}}.
		\end{align*}
		Since $\Big(\|b^{\varepsilon_k}-b^{\varepsilon_j}\|_{L^{\bar{r}}(B^{\beta-\delta'}_{p,q})}+\|\Div(b^{\varepsilon_k}-b^{\varepsilon_j})\|_{L^{\bar{r}}(B^{\beta-\delta'}_{p,q})}\Big)\xrightarrow[j,k\to+\infty]{}0$ by Proposition \ref{prop-b-conv}, $\frac{1}{\bar{r}'}-\gamma+\frac{\beta-\delta}{\alpha}>0$ and $S\leq1$, we deduce that $\big(\brho_{\mu}^{\varepsilon_k}\big)_{k\geq1}$ is a Cauchy sequence in $L^\infty_\gamma((0,S],B^{-\beta-\theta}_{p',1})$.

		\vspace{3mm}

		\textbullet\ Now, let us prove that $\big(\brho_{\mu}^{\varepsilon_k}\big)_{k\geq1}$ is a Cauchy sequence in $L^\infty((0,S],L^1)$. 

		\vspace{1mm}

	 	Thanks to the embedding $B^{\theta+\delta}_{1,1}\hookrightarrow B^0_{1,1}\hookrightarrow L^1$ given by \eqref{besov-prop-e1} for $\theta$ defined as before and $\delta>0$, we can work with the norm in $B^{\theta+\delta}_{1,1}$ to derive the desired estimate. We refer to Remark \ref{remark-mu-reg} below for the detailed explanation about why we prove the estimate in exactly this space and why the additional regularity assumption on the initial measure appears in the calculations below. Repeating the arguments from the first part of the proof for \eqref{prf-conv-ineq-1} with the integrability index $\ell=1$ for the product rule and applying first the Young inequality \eqref{besov-prop-y} and the heat kernel estimate \eqref{besov-prop-hk}, we get
		\begin{align}
			\label{prf-l1-quad}
			\begin{split}
			 & \Big\|\Big(\mathcal{B}^{\varepsilon_k}_{\brho^{\varepsilon_k}_{\mu}}(v,\cdot)\brho^{\varepsilon_k}_{\mu}(v,\cdot)-\mathcal{B}^{\varepsilon_j}_{\brho^{\varepsilon_j}_{\mu}}(v,\cdot)\brho^{\varepsilon_j}_{\mu}(v,\cdot)\Big)\ast\nabla p^\alpha_{s-v}\Big\|_{B^{\theta+\delta}_{1,1}}                                                              \\
			 & \leq C\|\mathcal{B}^{\varepsilon_k}_{\brho^{\varepsilon_k}_{\mu}}(v,\cdot)\brho^{\varepsilon_k}_{\mu}(v,\cdot)-\mathcal{B}^{\varepsilon_j}_{\brho^{\varepsilon_j}_{\mu}}(v,\cdot)\brho^{\varepsilon_j}_{\mu}(v,\cdot)\|_{B^{-\theta-\delta}_{1,1}}\|\nabla p^\alpha_{s-v}\|_{B^{2\theta+2\delta}_{1,1}} \\
			 & \leq C\Bigg(\|b^{\varepsilon_k}(v,\cdot)\|_{B^{\beta}_{p,q}}\Big(\|\big(\brho^{\varepsilon_k}_{\mu}-\brho^{\varepsilon_j}_{\mu}\big)(v,\cdot)\|_{B^{\theta+\delta}_{1,1}} v^{-\gamma}+\|\big(\brho^{\varepsilon_k}_{\mu}-\brho^{\varepsilon_j}_{\mu}\big)(v,\cdot)\|_{B^{-\beta-\theta}_{p',1}}\|\brho^{\varepsilon_j}_{\mu}(v,\cdot)\|_{B^{\theta+\delta}_{1,1}}\Big)\\                                                                               
			 & \quad+\|\big(b^{\varepsilon_k}-b^{\varepsilon_j}\big)(v,\cdot)\|_{B^{\beta-\delta'}_{p,q}}\|\brho^{\varepsilon_j}_{\mu}(v,\cdot)\|_{B^{\theta+\delta}_{1,1}} v^{-\gamma}\Bigg)(s-v)^{-\frac{1+2\theta+2\delta}{\alpha}}.
			\end{split}
		\end{align}

		Here, we used Lemma \ref{lemma-density-estimate} to bound $\|\brho^{\varepsilon_j}_{\mu}(v,\cdot)\|_{B^{-\beta-\theta}_{p',1}}$ by $v^{-\gamma}$. Thus,
		\begin{align*}
			& \|\big(\brho^{\varepsilon_k}_{\mu}-\brho^{\varepsilon_j}_{\mu}\big)(s,\cdot)\|_{B^{\theta+\delta}_{1,1}}                                                                                                                                                                                                        \\
			 & \leq C\int_0^s\Bigg(\|b^{\varepsilon_k}(v,\cdot)\|_{B^{\beta}_{p,q}}\Big(\|\big(\brho^{\varepsilon_k}_{\mu}-\brho^{\varepsilon_j}_{\mu}\big)(v,\cdot)\|_{B^{\theta+\delta}_{1,1}}+ v^{\gamma}\|\big(\brho^{\varepsilon_k}_{\mu}-\brho^{\varepsilon_j}_{\mu}\big)(v,\cdot)\|_{B^{-\beta-\theta}_{p',1}}\|\brho^{\varepsilon_j}_{\mu}(v,\cdot)\|_{B^{\theta+\delta}_{1,1}}\Big)\\                                                                               
			 & \quad+\|\big(b^{\varepsilon_k}-b^{\varepsilon_j}\big)(v,\cdot)\|_{B^{\beta-\delta'}_{p,q}}\|\brho^{\varepsilon_j}_{\mu}(v,\cdot)\|_{B^{\theta+\delta}_{1,1}}\Bigg) v^{-\gamma}(s-v)^{-\frac{1+2\theta+2\delta}{\alpha}}dv,
		\end{align*}
		and after applying the Hölder inequality and regrouping the terms, we obtain
		\begin{align*}
			\begin{split}
			& \|\big(\brho^{\varepsilon_k}_{\mu}-\brho^{\varepsilon_j}_{\mu}\big)(s,\cdot)\|_{B^{\theta+\delta}_{1,1}}                                                                                                                                                                                                        \\
			 & \leq C\|b^{\varepsilon_k}\|_{L^r(B^{\beta}_{p,q})}\Bigg(\Big(\int_0^s\|\big(\brho^{\varepsilon_k}_{\mu}-\brho^{\varepsilon_j}_{\mu}\big)(v,\cdot)\|_{B^{\theta+\delta}_{1,1}}^{r'} v^{-r'\gamma}(s-v)^{-r'\frac{1+2\theta+2\delta}{\alpha}}dv\Big)^{\frac{1}{r'}}\\
			 &\quad+\|\brho^{\varepsilon_k}_{\mu}-\brho^{\varepsilon_j}_{\mu}\|_{L^\infty_\gamma((0,S],B^{-\beta-\theta}_{p',1})}\Big(\int_0^s\|\brho^{\varepsilon_j}_{\mu}(v,\cdot)\|_{B^{\theta+\delta}_{1,1}}^{r'} v^{-r'\gamma}(s-v)^{-r'\frac{1+2\theta+2\delta}{\alpha}}dv\Big)^{\frac{1}{r'}}\Bigg)\\
			 &\quad+C\|b^{\varepsilon_k}-b^{\varepsilon_j}\|_{L^{\bar{r}}(B^{\beta-\delta'}_{p,q})}\Big(\int_0^s\|\brho^{\varepsilon_j}_{\mu}(v,\cdot)\|_{B^{\theta+\delta}_{1,1}}^{\bar{r}'} v^{-\bar{r}'\gamma}(s-v)^{-\bar{r}'\frac{1+2\theta+2\delta}{\alpha}}dv\Big)^{\frac{1}{\bar{r}'}}.
			 \end{split}
		\end{align*}

		We apply Lemma \ref{lemma-beta-function}:
		\begin{align}
			\label{prf-cauchy-density}
			\begin{split}
				& \|\big(\brho^{\varepsilon_k}_{\mu}-\brho^{\varepsilon_j}_{\mu}\big)(s,\cdot)\|_{B^{\theta+\delta}_{1,1}}                                                                                                                                                                                                        \\
				& \leq C\|b^{\varepsilon_k}\|_{L^r(B^{\beta}_{p,q})}\sup_{v\in(0,S]}\|\big(\brho^{\varepsilon_k}_{\mu}-\brho^{\varepsilon_j}_{\mu}\big)(v,\cdot)\|_{B^{\theta+\delta}_{1,1}} s ^{\frac{1}{r'}-\gamma-\frac{1+2\theta+2\delta}{\alpha}}\\
			 	&\quad+C\Big(\|b^{\varepsilon_k}\|_{L^r(B^{\beta}_{p,q})}\|\brho^{\varepsilon_k}_{\mu}-\brho^{\varepsilon_j}_{\mu}\|_{L^\infty_\gamma((0,S],B^{-\beta-\theta}_{p',1})} s ^{\frac{1}{r'}-\gamma-\frac{1+2\theta+2\delta}{\alpha}}\\
			 	&\quad+\|b^{\varepsilon_k}-b^{\varepsilon_j}\|_{L^{\bar{r}}(B^{\beta-\delta'}_{p,q})} s ^{\frac{1}{\bar{r}'}-\gamma-\frac{1+2\theta+2\delta}{\alpha}}\Big)\sup_{v\in(0,S]}\|\brho^{\varepsilon_j}_{\mu}(v,\cdot)\|_{B^{\theta+\delta}_{1,1}}.
			\end{split}
		\end{align}

		This is indeed possible because $\gamma<\frac{1}{2r'}$ from the choice of $\gamma$ in the proof of Lemma \ref{lemma-density-estimate} and $1+2\theta+2\delta<\alpha(1-\frac{1}{r})$ from the upper bound on $\theta$ given by \eqref{good-relation} for arbitrarily small $\delta$. Moreover, since $\beta<-1-2\theta$ from \eqref{C3} and $\frac{1}{r'}-\gamma+\frac{\beta-\delta}{\alpha}$ from the proof of Lemma \ref{lemma-density-estimate}, we also have
		\begin{align*}
			\frac{1}{r'}-\gamma-\frac{1+2\theta+2\delta}{\alpha}>0.
		\end{align*}
		
		The estimate \eqref{prf-cauchy-density} contains the norm of the mollified density that has not been studied yet. Once we have it, we can circularize \eqref{prf-cauchy-density} and the proof will be concluded. Therefore, we prove a uniform in $\varepsilon$ and time control on $\|\brho^{\varepsilon_j}_{\mu}(v,\cdot)\|_{B^{\theta+\delta}_{1,1}}$ by exploiting the additional regularity of $\mu$ in Besov norm with integrability index $1$, i.e., using $\mu\in B^{\theta+\delta}_{1,\infty}$. Write the Duhamel representation and apply the Young inequality \eqref{besov-prop-y} and the fact that $p^\alpha_{v}\in B^0_{1,\infty}$ which follows from \eqref{besov-prop-e3}:
		\begin{align*}
			\|\brho^{\varepsilon_j}_{\mu}(v,\cdot)\|_{B^{\theta+\delta}_{1,1}}&\leq\|\mu\ast p^\alpha_{ v}\|_{B^{\theta+\delta}_{1,1}}+\int_0^v\|(\mathcal{B}^{\varepsilon}_{\brho^{\varepsilon}_{\mu}}(r,\cdot)\brho^{\varepsilon}_{\mu}(r,\cdot))\ast\nabla p^\alpha_{v-r}\|_{B^{\theta+\delta}_{1,1}}dr\\
			&\leq C\|\mu\|_{B^{\theta+\delta}_{1,\infty}}\|p^\alpha_{ v}\|_{B^0_{1,\infty}}+C\int_0^v\|\mathcal{B}^{\varepsilon}_{\brho^{\varepsilon}_{\mu}}(r,\cdot)\brho^{\varepsilon}_{\mu}(r,\cdot)\|_{B^{-\theta}_{1,1}}\|\nabla p^\alpha_{v-r}\|_{B^{2\theta+\delta}_{1,1}}dr\\
			&\leq C\|\mu\|_{B^{\theta+\delta}_{1,\infty}}+C\int_0^v\|\mathcal{B}^{\varepsilon}_{\brho^{\varepsilon}_{\mu}}(r,\cdot)\brho^{\varepsilon}_{\mu}(r,\cdot)\|_{B^{-\theta}_{1,1}}\|p^\alpha_{v-r}\|_{B^{1+2\theta+\delta}_{1,1}}dr.
		\end{align*}
		
		Then, applying the product rule \eqref{ineq-product-rule} with parameters
		\begin{align*}
			\lambda,\lambda_1=\theta,\quad\lambda_2=\theta+\delta,\quad\ell=1,
		\end{align*}
		we derive
		\begin{align*}
			&\|\brho^{\varepsilon_j}_{\mu}(v,\cdot)\|_{B^{\theta+\delta}_{1,1}}\\
			&\leq C\|\mu\|_{B^{\theta+\delta}_{1,\infty}}+C\int_0^v\|\mathcal{B}^{\varepsilon}_{\brho^{\varepsilon}_{\mu}}(r,\cdot)\|_{B^{-\theta}_{\infty,\infty}}\|\brho^{\varepsilon}_{\mu}(r,\cdot)\|_{B^{\theta+\delta}_{1,1}}\|p^\alpha_{v-r}\|_{B^{1+2\theta+\delta}_{1,1}}dr\\
			&\stackrel{\eqref{besov-prop-y}, \eqref{besov-prop-hk}}{\leq} C\|\mu\|_{B^{\theta+\delta}_{1,\infty}}+C\sup_{r\in(0,v]}\|\brho^{\varepsilon}_{\mu}(r,\cdot)\|_{B^{\theta+\delta}_{1,1}}\int_0^v\|b^\varepsilon(r,\cdot)\|_{B^\beta_{p,q}}\|\brho^{\varepsilon}_{\mu}(r,\cdot)\|_{B^{-\beta-\theta}_{p',1}}(v-r)^{-\frac{1+2\theta+\delta}{\alpha}}dr\\
			&\leq C\|\mu\|_{B^{\theta+\delta}_{1,\infty}}+C\sup_{r\in(0,v]}\|\brho^{\varepsilon}_{\mu}(r,\cdot)\|_{B^{\theta+\delta}_{1,1}}\|b^\varepsilon\|_{L^r(B^\beta_{p,q})}\Big(\int_0^v r^{-r'\gamma}(v-r)^{-r'\frac{1+2\theta+\delta}{\alpha}}dr\Big)^{\frac{1}{r'}}\\
			&\stackrel{\text{Lemma } \ref{lemma-beta-function}}{\leq} C\|\mu\|_{B^{\theta+\delta}_{1,\infty}}+C\sup_{r\in(0,v]}\|\brho^{\varepsilon}_{\mu}(r,\cdot)\|_{B^{\theta+\delta}_{1,1}}\|b\|_{L^r(B^\beta_{p,q})} v^{\frac{1}{r'}-\gamma-\frac{1+2\theta+\delta}{\alpha}}.
		\end{align*}
		
		In the last inequality, we also used Proposition \ref{prop-b-conv} to bound the norm of $b^\varepsilon$. Now, taking the supremum over the time interval $(0,v]$ and regrouping the terms, in small time we obtain 
		\begin{align}
			\label{prf-unif-bound-theta}
			\begin{split}
			\sup_{t\in(0,v]}\|\brho^{\varepsilon}_{\mu}(t,\cdot)\|_{B^{\theta+\delta}_{1,1}}&\leq\frac{C\|\mu\|_{B^{\theta+\delta}_{1,1}}}{1-C\|b\|_{L^r(B^\beta_{p,q})} v^{\frac{1}{r'}-\gamma-\frac{1+2\theta+\delta}{\alpha}}}=:\bar{C}.
			\end{split}
		\end{align}

		The uniform estimate on $\sup_{t\in(0,v]}\|\brho^{\varepsilon}_{\mu}(t,\cdot)\|_{B^{\theta+\delta}_{1,1}}$ is thus proven. Let us finish showing that $\big(\brho_{\mu}^{\varepsilon_k}\big)_{k\geq1}$ is a Cauchy sequence in $L^\infty((0,S],L^1)$. Substituting \eqref{prf-unif-bound-theta} into \eqref{prf-cauchy-density} and taking the supremum in the l.h.s. of it, we have that
		\begin{align*}
			& \sup_{v\in(0,S]}\|\big(\brho^{\varepsilon_k}_{\mu}-\brho^{\varepsilon_j}_{\mu}\big)(v,\cdot)\|_{B^{\theta+\delta}_{1,1}}                                                                                                                                                                                                      \\
			 &\leq C\|b\|_{L^r(B^{\beta}_{p,q})}\sup_{v\in(0,S]}\|\big(\brho^{\varepsilon_k}_{\mu}-\brho^{\varepsilon_j}_{\mu}\big)(v,\cdot)\|_{B^{\theta+\delta}_{1,1}}S^{\frac{1}{r'}-\gamma-\frac{1+2\theta+2\delta}{\alpha}}\\
			 &\quad+\bar{C}\Big(\|b\|_{L^r(B^{\beta}_{p,q})}\|\brho^{\varepsilon_k}_{\mu}-\brho^{\varepsilon_j}_{\mu}\|_{L^\infty_\gamma((0,S],B^{-\beta-\theta}_{p',1})}S^{\frac{1}{r'}-\gamma-\frac{1+2\theta+2\delta}{\alpha}}+\|b^{\varepsilon_k}-b^{\varepsilon_j}\|_{L^{\bar{r}}(B^{\beta-\delta'}_{p,q})}S^{\frac{1}{\bar{r}'}-\gamma-\frac{1+2\theta+2\delta}{\alpha}}\Big).
		\end{align*}

		By the first part of the current lemma $\|\brho^{\varepsilon_k}_{\mu}-\brho^{\varepsilon_j}_{\mu}\|_{L^\infty_\gamma((0,S],B^{-\beta-\theta}_{p',1})}\xrightarrow[k,j\to\infty]{}0$ and by Proposition \ref{prop-b-conv} $\|b^{\varepsilon_k}-b^{\varepsilon_j}\|_{L^{\bar{r}}(B^{\beta-\delta'}_{p,q})}\xrightarrow[k,j\to\infty]{}0$. Hence, by the embedding \eqref{besov-prop-e1} we conclude that
		\begin{align*}
			\sup_{s\in(0,S]} & \|\big(\brho^{\varepsilon_k}_{\mu}-\brho^{\varepsilon_j}_{\mu}\big)(s,\cdot)\|_{L^1}\xrightarrow[k,j\to\infty]{}0.
		\end{align*}

		This completes the proof.
	\end{proof}
\end{lemma}
		
\begin{remark}[About the smoothness assumption on $\mu$]
	\label{remark-mu-reg}
	We highlight that in order to prove that the limit of the sequence formed by $\brho_{\mu}^{\varepsilon_k}$ is a probability measure, we need to assume that $\mu\in B^{\theta+\delta}_{1,\infty}$. The absence of the regularity assumption yields a time singularity that cannot be overcome. This condition on $\mu$ did not appear in \cite{deraynal2023multidimensionalstabledrivenmckeanvlasov}, and it somehow naturally appears in the current rather singular setting. Therein, the aforementioned result could be proven by simply using the fact that $\mu$ and $\brho^\varepsilon_{\mu}$ are probability measures by using the embedding \eqref{besov-prop-e1} $B^0_{1,1}\hookrightarrow L^1$ and by applying the product rule \eqref{besov-prop-pr1} which works for spaces of positive regularity. On the other hand, if in our distributional framework (with a more singular drift) we try to establish estimates in the space $B^0_{1,1}$, we will naturally obtain the norm $B^{\theta+\delta}_{1,1}$ of the mollified density on the r.h.s. because of the product rule \eqref{ineq-product-rule} which handles negative regularity (as in, for example, Lemma \ref{lemma-density-estimate}, see the controls \eqref{prf-quad-1}, \eqref{prf-quad-2}). This suggests that the norm has to be controlled uniformly in time. We thus control it again via the Duhamel expansion, and the absence of a time weight makes it impossible to compensate the time singularities of the heat kernel. Therefore, the appearing regularity $\theta+\delta$ has to be demanded from the initial measure and not from the heat kernel.
\end{remark}

\begin{lemma}[Existence of a solution for the limit Fokker-Planck equation]
	\label{lemma-density-existence}
	Let conditions \eqref{C3} and \eqref{mu-range} be satisfied. Let $(\varepsilon_k)_{k\geq1}$ be a decreasing sequence and $\brho_{\mu}$ be the limit point of $(\brho^{\varepsilon_k}_{\mu})_{k\geq1}$, which exists from Lemma \ref{lemma-density-conv}. Then $\brho_{\mu}$ satisfies the Fokker-Planck equation \eqref{pde-main} in a distributional sense, and it satisfies the Duhamel representation \eqref{duhamel-main}.

	\begin{proof}
		The proof of this statement is close to that of Lemma 10 in \cite{deraynal2022multidimensional}, up to balancing the regularity indices and up to choosing the right product rule in our singular framework. Therefore, we skip the details for the parts of the proof where the similarities with \cite{deraynal2022multidimensional} are strong, and we write more details where needed. Let $S<\mathcal{T}$. From the weak formulation for $\brho_{\mu}^{\varepsilon_k}$, we have that for any $\phi\in\mathcal{D}([0,S)\times\R^d)$ and for fixed $k\in\N$, $\brho_{\mu}$ satisfies
		\begin{align*}
			&\int_0^S\int_{\R^d}\brho_{\mu}(s,x)\big(\partial_s\phi+\mathcal{B}_{\brho_{\mu}}\cdot\nabla\phi+\mathcal{L}^*(\phi)\big)(s,x)dxds\\
			&=-\int_{\R^d}\phi(t,x)\mu(dx)+\Delta^1_{\brho_{\mu},\brho_{\mu}^{\varepsilon_k}}(\phi)+\Delta^2_{\brho_{\mu},\brho_{\mu}^{\varepsilon_k}}(\phi),
		\end{align*}
		where
		\begin{align*}
			\Delta^1_{\brho_{\mu},\brho_{\mu}^{\varepsilon_k}}(\phi)=\int_0^S\int_{\R^d}\big(\brho_{\mu}-\brho_{\mu}^{\varepsilon_k}\big)(s,x)\big(\partial_s\phi+\mathcal{L}^*(\phi)\big)(s,x)dxds
		\end{align*}
		and
		\begin{align*}
			\Delta^2_{\brho_{\mu},\brho_{\mu}^{\varepsilon_k}}(\phi)=\int_0^S\int_{\R^d}\big(\mathcal{B}_{\brho_{\mu}}\brho_{\mu}-\mathcal{B}^{\varepsilon_k}_{\brho^{\varepsilon_k}_{\mu}}\brho^{\varepsilon_k}_{\mu}\big)(s,x)\cdot\nabla\phi(s,x)dxds.
		\end{align*}

		We aim at showing that $|(\Delta^1_{\brho_{\mu},\brho_{\mu}^{\varepsilon_k}}+\Delta^2_{\brho_{\mu},\brho_{\mu}^{\varepsilon_k}})(\phi)|\xrightarrow[k\to+\infty]{}0$. 

		\vspace{3mm}

		First, since $(\partial_s+\mathcal{L}^*)(\phi)\in L^1([0,S),L^\infty)$, we get from Lemma \ref{lemma-density-conv} that $|\Delta^1_{\brho_{\mu},\brho_{\mu}^{\varepsilon_k}}(\phi)|\xrightarrow[k\to+\infty]{}0$.

		\vspace{3mm}

		For the second term $\Delta^2_{\brho_{\mu},\brho_{\mu}^{\varepsilon_k}}(\phi)$, we use the same techniques as in Lemma \ref{lemma-density-conv}. For $\varepsilon>0$ arbitrarily small, the duality \eqref{besov-prop-d} yields
		\begin{align*}
			 |\Delta^2_{\brho_{\mu},\brho_{\mu}^{\varepsilon_k}}(\phi)|& \leq \Big|\int_0^S\int_{\R^d}\big(\brho_{\mu}-\brho_{\mu}^{\varepsilon_k}\big)(s,x)\big(\mathcal{B}_{\brho_{\mu}}\cdot\nabla\phi\big)(s,x)dsdx\Big|                                                                                                                                                                                                              \\
			 & \quad+\Big|\int_0^S\int_{\R^d}\brho_{\mu}^{\varepsilon_k}(s,x)\big((\mathcal{B}_{\brho_{\mu}}-\mathcal{B}^{\varepsilon_k}_{\brho^{\varepsilon_k}_{\mu}})\cdot\nabla\phi\big)(s,x)dsdx\Big|                                                                                                                                                                             \\
			 & \leq\|\brho_{\mu}-\brho_{\mu}^{\varepsilon_k}\|_{L^\infty_\gamma((0,S],B^{\theta}_{1,1})}\int_0^S s ^{-\gamma}\|\big(\mathcal{B}_{\brho_{\mu}}\cdot\nabla\phi\big)(s,\cdot)\|_{B^{-\theta}_{\infty,\infty}}ds                                                                                                                                          \\
			 & \quad+\int_0^S\|\brho_{\mu}^{\varepsilon_k}(s,\cdot)\|_{B^{\theta+\varepsilon}_{1,1}}\|\big((\mathcal{B}_{\brho_{\mu}}-\mathcal{B}^{\varepsilon_k}_{\brho^{\varepsilon_k}_{\mu}})\cdot\nabla\phi\big)(s,\cdot)\|_{B^{-\theta-\varepsilon}_{\infty,\infty}}ds.
		\end{align*}

		Unlike in \cite[Lemma 10]{deraynal2022multidimensional}, we cannot simply take the supremum of $\nabla\phi$ out since we have a product of negative regularity. Therefore, applying the product rule \eqref{besov-prop-pr2} with $\rho>\theta$ implies
		 \begin{align*}
			 |\Delta^2_{\brho_{\mu},\brho_{\mu}^{\varepsilon_k}}(\phi)|& \leq\|\nabla\phi\|_{B^\rho_{\infty,\infty}}\Bigg(\|\brho_{\mu}-\brho_{\mu}^{\varepsilon_k}\|_{L^\infty_\gamma((0,S],B^{-\beta-\theta}_{p',1})}\int_0^S s ^{-\gamma}\|b(s,\cdot)\|_{B^\beta_{p,q}}\|\brho_{\mu}(s,\cdot)\|_{B^{-\beta-\theta}_{p',1}}ds                                                                                                               \\
			 & \quad+\int_0^S\Big(\|(b-b^{\varepsilon_k})(s,\cdot)\|_{B^{\beta-\varepsilon}_{p,q}}\|\brho^{\varepsilon_k}_{\mu}(s,\cdot)\|_{B^{-\beta-\theta}_{p',1}}^2\\
			 &\quad+\|b(s,\cdot)\|_{B^{\beta}_{p,q}}\|\big(\brho_{\mu}-\brho_{\mu}^{\varepsilon_k}\big)(s,\cdot)\|_{B^{-\beta-\theta}_{p',1}}\|\brho^{\varepsilon_k}_{\mu}(s,\cdot)\|_{B^{-\beta-\theta}_{p',1}}\Big)ds\Bigg),
		 \end{align*}
		 and the Hölder inequality in time together with Lemma \ref{lemma-beta-function}, gives
		 \begin{align*}
			|\Delta^2_{\brho_{\mu},\brho_{\mu}^{\varepsilon_k}}(\phi)|& \leq C\|\nabla\phi\|_{B^{\rho}_{\infty,\infty}}\Bigg(\|\brho_{\mu}-\brho_{\mu}^{\varepsilon_k}\|_{L^\infty_\gamma((0,S],B^{-\beta-\theta}_{p',1})}\|b\|_{L^r(B^\beta_{p,q})}\Big(\int_0^S s ^{-2r'\gamma}\Big)^{\frac{1}{r'}}ds                                                                                                                                          \\
			 & \quad+\|b-b^{\varepsilon_k}\|_{L^{\bar{r}}(B^{\beta-\varepsilon}_{p,q})}\Big(\int_0^S s ^{-2\bar{r}'\gamma}\Big)^{\frac{1}{\bar{r}'}}ds\\
			 &\quad+\|b\|_{L^r(B^\beta_{p,q})}\|\brho_{\mu}-\brho_{\mu}^{\varepsilon_k}\|_{L^\infty_\gamma((0,S],B^{-\beta-\theta}_{p',1})}\Big(\int_0^S s ^{-2r'\gamma}\Big)^{\frac{1}{r'}}\Bigg)                                       \\
			 & \leq C\|\nabla\phi\|_{B^\rho_{\infty,\infty}}\Big(\|b-b^{\varepsilon_k}\|_{L^{\bar{r}}(B^{\beta-\varepsilon}_{p,q})}+\|\brho_{\mu}-\brho_{\mu}^{\varepsilon_k}\|_{L^\infty_\gamma((0,S],B^{-\beta-\theta}_{p',1})}\Big)S^{\frac{1}{\bar{r}'}-2\gamma}.
		 \end{align*}

		We have that $\frac{1}{\bar{r}'}-2\gamma>0$ because $\gamma<\frac{1}{2\bar{r}'}$ from Lemma \ref{lemma-density-estimate}. Thus, from Proposition \ref{prop-b-conv} and Lemma \ref{lemma-density-conv}, we derive that $|\Delta^2_{\brho_{\mu},\brho_{\mu}^{\varepsilon_k}}(\phi)|\xrightarrow[k\to+\infty]{}0$. Finally, it follows that $\brho_{\mu}$ satisfies \eqref{pde-main} in a distributional sense.

		\vspace{3mm}

		The Duhamel representation \eqref{duhamel-main} is established similarly by formally replacing $(s,x)\mapsto\phi(s,x)$ with $(s,x)\mapsto p^\alpha_{s}(s,x)$. See Lemmas 3 and 9 of \cite{deraynal2022multidimensional} for more details.
	\end{proof}
\end{lemma}

\begin{lemma}[Uniqueness of the limit Fokker-Planck solution]
	\label{lemma-density-uniqueness}
	Under the conditions \eqref{C3} and \eqref{mu-range}, the Fokker-Planck equation \eqref{pde-main} admits at most one solution in $L^\infty_{\gamma}((0,S],B^{-\beta-\theta}_{p',1})$ for any $S<\mathcal{T}$.

	\begin{proof}
		The proof is of the same idea as in Lemma 9 of \cite{deraynal2023multidimensionalstabledrivenmckeanvlasov}. The main difference lies in how the non-linear product involving the density in the Duhamel representation is handled, i.e., a different product rule is applied here. The nature of the problem is the same as in the proof of Lemma \ref{lemma-density-estimate}. However, the techniques needed to prove the current lemma have already been developed in the previous lemmas, thus, we skip the details of the exact calculations. Let $\brho^1_{\mu}$ and $\brho^2_{\mu}$ be two possible solutions to \eqref{pde-main}. Then from Duhamel representation of a solution, for a.e. $s\in[0,S]$, $y\in\R^d$,
		\begin{align*}
			(\brho^1_{\mu}-\brho^2_{\mu})(s,y)=-\int_0^s\Big(\big(\mathcal{B}_{\brho^1_{\mu}}(v,\cdot)\brho^1_{\mu}(v,\cdot)-\mathcal{B}_{\brho^2_{\mu}}(v,\cdot)\brho^2_{\mu}(v,\cdot)\big)\ast\nabla p^\alpha_{s-v}\Big)(y)dv.
		\end{align*}

		Similarly to the computations made in Lemma \ref{lemma-density-conv}, performing first an integration by parts, we obtain
		\begin{align*}
			 & \|\big(\mathcal{B}_{\brho^1_{\mu}}(v,\cdot)\brho^1_{\mu}(v,\cdot)-\mathcal{B}_{\brho^2_{\mu}}(v,\cdot)\brho^2_{\mu}(v,\cdot)\big)\ast\nabla p^\alpha_{s-v}\|_{B^{-\beta-\theta}_{p',1}}                                                                                                                                         \\
			 & \leq  C\Big(\|\brho^1_{\mu}(v,\cdot)\|_{B^{-\beta-\theta}_{p',1}}+\|\brho^2_{\mu}(v,\cdot)\|_{B^{-\beta-\theta}_{p',1}}\Big)\Big(\|b(v,\cdot)\|_{B^{\beta}_{p,q}}\\
			 &\quad+\|\Div(b(v,\cdot))\|_{B^{\beta}_{p,q}}\Big)\|(\brho^1_{\mu}-\brho^2_{\mu})(v,\cdot)\|_{B^{-\beta-\theta}_{p',1}}\|p^\alpha_{s-v}\|_{B^{-\beta+\delta}_{1,1}}.
		\end{align*}

		From the Hölder inequality, we derive
		\begin{align*}
			 & \|(\brho^1_{\mu}-\brho^2_{\mu})(s,\cdot)\|_{B^{-\beta-\theta}_{p',1}}                                                                                                                                                                                           \\
			 & \leq C\Big(\|b\|_{L^r(B^{\beta}_{p,q})}+\|\Div(b)\|_{L^r(B^{\beta}_{p,q})}\Big)\|\brho^1_{\mu}-\brho^2_{\mu}\|_{L^\infty_\gamma((0,S],B^{-\beta-\theta}_{p',1})}\Big(\int_0^s v^{-2r'\gamma}(s-v)^{r'\frac{\beta-\delta}{\alpha}}dv\Big)^{\frac{1}{r'}},
		\end{align*}
		which implies
		\begin{align*}
			& s ^\gamma\|(\brho^1_{\mu}-\brho^2_{\mu})(s,\cdot)\|_{B^{-\beta-\theta}_{p',1}}\\
			&\leq C\Big(\|b\|_{L^r(B^{\beta}_{p,q})}+\|\Div(b)\|_{L^r(B^{\beta}_{p,q})}\Big)\|\brho^1_{\mu}-\brho^2_{\mu}\|_{L^\infty_\gamma((0,S],B^{-\beta-\theta}_{p',1})}S^{\frac{1}{r'}-\gamma+\frac{\beta-\delta}{\alpha}}.
		\end{align*}

		Since $\frac{1}{r'}-\gamma+\frac{\beta-\delta}{\alpha}>0$ by \eqref{C3} and $S$ can be taken small enough, we obtain the claim.
	\end{proof}
\end{lemma}

\section{Well-posedness of the McKean-Vlasov SDE in short time}
\label{section-sde}

To use the previous results, we assume here that the time horizon $T>0$ is sufficiently small. We prove an important result about the mollified non-linear convolution drift.

\begin{lemma}[Integrability property of the mollified drift]
	\label{lemma-mollified-drift-bound}
	Assume that conditions \eqref{C3} and \eqref{mu-range} hold. Assume that \eqref{good-relation} or \eqref{good-relation-enhaced} is satisfied. Then for any $S\in(0,\mathcal{T})$, $r_\theta\in(\alpha,r)$,
	\begin{align*}
		\mathcal{B}^\varepsilon_{\brho^\varepsilon_{\mu}}\in L^{r_\theta}((0,S],B^{-\theta}_{\infty,\infty}),
	\end{align*}
	and
	\begin{align}
		\label{ineq-B-eps}
		\|\mathcal{B}^\varepsilon_{\brho^\varepsilon_{\mu}}(s,\cdot)\|_{B^{-\theta}_{\infty,\infty}} & \leq C\|b\|_{L^r(B^\beta_{p,q})}S^{\frac{1}{r_\theta}-\frac{1}{r}-\gamma}.
	\end{align}

	\begin{proof}
		For any $s\in(0,S]$, by \eqref{besov-prop-y}
		\begin{align*}
			\|\mathcal{B}^\varepsilon_{\brho^\varepsilon_{\mu}}(s,\cdot)\|_{B^{-\theta}_{\infty,\infty}} & \leq C\|b^\varepsilon(s,\cdot)\|_{B^\beta_{p,q}}\|\brho^\varepsilon_{\mu}(s,\cdot)\|_{B^{-\beta-\theta}_{p',1}}.
		\end{align*}

		Then, applying the Hölder inequality, Proposition \ref{prop-b-conv} and Proposition \ref{prop-main-pde}, we get
		\begin{align*}
			\|\mathcal{B}^\varepsilon_{\brho^\varepsilon_{\mu}}\|_{L^{r_\theta}(B^{-\theta}_{\infty,\infty})} & \leq C\|b^\varepsilon\|_{L^r(B^\beta_{p,q})}\Big(\int_0^S\|\brho^\varepsilon_{\mu}(s,\cdot)\|_{B^{-\beta-\theta}_{p',1}}^{\frac{r_\theta r}{r-r_\theta}}ds\Big)^{\frac{r-r_\theta}{r_\theta r}}                       \\
			                                                                                                    & \leq C\|b^\varepsilon\|_{L^r(B^\beta_{p,q})}\|\brho^\varepsilon_{\mu}\|_{L^\infty_{\gamma}(B^{-\beta-\theta}_{p',1})}\Big(\int_0^S s^{-\gamma\frac{r_\theta r}{r-r_\theta}}ds\Big)^{\frac{r-r_\theta}{r_\theta r}} \\
			                                                                                                    & \leq C\|b\|_{L^r(B^\beta_{p,q})}\|\brho^\varepsilon_{\mu}\|_{L^\infty_{\gamma}(B^{-\beta-\theta}_{p',1})}S^{\frac{1}{r_\theta}-\frac{1}{r}-\gamma}.
		\end{align*}

		Above, $\gamma\frac{r_\theta r}{r-r_\theta}<1\iff r_\theta<\frac{r}{1+r\gamma}<r$, and the result follows.
	\end{proof}

\end{lemma}

\subsection{Martingale problem well-posedness}

To prove well-posedness of the non-linear martingale problem, we first prove the existence of a measure satisfying Definition \ref{def-martingale-problem} by using standard machinery of showing tightness of a sequence of measures and identifying its limit. The main result states as follows.

\begin{proposition}
	\label{prop-martingale-problem-wp}
	Let condition \eqref{C3} be satisfied. Let $\theta,r$ satisfy \eqref{good-relation}. Denote by $(\Prob^\alpha_\varepsilon)_{\varepsilon>0}$ the solution to the non-linear martingale problem associated with \eqref{main-sde-mollified}. Then there exists a unique limit point of a converging subsequence $(\Prob^\alpha_{\varepsilon_k})_{k\geq1}$, $\varepsilon_k\downarrow_{k\to\infty}0$, in $\mathcal{P}(\Omega_\alpha)$ equipped with its weak topology, such that it solves the non-linear martingale problem associated with \eqref{main-sde} in the sense of Def. \ref{def-martingale-problem}.

	\begin{proof}
		\textbf{Tightness.}
		Recall that the mollified McKean-Vlasov SDE \eqref{main-sde-mollified} is well-posed in a weak sense for any $\varepsilon>0$, which implies well-posedness of the corresponding martingale problem. For each $\varepsilon>0$, let $\Prob^\alpha_\varepsilon$ denote its solution on $\Omega_\alpha$, and $(X^{\varepsilon,\mu}_s)_{s\in[0,S]}$ its canonical process. Writing
		\begin{align}
			\label{mart-pr-X-diff-1}
			X^{\varepsilon,\mu}_s-X^{\varepsilon,\mu}_v=\int_v^s\mathcal{B}^\varepsilon_{\brho^\varepsilon_{\mu}}(r,X^{\varepsilon,\mu}_r)dr+\mathcal{W}_s-\mathcal{W}_v
		\end{align}
		we see that although for the mollified drift $\mathcal{B}^\varepsilon_{\brho^\varepsilon_{\mu}}(r,\cdot)$ the above integral is well-defined, in order to get uniform in $\varepsilon$ estimates, we handle $\mathcal{B}^\varepsilon_{\brho^\varepsilon_{\mu}}(r,\cdot)$ as a distribution in $B^{-\theta}_{\infty,\infty}$. To this end, we apply the Zvonkin transform as in \cite{deraynal2022-linear-besov}. Namely, we introduce the PDE
		\begin{align}
			\label{mart-pr-zvonkin}
			(\partial_s+\mathcal{B}^\varepsilon_{\brho^\varepsilon_{\mu}}\cdot D+\mathcal{L}^\alpha)u^\varepsilon_k(s,x)=-\big(\mathcal{B}^\varepsilon_{\brho^\varepsilon_{\mu}}(s,x)\big)_k\ \text{ on }[0,S)\times\R^d,\quad u^\varepsilon_k(S,\cdot)=0\text{ on }\R^d,
		\end{align}
		where $\big(\mathcal{B}^\varepsilon_{\brho^\varepsilon_{\mu}}(s,x)\big)_k$ is the $k$-th component of $\mathcal{B}^\varepsilon_{\brho^\varepsilon_{\mu}}(s,x)$, $k=1,...,d$. Clearly, this PDE admits a unique classical solution for each $k=1,...,d$, $\varepsilon>0$. Denote for fixed $\varepsilon>0$, $u^\varepsilon=(u^\varepsilon_1,...,u^\varepsilon_d)$ the vector of solutions to \eqref{mart-pr-zvonkin}. Then applying Itô's formula to $(X^{\varepsilon,\mu}_r+u^\varepsilon(r,X^{\varepsilon,\mu}_r))_{r\in[v,s]}$ we get
		\begin{align}
			\label{mart-pr-X-diff-2}
			X^{\varepsilon,\mu}_s-X^{\varepsilon,\mu}_v=M_{v,s}(\alpha,u^\varepsilon,X^{\varepsilon,\mu})+\big(u^\varepsilon(v,X^{\varepsilon,\mu}_v)-u^\varepsilon(s,X^{\varepsilon,\mu}_s)\big)+\mathcal{W}_s-\mathcal{W}_v,
		\end{align}
		where $M_{v,s}(\alpha,u^\varepsilon,X^{\varepsilon,\mu})$ is a martingale given by
		\begin{align*}
			&M_{v,s}(\alpha,u^\varepsilon,X^{\varepsilon,\mu})\\
			&=
			\begin{cases}
				\begin{aligned}
				&\int_v^sDu^\varepsilon(r,X_r^{\varepsilon,\mu})\cdot dW_r,\text{ where }W\text{ is a Brownian motion, if } \alpha=2,\\
				&\int_v^s\int_{\R^d\setminus\{0\}}(u^\varepsilon(r,X^{\varepsilon,\mu}_{r-}+x)-u^\varepsilon(r,X^{\varepsilon,\mu}_{r-}))\tilde{N}(dr,dx),\begin{aligned}[t]
					\text{ where }\tilde{N}\text{ is the compensated} \\
					\text{Poisson measure, if } \alpha<2.
				\end{aligned}
				\end{aligned}
			\end{cases}
		\end{align*}
		Therefore, instead of \eqref{mart-pr-X-diff-1} we obtained the representation \eqref{mart-pr-X-diff-2}. Write for some $\lambda>0$
		\begin{align*}
			 & \Exp^{\Prob^\alpha_\varepsilon}\big[|X^{\varepsilon,\mu}_s-X^{\varepsilon,\mu}_v|^{\lambda}\big]\leq\Exp\big[|M_{v,s}(\alpha,u^\varepsilon,X^{\varepsilon,\mu})|^\lambda\big]\\
			 &\quad+\Exp^{\Prob^\alpha_\varepsilon}\big[|u^\varepsilon(s,X^{\varepsilon,\mu}_s)-u^\varepsilon(v,X^{\varepsilon,\mu}_s)|^{\lambda}\big]+\Exp^{\Prob^\alpha_\varepsilon}\big[|u^\varepsilon(s,X^{\varepsilon,\mu}_v)-u^\varepsilon(s,X^{\varepsilon,\mu}_s)|^{\lambda}\big]+\Exp^{\Prob^\alpha_\varepsilon}\big[|\mathcal{W}_s-\mathcal{W}_v|^{\lambda}\big].
		\end{align*}
		To conclude, we use the regularity estimates on $u^\varepsilon$ given by Proposition 9, \cite{deraynal2022-linear-besov}. More precisely, from \cite{deraynal2022-linear-besov}, we have
		\begin{align*}
			\Exp^{\Prob^\alpha_\varepsilon}\big[|u^\varepsilon(s,X^{\varepsilon,\mu}_s)-u^\varepsilon(v,X^{\varepsilon,\mu}_s)|^{\lambda}\big]&\leq C|s-v|^{\frac{\lambda}{\alpha}(-\theta+\alpha-\frac{\alpha}{r_\theta})},
		\end{align*}
		and
		\begin{align*}
			\Exp^{\Prob^\alpha_\varepsilon}\big[|u^\varepsilon(s,X^{\varepsilon,\mu}_v)-u^\varepsilon(s,X^{\varepsilon,\mu}_s)|^{\lambda}\big]&\leq\|Du^\varepsilon\|_{L^\infty\big(B^{-\theta-1+\alpha-\frac{\alpha}{r_\theta}}_{\infty,\infty}\big)}\Exp^{\Prob^\alpha_\varepsilon}\big[|X^{\varepsilon,\mu}_v-X^{\varepsilon,\mu}_s|^{\lambda}\big]\\
			&\leq C_S\Exp^{\Prob^\alpha_\varepsilon}\big[|X^{\varepsilon,\mu}_v-X^{\varepsilon,\mu}_s|^{\lambda}\big],
		\end{align*}
		where $C_S>0$ is small enough in small time. Let us handle the martingale part. For $\alpha=2$,
		\begin{align*}
			\Exp\big[|M_{v,s}(\alpha,u^\varepsilon,X^{\varepsilon,\mu})|^\lambda\big]&=\Exp\big[|\int_v^s|Du^\varepsilon(r,X_r^{\varepsilon,t,\mu})|^2dr|^{\frac{\lambda}{2}}\big]\\
			&\leq\|Du^\varepsilon\|_{L^\infty(L^\infty)}^\lambda(s-v)^{\frac{\lambda}{2}}.
		\end{align*}
		Here, $\|Du^\varepsilon\|_{L^\infty(L^\infty)}<+\infty$ as $\|Du^\varepsilon\|_{L^\infty\big(B^{-\theta-1+\alpha-\frac{\alpha}{r_\theta}}_{\infty,\infty}\big)}<+\infty$ from \cite{deraynal2022-linear-besov}, and we can take $\lambda>2$. Similarly, when $\alpha\in(1,2)$,
		\begin{align*}
			\Exp\big[|M_{t,s}(\alpha,u^\varepsilon,X^{\varepsilon,\mu})|^\lambda\big]&=\Exp\Big[\Big|\int_0^s\int_{\R^d\setminus\{0\}}|u^\varepsilon(r,X_{r-}^{\varepsilon,t,\mu}+x)-u^\varepsilon(r,X_{r-}^{\varepsilon,t,\mu})|^2\nu(dx)dr\Big|^{\frac{\lambda}{2}}\Big]\\
			&\leq\|Du^\varepsilon\|_{L^\infty\big(B^{-\theta-1+\alpha-\frac{\alpha}{r_\theta}}_{\infty,\infty}\big)}^\lambda\Big|\int_0^s\int_{\R^d\setminus\{0\}}x\nu(dx)dr\Big|^{\frac{\lambda}{2}}\\
			&\leq C\|Du^\varepsilon\|_{L^\infty\big(B^{-\theta-1+\alpha-\frac{\alpha}{r_\theta}}_{\infty,\infty}\big)}^\lambda s ^{\frac{\lambda}{2}}.
		\end{align*}

		It is then enough to apply the Kolmogorov tightness criterion (see Theorem 7.3, \cite{billing_prob_convergenc}) for $\alpha=2$, which yields from the above estimates
		\begin{align*}
			\Exp^{\Prob^2_\varepsilon}\big[|X^{\varepsilon,\mu}_s-X^{\varepsilon,\mu}_v|^{\lambda}\big]\leq C|s-v|^{1+\xi}, \quad\forall s,v\in[0,S],\quad\alpha=2,
		\end{align*}
		and Aldous criterion (see Proposition 34.9, \cite{Bass_2011-stoch-pr}) for $\alpha\in(1,2)$, which also yields from the above estimates and tightness of symmetric $\alpha$-stable processes in $\Omega_\alpha$,
		\begin{align*}
			\Exp^{\Prob^\alpha_\varepsilon}\big[|X^{\varepsilon,\mu}_s-X^{\varepsilon,\mu}_t|^{\lambda}\big]\leq C|S|^{\xi}, \quad\forall s\in[0,S],\quad\alpha\in(1,2),
		\end{align*}
		for some $\xi>0$, $C>0$.

		\vspace{2mm}

		\textbf{Limit points.}
		For a decreasing subsequence $(\varepsilon_k)_{k\geq0}$, let $(\Prob^\alpha_{\varepsilon_k})_{k\geq0}$ be a converging subsequence of $(\Prob^\alpha_{\varepsilon})_{\varepsilon>0}$ and let $\Prob^\alpha$ be its limit. We aim to show that the limit $\Prob^\alpha$ is a solution to the non-linear martingale problem associated with \eqref{main-sde}. First, note that for $\varepsilon>0,k\geq0$, from the Young inequality \eqref{besov-prop-y},
		\begin{align*}
			&\|(\mathcal{B}^{\varepsilon_k}_{\brho^{\varepsilon_k}_{\mu}}-\mathcal{B}_{\brho_{\mu}})(s,\cdot)\|_{B^{-\theta-\varepsilon}_{\infty,\infty}}\\
			&\leq C\|\brho^{\varepsilon_k}_{\mu}(s,\cdot)\|_{B^{-\beta-\theta}_{p',1}}\|(b^{\varepsilon_k}-b)(s,\cdot)\|_{B^{\beta-\varepsilon}_{p,q}}+\|b(s,\cdot)\|_{B^\beta_{p,q}}\|(\brho_{\mu}^{\varepsilon_k}-\brho_{\mu})(s,\cdot)\|_{B^{-\beta-\theta}_{p',1}},
		\end{align*}
		and thus
		\begin{align*}
			\|\mathcal{B}^{\varepsilon_k}_{\brho^{\varepsilon_k}_{\mu}}-\mathcal{B}_{\brho_{\mu}}\|_{L^{r_\theta}((0,S],B^{-\theta-\varepsilon}_{\infty,\infty})} & \leq C\Big(\|\brho^{\varepsilon_k}_{\mu}\|_{L^{\infty}_\gamma((0,S],B^{-\beta-\theta}_{p',1})}\|b^{\varepsilon_k}-b\|_{L^{\bar{r}}((0,S],B^{\beta-\varepsilon}_{p,q})}S^{\frac{1}{r_\theta}-\frac{1}{\bar{r}}-\gamma} \\
			                                                                                                                                                       & \quad+\|b\|_{L^{r}((0,S],B^\beta_{p,q})}\|\brho_{\mu}^{\varepsilon_k}-\brho_{\mu}\|_{L^\infty_\gamma((0,S],B^{-\beta-\theta}_{p',1})}S^{\frac{1}{r_\theta}-\frac{1}{r}-\gamma}\Big)                        \\
			                                                                                                                                                       & \xrightarrow[k\to+\infty]{}0,
		\end{align*}
		where $\bar{r}\in[1,+\infty]$ is defined in Proposition \ref{prop-b-conv} and the convergence follows by Propositions \ref{prop-main-pde} and \ref{prop-b-conv}. Now, thanks to \eqref{besov-prop-e4}, we have that for any $\phi\in C_0^\infty((0,S)\times\R^d,\R^d)$,
		\begin{align*}
			\int_0^s\int_{\R^d}\mathcal{B}^{\varepsilon_k}_{\brho^{\varepsilon_k}_{\mu}}(v,x)\cdot\phi(s,x)dxds\xrightarrow[k\to+\infty]{}\int_0^s\int_{\R^d}\mathcal{B}_{\brho_{\mu}}(v,x)\cdot\phi(s,x)dxds,
		\end{align*}
		meaning that $\int_0^s\mathcal{B}^{\varepsilon_k}_{\brho^{\varepsilon_k}_{\mu}}(v,\cdot)dv$ converges weakly to $\int_0^s\mathcal{B}_{\brho_{\mu}}(v,\cdot)dv$. By Lemma 5.1 in \cite{markov-processes-conv}, this means that $\Prob^\alpha$ solves the non-linear martingale problem associated with \eqref{main-sde}.

		\vspace{2mm}

		\textbf{Uniqueness of a limit.}
		First, note that whenever a solution to the martingale problem associated with \eqref{main-sde} exists, its marginal distributions are absolutely continuous for almost every time $s\in(0,S)$ and the density $\brho_{\mu}(s,\cdot)$ belongs to the Lebesgue-Besov space $L^\infty_\gamma((0,S],B^{-\beta-\theta}_{p',1})$ by Proposition \ref{prop-main-pde}. Moreover, $\brho_{\mu}$ is a unique solution to the non-linear Fokker-Planck equation \eqref{pde-main} belonging to $L^\infty_\gamma((0,S],B^{-\beta-\theta}_{p',1})$. This means that for a.e time $s\in(0,S]$, the non-linear drift $\mathcal{B}_{\brho_{\mu}}(s,\cdot)$ is uniquely given by $b(s,\cdot)\ast\brho_{\mu}(s,\cdot)$ and
		\begin{align*}
			\|\mathcal{B}_{\brho_{\mu}}\|_{L^{r_\theta}((0,S],B^{-\theta}_{\infty,\infty})}\leq C\|b\|_{L^r(B^\beta_{p,q})}\|\brho_{\mu}\|_{L^{\infty}_\gamma((0,S],B^{-\beta-\theta}_{p',1})}<+\infty.
		\end{align*}

		Thus, by Theorem 3 of \cite{deraynal2022-linear-besov}, understanding $\mathcal{B}_{\brho_{\mu}}$ as a linear drift and noting that $\theta,r_\theta$ satisfy a good relation \eqref{good-relation}, we deduce that the martingale solution (in the sense of Definition \ref{def-martingale-problem}) related to the non-linear McKean-Vlasov equation \eqref{main-sde} is unique.
	\end{proof}

\end{proposition}

\subsection{Weak well-posedness}

The main result on weak well-posedness of \eqref{main-sde} is as follows.
\begin{proposition}
	\label{prop-weak-wp}
	Let condition \eqref{C3} be satisfied. Let $\theta,r$ satisfy \eqref{good-relation-enhaced}. Then for any initial measure $\mu\in\mathcal{P}(\R^d)\cap B^{\beta_0}_{p_0,q_0}\cap B^{\bar{\theta}}_{1,\infty}$, $\bar{\theta}>\theta$, and any $S\in(0,\mathcal{T})$, the formal McKean-Vlasov SDE \eqref{main-sde} admits a unique weak solution in the sense of Def. \ref{def-weak-solution} such that its time marginal laws $(\bmu^{\mu}_{s})_{s\in[0,S]}$ have a density $\brho_{\mu}(s,x)$ for all $x\in\R^d$ and a.a. $s\in(0,S]$. Moreover, $\brho_{\mu}\in L^\infty_\gamma((0,S],B^{-\beta-\theta}_{p',1})$.

	\begin{proof}
		First note that if $\theta=0$, the existence and uniqueness of a weak solution is straightforward from the well-posedness of the martingale problem. In this case, the Young integral coincides with the usual integrated drift since the drift is a function (see Proposition 7, \cite{deraynal2022-linear-besov}).

		\vspace{2mm}

		If $\theta\in(0,\frac{1}{2})$ and $r$ satisfy condition \eqref{good-relation-enhaced}, the claim follows from Theorem 6 (i), \cite{deraynal2022-linear-besov}.
	\end{proof}
\end{proposition}

\section{Global well-posedness of the McKean-Vlasov SDE}
\label{section-sde-lt}

Recall that Theorem \ref{thm-main-sde} is formulated only in small time, i.e., on the time interval the length of which is at most one. Here, we explain how this result can be extended to any finite time horizon when the initial data $\mu$ is small enough in some Besov norm and when some additional conditions are satisfied. In the usual framework, i.e., when the initial data does not play a major role, the global well-posedness of a system can be proven by reiterating the machinery once it is shown to be well-posed on a small time interval. In our framework, however, it cannot be done because of the strong assumptions on the initial probability measure. The unique solution to the Fokker-Planck equation \eqref{pde-main} may not necessarily possess the same regularity structure as the initial distribution.

\vspace{1mm}

As before, the analysis starts by establishing the appropriate theory for the corresponding PDE satisfied by the mollified density. The main difference is that the estimates on the density $\brho^\varepsilon_{\mu}$ are sought in some time weighted space $L^\infty_{\mathbf{w}}(B^\gamma_{\ell,m})$ where $\mathbf{w}$ is the weight that takes into account both short and long time intervals. This space can be seen as an extension of the space $L^\infty_\gamma(B^\beta_{\ell,m})$ defined in \eqref{def-leb-bes-weight}. Additionally, in the long time we do not have time singularities of the heat kernel, which simplifies the analysis somewhat. Nonetheless, short-time terms appear in analysis and singularities have to be balanced carefully. This leads to the need to modify the previously considered weighted Lebesgue-Besov space. 

\vspace{1mm}

First, let us introduce the tools we will need later on.

\vspace{2mm}

\noindent
\textbf{Global weighted Lebesgue-Besov space.} Let $r,\ell,m\in[1,+\infty]$, $\gamma,\lambda_1,\lambda_2\in\R$, $S\in[0,T]$, $T<+\infty$. Let $w_{\lambda_1}^{\lambda_2}:\R_+\to\R_+$ be the \emph{weight function} given by
\begin{align}
	\label{def-weight-lt}
	w_{\lambda_1}^{\lambda_2}(s):=\big(s\wedge1\big)^{\lambda_1}\big(s\vee1\big)^{\lambda_2}
\end{align}
The index $\lambda_1$ of the weight corresponds to the exponent of the \emph{short time weight} and $\lambda_2$ corresponds to the exponent of the \emph{long time weight}.

\vspace{1mm}

For $r<+\infty$, we define
\begin{align*}
	&L^r_{w_{\lambda_1}^{\lambda_2}}((0,S],B^\gamma_{\ell,m})\\
	&:=\big\{f:s\in[0,S]\mapsto f(s,\cdot)\in B^\gamma_{\ell,m}\text{ measurable and s.t. }\int_0^Sw_{r\lambda_1}^{r\lambda_2}(s)\|f(s,\cdot)\|^r_{B^\gamma_{\ell,m}}ds<+\infty\big\}.
\end{align*}
If $r=+\infty$,
\begin{align*}
	&L^\infty_{w_{\lambda_1}^{\lambda_2}}((0,S],B^\gamma_{\ell,m})\\
	&:=\big\{f:s\in[0,S]\mapsto f(s,\cdot)\in B^\gamma_{\ell,m}\text{ measurable and s.t. }\esssup_{s\in(0,S]}w_{\lambda_1}^{\lambda_2}(s)\|f(s,\cdot)\|_{B^\gamma_{\ell,m}}<+\infty\big\}.
\end{align*}

Note that in short time, we followed the notation $L^r_{\omega}((0,S],B^\gamma_{\ell,m})$ defined in Section \ref{section-tools}, which in the current notation can be understood as $L^r_{w_{\omega}^{0}}((0,S],B^\gamma_{\ell,m})$. We also emphasize that, here, the weight is in the backward time variable, which as before allows us to compensate the singularities around the time origin.

\vspace{2mm}

We also need extended heat kernel estimates that give a more precise bound in long time, see, e.g., \cite{deraynal2023multidimensionalstabledrivenmckeanvlasov} for a proof.

\vspace{2mm}

\noindent
\textbf{Heat kernel Besov norm bound in long time.} For any $\ell,m\in[1,+\infty]$ and $\gamma\in\R$ such that $\gamma\neq-d(1-\frac{1}{\ell})$, there exists a constant $c=c(\alpha,\ell,m,\gamma,d)>0$ s.t. for any multi-index $\mathbf{a}\in\N^d$ with $|\mathbf{a}|\leq1$ and $0<v<s<+\infty$,
\begin{align}
	\label{besov-prop-hk-lt}
	\|\partial^\mathbf{a}p^\alpha_{s-v}\|_{B^\gamma_{\ell,m}}\leq\frac{c}{\big((s-v)\wedge1\big)^{(\frac{\gamma}{\alpha}+\frac{d}{\alpha}(1-\frac{1}{\ell}))_++\frac{|\mathbf{a}|}{\alpha}}\big((s-v)\vee1\big)^{\frac{d}{\alpha}(1-\frac{1}{\ell})+\frac{|\mathbf{a}|}{\alpha}}},
\end{align}
where $\partial^\mathbf{a}=\frac{\partial^{|\mathbf{a}|}}{\partial x_1^{\alpha_1}...\partial x_d^{\alpha_d}}$, $|\mathbf{a}|=a_1+...+a_d$. We see that in small and long times the heat kernel behaves differently, meaning that in long time the regularity index of the norm no longer plays a role and only the integrability index is felt.

\subsection{Auxiliary tools and lemmas}

We state this useful lemma proved in \cite{deraynal2023multidimensionalstabledrivenmckeanvlasov}.
\begin{lemma}[Integration of time singularities, \cite{deraynal2023multidimensionalstabledrivenmckeanvlasov}]
	\label{lemma-beta-function-lt}
	For $0\leq a_1,a_2,b_1,b_2<1$, for $t<s$ in $\R_+$ define
	\begin{align*}
		I_{t,s}:=\int_t^sw^{-a_2}_{-a_1}(s-v)w^{-b_2}_{-b_1}(v-t)dv.
	\end{align*}
	Then there exists $C=C(a,b_1,b_2)>0$ such that for $t<s$ in $\R_+^2$,
	\begin{align*}
		I_{t,s}\leq Cw^{1-b_2-a_2}_{1-b_1-a_1}(s-t).
	\end{align*}
\end{lemma}

\begin{lemma}[General initial data estimate]
	\label{lemma-initial-data-lt}
	For $\rho\in\R$, $p\in(1,+\infty]$ such that $\rho+\frac{d}{p}\neq0$ and $\rho+\frac{d}{p}-\bar{\beta}_0\neq0$, $\mu\in B^{\beta_0}_{p_0,q_0}$ and $v\in(0,T]$, it holds
	\begin{align*}
		\|\mu\ast p^{\alpha}_{v}\|_{B^{\rho}_{p',1}}\leq C\|\mu\|_{B^{\bar{\beta}_0}_{\bar{p}_0,\bar{q}_0}}w^{-\frac{d}{\alpha p}}_{-\gamma_0}(v),
	\end{align*}
	where $\gamma_0=\frac{1}{\alpha}\Big(\rho+\frac{d}{p}-\zeta_0\Big)_+$ and $\zeta_0$ is given by \eqref{def-zeta0}.

	\begin{proof}
		If $v$ is such that $v\leq1$, then the estimation for the initial data is performed as in Lemma \ref{lemma-density-estimate} and results in
		\begin{align*}
			\|\mu\ast p^{\alpha}_{v}\|_{B^{\rho}_{p',1}} & \leq C\|\mu\|_{B^{\bar{\beta}_0}_{\bar{p}_0,\bar{q}_0}}v^{-\gamma_0}.
		\end{align*}

		Then, if $v>1$, one can check that the approach in Lemma \ref{lemma-density-estimate} yields a cruder estimate. Instead, we use the fact that $\mu\in\mathcal{P}(\R^d)$, the Young inequality \eqref{besov-prop-y} and the heat kernel estimate \eqref{besov-prop-hk-lt} to derive
		\begin{align*}
			\|\mu\ast p^{\alpha}_{v}\|_{B^{\rho}_{p',1}} & \leq C\|\mu\|_{B^0_{1,\infty}}\|p^{\alpha}_{v}\|_{B^{\rho}_{p',1}} \\
			                                             & \leq C\|\mu\|_{B^0_{1,\infty}}v^{-\frac{d}{\alpha p}},
		\end{align*}
		when $v>1$. Thus, combining this estimate with \eqref{prf-ineq-initial-data} we get
		\begin{align*}
			\|\mu\ast p^{\alpha}_{v}\|_{B^{\rho}_{p',1}} & \leq C\|\mu\|_{B^{\bar{\beta}_0}_{\bar{p}_0,\bar{q}_0}}\big(v\wedge1\big)^{-\gamma_0}\big(v\vee1\big)^{-\frac{d}{\alpha p}} \\
			                                             & =C\|\mu\|_{B^{\bar{\beta}_0}_{\bar{p}_0,\bar{q}_0}}w^{-\frac{d}{\alpha p}}_{-\gamma_0}(v),
		\end{align*}
		where as before, $\gamma_0=\frac{1}{\alpha}\Big(\rho+\frac{d}{p}-\zeta_0\Big)_+$.
	\end{proof}
\end{lemma}

\subsection{The non-linear Fokker-Planck equation well-posedness and regularity in global time}

In this section, we adapt the result obtained in Section \ref{section-pde} for the space $L^\infty_{w_{\lambda_1}^{\lambda_2}}((0,S],B^\gamma_{\ell,m})$ where $S>1$. We omit the details of proofs that were carefully written in the aforementioned section and focus primarily on the nuances of the long time regime. The main result for the Fokker-Planck equation \eqref{pde-main} in global time reads as follows.
\begin{proposition}
	Let condition \eqref{C3} and \eqref{C3LT} hold. Define
	\begin{align}
		\gamma_1:=\gamma=\frac{\eta}{\alpha}+\frac{1}{\alpha}\Big(-\theta-\beta+\frac{d}{p}-\zeta_0\Big)_+,\quad\gamma_2=\frac{1}{r'}-\frac{1}{\alpha}.
	\end{align}
	Then, for any $\mu\in\mathcal{P}(\R^d)\cap B^{\beta_0}_{p_0,q_0}\cap B^{\bar{\theta}}_{1,\infty}$ for $\bar{\theta}>\theta$, there exists a time horizon $\mathcal{T}_1\in(0,T]$ such that the non-linear Fokker-Planck equation \eqref{pde-main} admits a unique solution $\brho_{\mu}$ belonging to $L^\infty_{w^{\gamma_2}_{\gamma_1}}((0,S],B^{-\beta-\theta}_{p',1})$ for any $S\in(0,\mathcal{T}_1]$.

	\vspace{1mm}
	\noindent
	Moreover, for all $s\in[0,S]$, $\brho_{\mu}(s,\cdot)\in\mathcal{P}(\R^d)$, and for a.e. $s\in(0,S]$, $\brho_{\mu}(s,\cdot)$ is absolutely continuous w.r.t. the Lebesgue measure and satisfies the Duhamel representation \eqref{duhamel-main}.
\end{proposition}

The rest of the section is devoted to proving the above proposition. The proof consists of a chain of Lemmas and has the same structure as Section \ref{section-pde}.

\begin{lemma}[A priori estimates on the mollified density in long time]
	\label{lemma-density-estimate-lt}
	Let $T>1$. Let \eqref{mu-range}, \eqref{C3} and \eqref{C3LT} hold. Then there exists $C>0$ such that for all $S\leq T$,
	\begin{align}
		\label{ineq-quadratic-density-long-time}
		\begin{split}
			 & \sup_{v\in(0,S]}w^{\gamma_2}_{\gamma_1}(v)\|\brho^{\varepsilon}_{\mu}(v,\cdot)\|_{B^{-\beta-\theta}_{p',1}}\leq C\|\mu\|_{B^{\bar{\beta}_0}_{\bar{p}_0,\bar{q}_0}}w^{\frac{1}{r'}-\frac{1}{\alpha}-\frac{d}{\alpha p}}_{\frac{\eta}{\alpha}}(S)                   \\
			 & \quad+C\Big(\|b\|_{L^r(B^\beta_{p,q})}+\|\Div(b)\|_{L^r(B^\beta_{p,q})}\Big)\big(\sup_{v\in(0,S]}w^{\gamma_2}_{\gamma_1}(v)\|\brho^{\varepsilon}_{\mu}(v,\cdot)\|_{B^{-\beta-\theta}_{p',1}}\big)^2w^{0}_{\frac{1}{r'}-\gamma_1+\frac{\beta-\delta}{\alpha}}(S).
		\end{split}
	\end{align}

	\begin{proof}
		As in Lemma \ref{lemma-density-estimate}, for fixed $v\in(0,S]$ write the Duhamel expansion \eqref{duhamel-main-mollified}
		\begin{align}
			\label{prf-3-duhamel-lt}
			\|\brho^{\varepsilon}_{\mu}(v,\cdot)\|_{B^{-\beta-\theta}_{p',1}}\leq\|\mu\ast p^{\alpha}_{ v}\|_{B^{-\beta-\theta}_{p',1}}+\int_0^v\|(\mathcal{B}^{\varepsilon}_{\brho^{\varepsilon}_{\mu}}(s,\cdot)\brho^{\varepsilon}_{\mu}(s,\cdot))\ast\nabla p^\alpha_{v-s}\|_{B^{-\beta-\theta}_{p',1}}ds.
		\end{align}

		\textbf{Initial condition.} Here, we use Lemma \ref{lemma-initial-data-lt} with $\rho=-\beta-\theta$, and we obtain that
		\begin{align*}
			\|\mu\ast p^{\alpha}_{ v}\|_{B^{-\beta-\theta}_{p',1}} & \leq C\|\mu\|_{B^{\bar{\beta}_0}_{\bar{p}_0,\bar{q}_0}}w^{-\frac{d}{\alpha p}}_{-\gamma_0}(v),
		\end{align*}
		where $\gamma_0=\frac{1}{\alpha}\Big(-\theta-\beta+\frac{d}{p}-\zeta_0\Big)_+$.

		\vspace{2mm}

		\textbf{Bilinear operator.} Without loss of generality, we can assume that $ v>1$, otherwise, we are in the framework of Lemma \ref{lemma-density-estimate}. To get a sharp bound for $ v>1$, we split the integral in \eqref{prf-3-duhamel-lt} into two integrals over the intervals $[0,v-1]$ and $[v-1,v]$:
		\begin{align*}
			 & \int_0^v\|(\mathcal{B}^{\varepsilon}_{\brho^{\varepsilon}_{\mu}}(s,\cdot)\brho^{\varepsilon}_{\mu}(s,\cdot))\ast\nabla p^\alpha_{v-s}\|_{B^{-\beta-\theta}_{p',1}}ds                                                                                                                                                                                      \\
			 & =\int_0^{v-1}\|(\mathcal{B}^{\varepsilon}_{\brho^{\varepsilon}_{\mu}}(s,\cdot)\brho^{\varepsilon}_{\mu}(s,\cdot))\ast\nabla p^\alpha_{v-s}\|_{B^{-\beta-\theta}_{p',1}}ds+\int_{v-1}^{v}\|(\mathcal{B}^{\varepsilon}_{\brho^{\varepsilon}_{\mu}}(s,\cdot)\brho^{\varepsilon}_{\mu}(s,\cdot))\ast\nabla p^\alpha_{v-s}\|_{B^{-\beta-\theta}_{p',1}}ds.
		\end{align*}

		The second integral is handled similarly to the one in Lemma \ref{lemma-density-estimate} as we are on a small time interval and have a time singularity at $s=v$. After an integration by parts and the Young inequality \eqref{besov-prop-y} as in \eqref{prf-quad-term}, and the product rule \eqref{ineq-product-rule} as in the controls \eqref{prf-quad-1}, \eqref{prf-quad-2}, we get
		\begin{align*}
			 & \int_{v-1}^{v}\|(\mathcal{B}^{\varepsilon}_{\brho^{\varepsilon}_{\mu}}(s,\cdot)\brho^{\varepsilon}_{\mu}(s,\cdot))\ast\nabla p^\alpha_{v-s}\|_{B^{-\beta-\theta}_{p',1}}ds                                                                     \\
			 & \leq C\int_{v-1}^v\Big(\|b^\varepsilon(s,\cdot)\|_{B^{\beta}_{p,q}}+\|\Div(b^\varepsilon(s,\cdot))\|_{B^{\beta}_{p,q}}\Big)\|\boldsymbol{\rho}^{\varepsilon}_{\mu}(s,\cdot)\|^2_{B^{-\beta-\theta}_{p',1}}(v-s)^{\frac{\beta}{\alpha}}ds.
		\end{align*}

		For handling the first integral, integration by parts is not needed since the heat kernel has no time singularities away from zero. Therefore, similarly to \eqref{prf-quad-1} of Lemma \ref{lemma-density-estimate}, we have
		\begin{align*}
			 & \int_0^{v-1}\|(\mathcal{B}^{\varepsilon}_{\brho^{\varepsilon}_{\mu}}(s,\cdot)\brho^{\varepsilon}_{\mu}(s,\cdot))\ast\nabla p^\alpha_{v-s}\|_{B^{-\beta-\theta}_{p',1}}ds                                                                                                                 \\
			 & \stackrel{\eqref{besov-prop-y}}{\leq} C\int_0^{v-1}\|(\mathcal{B}^{\varepsilon}_{\brho^{\varepsilon}_{\mu}}(s,\cdot)\brho^{\varepsilon}_{\mu}(s,\cdot))\|_{B^{-\theta}_{p',\infty}}\|\nabla p^\alpha_{v-s}\|_{B^{-\beta}_{1,1}}ds                                                                           \\
			 & \stackrel{\eqref{besov-prop-pr2}, \eqref{besov-prop-hk-lt}}{\leq} C\int_0^{v-1}\|\mathcal{B}^{\varepsilon}_{\brho^{\varepsilon}_{\mu}}(s,\cdot)\|_{B^{-\theta}_{\infty,\infty}}\|\brho^{\varepsilon}_{\mu}(s,\cdot)\|_{B^{-\beta-\theta}_{p',\infty}}\big((v-s)\wedge1\big)^{\frac{\beta-1}{\alpha}}\big((v-s)\vee1\big)^{-\frac{1}{\alpha}}ds \\
			 & \stackrel{\eqref{besov-prop-y}}{\leq} C\int_0^{v-1}\|b^\varepsilon(s,\cdot)\|_{B^\beta_{p,q}}\|\brho^{\varepsilon}_{\mu}(s,\cdot)\|^2_{B^{-\beta-\theta}_{p',\infty}}(v-s)^{-\frac{1}{\alpha}}ds,
		\end{align*}
		where in the last inequality we used the fact that $v-s\geq1$. As in Lemma \ref{lemma-density-estimate}, this estimate holds when $\beta>-1-2\theta$. Hence,
		\begin{align*}
			 & \int_0^v\|(\mathcal{B}^{\varepsilon}_{\brho^{\varepsilon}_{\mu}}(s,\cdot)\brho^{\varepsilon}_{\mu}(s,\cdot))\ast\nabla p^\alpha_{v-s}\|_{B^{-\beta-\theta}_{p',1}}ds                                                                                                                                            \\
			 & \leq C\int_{v-1}^v\big(\|b^\varepsilon(s,\cdot)\|_{B^{\beta}_{p,q}}+\|\Div(b^\varepsilon(s,\cdot))\|_{B^{\beta}_{p,q}}\big)\|\boldsymbol{\rho}^{\varepsilon}_{\mu}(s,\cdot)\|^2_{B^{-\beta-\theta}_{p',1}}(v-s)^{\frac{\beta}{\alpha}}ds                                                                   \\
			 & \quad+C\int_0^{v-1}\|b^\varepsilon(s,\cdot)\|_{B^\beta_{p,q}}\|\brho^{\varepsilon}_{\mu}(s,\cdot)\|^2_{B^{-\beta-\theta}_{p',\infty}}(v-s)^{-\frac{1}{\alpha}}ds                                                                                                                                                  \\
			 & \leq C\int_{0}^v\Big(\|b^\varepsilon(s,\cdot)\|_{B^{\beta}_{p,q}}+\|\Div(b^\varepsilon(s,\cdot))\|_{B^{\beta}_{p,q}}\Big)\|\boldsymbol{\rho}^{\varepsilon}_{\mu}(s,\cdot)\|^2_{B^{-\beta-\theta}_{p',1}}\big((v-s)\wedge1\big)^{\frac{\beta}{\alpha}}\big((v-s)\vee1\big)^{-\frac{1}{\alpha}}ds            \\
			 & \leq C\Big(\|b\|_{L^r((0,T],B^{\beta}_{p,q})}+\|\Div(b)\|_{L^r((0,T],B^{\beta}_{p,q})}\Big)\\
			 &\quad\cdot\Big(\int_0^v\|\boldsymbol{\rho}^{\varepsilon}_{\mu}(s,\cdot)\|^{2r'}_{B^{-\beta-\theta}_{p',1}}\big((v-s)\wedge1\big)^{\frac{\beta}{\alpha}r'}\big((v-s)\vee1\big)^{-\frac{r'}{\alpha}}ds\Big)^{\frac{1}{r'}}.
		\end{align*}

		\vspace{2mm}

		\textbf{Compensating the time singularity.} The final estimate is
		\begin{align*}
			 & \|\brho^{\varepsilon}_{\mu}(v,\cdot)\|_{B^{-\beta-\theta}_{p',1}}\leq C\|\mu\|_{B^{\bar{\beta}_0}_{\bar{p}_0,\bar{q}_0}}w^{-\frac{d}{\alpha p}}_{-\gamma_0}(v)                                                                                                                                                   \\
			 & \quad+C\Big(\|b\|_{L^r((0,T],B^{\beta}_{p,q})}+\|\Div(b)\|_{L^r((0,T],B^{\beta}_{p,q})}\Big)\\
			 &\quad\cdot\Big(\int_0^v\|\boldsymbol{\rho}^{\varepsilon}_{\mu}(s,\cdot)\|^{2r'}_{B^{-\beta-\theta}_{p',1}}\big((v-s)\wedge1\big)^{\frac{\beta}{\alpha}r'}\big((v-s)\vee1\big)^{-\frac{r'}{\alpha}}ds\Big)^{\frac{1}{r'}}.
		\end{align*}

		As in Lemma \ref{lemma-density-estimate}, we now multiply the above inequality by the factor $w_{\gamma_2}^{\gamma_1} v$ introduced in \eqref{def-weight-lt} (in comparison with $ v^\gamma$ in the mentioned Lemma), with $\gamma_1,\gamma_2>0$,
		\begin{align*}
			 & w_{\gamma_1}^{\gamma_2}(v)\|\brho^{\varepsilon}_{\mu}(v,\cdot)\|_{B^{-\beta-\theta}_{p',1}}\leq C\|\mu\|_{B^{\bar{\beta}_0}_{\bar{p}_0,\bar{q}_0}}w^{\gamma_2-\frac{d}{\alpha p}}_{\gamma_1-\gamma_0}(v)                                                                                               \\
			 & \quad+Cw_{\gamma_1}^{\gamma_2}(v)\Big(\|b\|_{L^r((0,T],B^{\beta}_{p,q})}+\|\Div(b)\|_{L^r((0,T],B^{\beta}_{p,q})}\Big)\Big(\int_0^v\|\boldsymbol{\rho}^{\varepsilon}_{\mu}(s,\cdot)\|^{2r'}_{B^{-\beta-\theta}_{p',1}}w_{\frac{\beta}{\alpha}r'}^{-\frac{r'}{\alpha}}(v-s)ds\Big)^{\frac{1}{r'}},
		\end{align*}
		and using Lemma \ref{lemma-beta-function-lt},
		\begin{align*}
			&\int_0^v\|\boldsymbol{\rho}^{\varepsilon}_{\mu}(s,\cdot)\|^{2r'}_{B^{-\beta-\theta}_{p',1}}(v-s)^{\frac{\beta}{\alpha}r'}ds \\
			& \leq\big(\sup_{s\in(0,S]}w_{\gamma_1}^{\gamma_2}(s)\|\boldsymbol{\rho}^{\varepsilon}_{\mu}(s,\cdot)\|_{B^{-\beta-\theta}_{p',1}}\big)^{2r'}\int_0^vw_{-2r'\gamma_1}^{-2r'\gamma_2}(s) w_{\frac{\beta}{\alpha}r'}^{-\frac{r'}{\alpha}}(v-s)ds \\
			                                                                                                                                     & \leq\big(\sup_{s\in(0,S]}w_{\gamma_1}^{\gamma_2}(s)\|\boldsymbol{\rho}^{\varepsilon}_{\mu}(s,\cdot)\|_{B^{-\beta-\theta}_{p',1}}\big)^{2r'}w_{1-2r'\gamma_1-\frac{r'}{\alpha}}^{1-2r'\gamma_2+\frac{\beta}{\alpha}r'}(v),
		\end{align*}
		if $\gamma_1,\gamma_2$ are such that
		\begin{align*}
			\gamma_1<\frac{1}{2r'},\quad\gamma_2<\frac{1}{2r'},
		\end{align*}
		where the condition on $\gamma_2$ is needed to invoke Lemma \ref{lemma-beta-function-lt}, and it holds
		\begin{align*}
			\beta>-\alpha\Big(1-\frac{1}{r}\Big),
		\end{align*}
		\begin{align*}
			\alpha\Big(1-\frac{1}{r}\Big)>1\iff\frac{1}{r'}-\frac{1}{\alpha}>0,
		\end{align*}
		to have integrable singularities, where the last inequality follows from \eqref{C3} for $\beta<-1$. It then follows
		\begin{align*}
			 & \sup_{v\in(0,S]}w_{\gamma_1}^{\gamma_2}(v)\|\brho^{\varepsilon}_{\mu}(v,\cdot)\|_{B^{-\beta-\theta}_{p',1}}\leq C\|\mu\|_{B^{\bar{\beta}_0}_{\bar{p}_0,\bar{q}_0}}w^{\gamma_2-\frac{d}{\alpha p}}_{\gamma_1-\gamma_0}(S)                                                                                                             \\
			 & \quad+C\Big(\|b\|_{L^r((0,T],B^{\beta}_{p,q})}+\|\Div(b)\|_{L^r((0,T],B^{\beta}_{p,q})}\Big)\big(\sup_{v\in(0,S]}w_{\gamma_1}^{\gamma_2}(v)\|\boldsymbol{\rho}^{\varepsilon}_{\mu}(v,\cdot)\|_{B^{-\beta-\theta}_{p',1}}\big)^{2}w_{\frac{1}{r'}-\gamma_1+\frac{\beta}{\alpha}}^{\frac{1}{r'}-\gamma_2-\frac{1}{\alpha}}(S).
		\end{align*}

		Contrary to the estimates in Lemma \ref{lemma-density-estimate}, $S$ can be greater than $1$, thus to avoid singularities we want $\gamma_2$ (exponent in large time) to be such that
		\begin{align}
			\label{ineq-gamma2-condition}
			\gamma_2-\frac{d}{\alpha p}\leq0,\quad\frac{1}{r'}-\gamma_2-\frac{1}{\alpha}\leq0,
		\end{align}
		and $\gamma_1$ (exponent in small time) to satisfy
		\begin{align}
			\label{ineq-gamma1-condition}
			\gamma_1-\gamma_0>0,\quad\frac{1}{r'}-\gamma_1+\frac{\beta}{\alpha}>0.
		\end{align}

		For a proper $\gamma_2$ to exist, it should hold
		\begin{align}
			\label{gamma2-cond1}
			\begin{split}
				0<\frac{1}{r'}-\frac{1}{\alpha}\leq\gamma_2\leq\frac{d}{\alpha p}\iff1-\alpha+\frac{\alpha}{r}+\frac{d}{p}\geq0\iff\alpha\Big(1-\frac{1}{r}\Big)\leq1+\frac{d}{p},
			\end{split}
		\end{align}
		and
		\begin{align}
			\label{gamma2-cond2}
			\begin{split}
				0<\frac{1}{r'}-\frac{1}{\alpha}\leq\gamma_2<\frac{1}{2r'}\iff\alpha\Big(1-\frac{1}{r}\Big)<2.
			\end{split}
		\end{align}

		Similarly, for $\gamma_1$, recalling that $\gamma_0=\frac{1}{\alpha}\Big(-\beta-\theta+\frac{d}{p}-\zeta_0\Big)_+$,
		\begin{align*}
			\frac{1}{\alpha}\Big(-\beta-\theta+\frac{d}{p}-\zeta_0\Big)_+<\gamma_1<\frac{1}{r'}+\frac{\beta}{\alpha} \\
			\iff\beta>-\alpha\Big(1-\frac{1}{r}\Big)+\Big(-\theta-\beta+\frac{d}{p}-\zeta_0\Big)_+.
		\end{align*}

		As in Lemma \ref{lemma-density-estimate}, we can take $\gamma_1=\gamma_0+\frac{\eta}{p}$, and then $\beta$ has to satisfy
		\begin{align*}
			\beta>-\alpha\Big(1-\frac{1}{r}\Big)+\Big(-\theta-\beta+\frac{d}{p}-\zeta_0\Big)_+.
		\end{align*}

		Therefore, for the long time, we can choose $\gamma_2=\frac{1}{r'}-\frac{1}{\alpha}$ so that \eqref{gamma2-cond1} and \eqref{gamma2-cond2} are satisfied. With the chosen values of $\gamma_1$ and $\gamma_2$, we exactly obtain \eqref{ineq-quadratic-density-long-time}, which concludes the proof.
	\end{proof}
\end{lemma}

\begin{lemma}[A priori control through a Grönwall type inequality in long time]
	\label{lemma-gronwall-control-lt}
	If \eqref{mu-range}, \eqref{C3} and \eqref{C3LT} hold, then there exist a time $\mathcal{T}\in(0,T)$ and constants $\mathcal{C}_0>0$, $C>0$ such that, setting
	\begin{align*}
		\mathcal{T}_1:=\mathcal{T}\ind_{\{\|\mu\|_{B^{\bar{\beta}_0}_{\bar{p}_0,\bar{q}_0}}\geq \mathcal{C}_0\}}+T\ind_{\{\|\mu\|_{B^{\bar{\beta}_0}_{\bar{p}_0,\bar{q}_0}}<\mathcal{C}_0\}},
	\end{align*}
	we have that for any $S<\mathcal{T}_1$,
	\begin{align*}
		\sup_{v\in(0,S]}w^{\gamma_2}_{\gamma_1}(v)\|\brho^{\varepsilon}_{\mu}(v,\cdot)\|_{B^{-\beta-\theta}_{p',1}}\leq C,
	\end{align*}
	for some $C>0$.

	\begin{proof}
		The proof of Lemma \ref{lemma-gronwall-control} is modified in the following way. Function $f^\varepsilon_{(0,S]}$ is replaced by
		\begin{align*}
			\tilde{f}^\varepsilon_{(0,S]}(s):=\sup_{v\in(0,s]}w^{\gamma_2}_{\gamma_1}(v)\|\brho^{\varepsilon}_{\mu}(v,\cdot)\|_{B^{-\beta-\theta}_{p',1}}, \quad s\in(0,S],
		\end{align*}
		and functions $C_{\mu},C_{b}$ are replaced respectively by
		\begin{align*}
			C_{\mu}(s) & =C\|\mu\|_{B^{\bar{\beta}_0}_{\bar{p}_0,\bar{q}_0}}w^{0}_{\frac{\eta}{\alpha}}(s)=:C_0w^{0}_{\frac{\eta}{\alpha}}
			(s),                                                                                   \\
			C_{b}(s)   & =C\big(\|b\|_{L^r(B^\beta_{p,q})}+\|\Div(b)\|_{L^r(B^\beta_{p,q})}\big)w^{0}_{-\gamma_1+\frac{1}{r'}+\frac{\beta-\delta}{\alpha}}(s)=:C_bw^{0}_{-\gamma_1+\frac{1}{r'}+\frac{\beta-\delta}{\alpha}}(s).
		\end{align*}

		Note that in Lemma \ref{lemma-density-estimate-lt}, we had the bound $C\|\mu\|_{B^{\bar{\beta}_0}_{\bar{p}_0,\bar{q}_0}}w^{\frac{1}{r'}-\frac{1}{\alpha}-\frac{d}{\alpha p}}_{\frac{\eta}{\alpha}}(s) $ for the initial data. But since $\frac{1}{r'}-\frac{1}{\alpha}-\frac{d}{\alpha p}\leq0$ by \eqref{C3LT}, we have $C\|\mu\|_{B^{\bar{\beta}_0}_{\bar{p}_0,\bar{q}_0}}w^{\frac{1}{r'}-\frac{1}{\alpha}-\frac{d}{\alpha p}}_{\frac{\eta}{\alpha}}(s) \leq C_{\mu}(s)$.
		Then following the same argument as in short time, $\tilde{f}^\varepsilon_{(0,S]}$ satisfies
		\begin{align*}
			\tilde{f}^\varepsilon_{(0,S]}(s)\leq\frac{1-\sqrt{1-4C_{\mu}(s)C_{b}(s)}}{2C_{b}(s)},
		\end{align*}
		and
		\begin{align*}
			C_{\mu}(s)C_{b}(s)<\frac{1}{4}\iff w^{0}_{\frac{\eta}{\alpha}-\gamma_1+\frac{1}{r'}+\frac{\beta-\delta}{\alpha}}(s)<(4C_0C_b)^{-1}.
		\end{align*}

		Then it holds for any $s\in(0,\mathcal{T})$, where $\mathcal{T}:=(4C_0C_b)^{1/(\gamma-1+\frac{1}{r}-\frac{\beta-\delta+\eta}{\alpha})}$ (in this case $s\leq1$: small time independently of initial data) or for any $S\in(0,\mathcal{T})$, $s<S$, as soon as $\|\mu\|_{B^{\bar{\beta}_0}_{\bar{p}_0,\bar{q}_0}}<(4CC_b)^{-1}=:\mathcal{C}_0$ (in this case $s>1$: global time under small initial data condition). Under either condition, we obtain the uniform bound
		\begin{align*}
			\sup_{v\in(0,S]}w^{\gamma_2}_{\gamma_1}(v)\|\brho^\varepsilon_{\mu}(v,\cdot)\|_{B^{-\beta-\theta}_{p',1}}\leq\frac{1}{2C_b}(4C_0C_b)^{1-\frac{\eta}{\alpha(-\gamma+1-\frac{1}{r}+\frac{\beta-\delta+\eta}{\alpha})}}.
		\end{align*}
	\end{proof}

\end{lemma}

\begin{lemma}[Convergence of the mollified densities in long time]
	\label{lemma-density-conv-lt}
	Let conditions \eqref{mu-range}, \eqref{C3} and \eqref{C3LT} be satisfied. Then, for all $S\in(0,\mathcal{T}_1)$, with $\mathcal{T}_1$ defined as in Lemma \ref{lemma-gronwall-control-lt}, for any decreasing sequence $(\varepsilon_k)_{k\geq1}$ s.t. $\varepsilon_k\xrightarrow[k\to\infty]{}0$, the sequence $\big(\brho_{\mu}^{\varepsilon_k}\big)_{k\geq1}$ is \emph{Cauchy} in $L^\infty_{w^{\gamma_2}_{\gamma_1}}((0,S],B^{-\beta-\theta}_{p',1})\cap L^\infty((0,S],L^1)$. In particular, there exists $\brho_{\mu}\in L^\infty_{w^{\gamma_2}_{\gamma_1}}((0,S],B^{-\beta-\theta}_{p',1})$ s.t.
	\begin{align}
		\label{conv-mol-dens-lt}
		\sup_{s\in(0,S]}w^{\gamma_2}_{\gamma_1}(s)\|(\brho^{\varepsilon_k}_{\mu}-\brho_{\mu})(s,\cdot)\|_{B^{-\beta-\theta}_{p',1}}+\sup_{s\in(0,S]}\|(\brho^{\varepsilon_k}_{\mu}-\brho_{\mu})(s,\cdot)\|_{L^1}\xrightarrow[k\to+\infty]{}0.
	\end{align}

	\begin{proof}
		Fix $k,j\in\N$ s.t. $k\geq j$. For fixed $s\in(0,S]$ such that $s>1$, $y\in\R^d$, from the Duhamel representation \eqref{duhamel-main-mollified} we write
		\begin{align*}
			\big(\brho^{\varepsilon_k}_{\mu}-\brho^{\varepsilon_j}_{\mu}\big)(s,y) & =-\int_0^s\Big(\mathcal{B}^{\varepsilon_k}_{\brho^{\varepsilon_k}_{\mu}}(v,\cdot)\brho^{\varepsilon_k}_{\mu}(v,\cdot)-\mathcal{B}^{\varepsilon_j}_{\brho^{\varepsilon_j}_{\mu}}(v,\cdot)\brho^{\varepsilon_j}_{\mu}(v,\cdot)\Big)\ast\nabla p^\alpha_{s-v}(y)dv          \\
			                                                                           & =-\int_0^{s-1}\Big(\mathcal{B}^{\varepsilon_k}_{\brho^{\varepsilon_k}_{\mu}}(v,\cdot)\brho^{\varepsilon_k}_{\mu}(v,\cdot)-\mathcal{B}^{\varepsilon_j}_{\brho^{\varepsilon_j}_{\mu}}(v,\cdot)\brho^{\varepsilon_j}_{\mu}(v,\cdot)\Big)\ast\nabla p^\alpha_{s-v}(y)dv      \\
			                                                                           & \quad-\int_{s-1}^s\Big(\mathcal{B}^{\varepsilon_k}_{\brho^{\varepsilon_k}_{\mu}}(v,\cdot)\brho^{\varepsilon_k}_{\mu}(v,\cdot)-\mathcal{B}^{\varepsilon_j}_{\brho^{\varepsilon_j}_{\mu}}(v,\cdot)\brho^{\varepsilon_j}_{\mu}(v,\cdot)\Big)\ast\nabla p^\alpha_{s-v}(y)dv,
		\end{align*}
		and
		\begin{align*}
			\|(\brho^{\varepsilon_k}_{\mu}-\brho^{\varepsilon_j}_{\mu})(s,\cdot)\|_{B^{-\beta-\theta}_{p',1}} & \leq\int_0^{s-1}\Big\|\Big(\mathcal{B}^{\varepsilon_k}_{\brho^{\varepsilon_k}_{\mu}}(v,\cdot)\brho^{\varepsilon_k}_{\mu}(v,\cdot)-\mathcal{B}^{\varepsilon_j}_{\brho^{\varepsilon_j}_{\mu}}(v,\cdot)\brho^{\varepsilon_j}_{\mu}(v,\cdot)\Big)\ast\nabla p^\alpha_{s-v}\Big\|_{B^{-\beta-\theta}_{p',1}}dv      \\
			                                                                                                      & \quad+\int_{s-1}^{s}\Big\|\Big(\mathcal{B}^{\varepsilon_k}_{\brho^{\varepsilon_k}_{\mu}}(v,\cdot)\brho^{\varepsilon_k}_{\mu}(v,\cdot)-\mathcal{B}^{\varepsilon_j}_{\brho^{\varepsilon_j}_{\mu}}(v,\cdot)\brho^{\varepsilon_j}_{\mu}(v,\cdot)\Big)\ast\nabla p^\alpha_{s-v}\Big\|_{B^{-\beta-\theta}_{p',1}}dv.
		\end{align*}

		\textbullet\ To show that $\big(\brho_{\mu}^{\varepsilon_k}\big)_{k\geq1}$ is a Cauchy sequence in $L^\infty_{w^{\gamma_2}_{\gamma_1}}((0,S],B^{-\beta-\theta}_{p',1})$, we handle the second integral (short time integration interval) as in Lemma \ref{lemma-gronwall-control} and we obtain
		\begin{align}
			\label{ineq-cauchy-1-lt}
			\begin{split}
			 & \int_{s-1}^s\Big\|\Big(\mathcal{B}^{\varepsilon_k}_{\brho^{\varepsilon_k}_{\mu}}(v,\cdot)\brho^{\varepsilon_k}_{\mu}(v,\cdot)-\mathcal{B}^{\varepsilon_j}_{\brho^{\varepsilon_j}_{\mu}}(v,\cdot)\brho^{\varepsilon_j}_{\mu}(v,\cdot)\Big)\ast\nabla p^\alpha_{s-v}\Big\|_{B^{-\beta-\theta}_{p',1}}dv                           \\
			 & \leq C\int_0^s\Bigg(\Big(\|b^{\varepsilon_k}(v,\cdot)\|_{B^{\beta}_{p,q}}+\|\Div(b^{\varepsilon_k}(v,\cdot))\|_{B^{\beta}_{p,q}}\Big)\|\big(\brho^{\varepsilon_k}_{\mu}-\brho^{\varepsilon_j}_{\mu}\big)(v,\cdot)\|_{B^{-\beta-\theta}_{p',1}}\|\brho^{\varepsilon_j}_{\mu}(v,\cdot)\|_{B^{-\beta-\theta}_{p',1}}           \\
			 & \quad+\Big(\|\big(b^{\varepsilon_k}-b^{\varepsilon_j}\big)(v,\cdot)\|_{B^{\beta-\delta'}_{p,q}}+\|\Div\big(b^{\varepsilon_k}-b^{\varepsilon_j}\big)(v,\cdot)\|_{B^{\beta-\delta'}_{p,q}}\Big)\|\brho^{\varepsilon_j}_{\mu}(v,\cdot)\|^2_{B^{-\beta-\theta}_{p',1}}\Bigg)\big((s-v)\wedge1\big)^{\frac{\beta-\delta}{\alpha}}dv.
			\end{split}
		\end{align}

		To handle the integral on the time interval $[0,s-1]$, we note that by the heat kernel estimate in long time \eqref{besov-prop-hk-lt} for $s-v\geq1$ the integrand is controlled as follows:
		\begin{align*}
			 & \Big\|\Big(\mathcal{B}^{\varepsilon_k}_{\brho^{\varepsilon_k}_{\mu}}(v,\cdot)\brho^{\varepsilon_k}_{\mu}(v,\cdot)-\mathcal{B}^{\varepsilon_j}_{\brho^{\varepsilon_j}_{\mu}}(v,\cdot)\brho^{\varepsilon_j}_{\mu}(v,\cdot)\Big)\ast\nabla p^\alpha_{s-v}\Big\|_{B^{-\beta-\theta}_{p',1}}                              \\
			 & \leq C\|\mathcal{B}^{\varepsilon_k}_{\brho^{\varepsilon_k}_{\mu}}(v,\cdot)\brho^{\varepsilon_k}_{\mu}(v,\cdot)-\mathcal{B}^{\varepsilon_j}_{\brho^{\varepsilon_j}_{\mu}}(v,\cdot)\brho^{\varepsilon_j}_{\mu}(v,\cdot)\|_{B^{-\beta-\theta}_{p',\infty}}\|\nabla p^\alpha_{s-v}\|_{B^{0}_{1,1}}          \\
			 & \leq C\Big(\|b^{\varepsilon_k}(v,\cdot)\|_{B^{\beta}_{p,q}}\|\big(\brho^{\varepsilon_k}_{\mu}-\brho^{\varepsilon_j}_{\mu}\big)(v,\cdot)\|_{B^{-\beta-\theta}_{p',1}}\|\brho^{\varepsilon_j}_{\mu}(v,\cdot)\|_{B^{-\beta-\theta}_{p',1}}                                                                   \\
			 & \quad+\|\big(b^{\varepsilon_k}-b^{\varepsilon_j}\big)(v,\cdot)\|_{B^{\beta-\delta'}_{p,q}}\|\brho^{\varepsilon_j}_{\mu}(v,\cdot)\|^2_{B^{-\beta-\theta}_{p',1}}\Big)(s-v)^{-\frac{1}{\alpha}}                                                                                                                      \\
			 & \leq C\Bigg(\Big(\|b^{\varepsilon_k}(v,\cdot)\|_{B^{\beta}_{p,q}}+\|\Div(b^{\varepsilon_k}(v,\cdot))\|_{B^{\beta}_{p,q}}\Big)\|\big(\brho^{\varepsilon_k}_{\mu}-\brho^{\varepsilon_j}_{\mu}\big)(v,\cdot)\|_{B^{-\beta-\theta}_{p',1}}\|\brho^{\varepsilon_j}_{\mu}(v,\cdot)\|_{B^{-\beta-\theta}_{p',1}} \\
			 & \quad+\Big(\|\big(b^{\varepsilon_k}-b^{\varepsilon_j}\big)(v,\cdot)\|_{B^{\beta-\delta'}_{p,q}}+\|\Div\big(b^{\varepsilon_k}-b^{\varepsilon_j}\big)(v,\cdot)\|_{B^{\beta-\delta'}_{p,q}}\Big)\|\brho^{\varepsilon_j}_{\mu}(v,\cdot)\|^2_{B^{-\beta-\theta}_{p',1}}\Bigg)(s-v)^{-\frac{1}{\alpha}},
		\end{align*}
		and so
		\begin{align}
			\label{ineq-cauchy-2-lt}
			\begin{split}
			 & \int_0^{s-1}\Big\|\Big(\mathcal{B}^{\varepsilon_k}_{\brho^{\varepsilon_k}_{\mu}}(v,\cdot)\brho^{\varepsilon_k}_{\mu}(v,\cdot)-\mathcal{B}^{\varepsilon_j}_{\brho^{\varepsilon_j}_{\mu}}(v,\cdot)\brho^{\varepsilon_j}_{\mu}(v,\cdot)\Big)\ast\nabla p^\alpha_{s-v}dv\Big\|_{B^{-\beta-\theta}_{p',1}}dv                        \\
			 & \leq C\int_0^{s}\Bigg(\Big(\|b^{\varepsilon_k}(v,\cdot)\|_{B^{\beta}_{p,q}}+\|\Div(b^{\varepsilon_k}(v,\cdot))\|_{B^{\beta}_{p,q}}\Big)\|\big(\brho^{\varepsilon_k}_{\mu}-\brho^{\varepsilon_j}_{\mu}\big)(v,\cdot)\|_{B^{-\beta-\theta}_{p',1}}\|\brho^{\varepsilon_j}_{\mu}(v,\cdot)\|_{B^{-\beta-\theta}_{p',1}} \\
			 & \quad+\Big(\|\big(b^{\varepsilon_k}-b^{\varepsilon_j}\big)(v,\cdot)\|_{B^{\beta-\delta'}_{p,q}}+\|\Div\big(b^{\varepsilon_k}-b^{\varepsilon_j}\big)(v,\cdot)\|_{B^{\beta-\delta'}_{p,q}}\Big)\|\brho^{\varepsilon_j}_{\mu}(v,\cdot)\|^2_{B^{-\beta-\theta}_{p',1}}\Bigg)\big((s-v)\vee1\big)^{-\frac{1}{\alpha}}dv.
			\end{split}
		\end{align}

		Thus, combining \eqref{ineq-cauchy-1-lt} and \eqref{ineq-cauchy-2-lt} and using Lemma \ref{lemma-density-estimate-lt}, we have
		\begin{align*}
			&\|(\brho^{\varepsilon_k}_{\mu}-\brho^{\varepsilon_j}_{\mu})(s,\cdot)\|_{B^{-\beta-\theta}_{p',1}} \\
			&\leq C\int_0^{s}\Bigg(\Big(\|b^{\varepsilon_k}(v,\cdot)\|_{B^{\beta}_{p,q}}+\|\Div(b^{\varepsilon_k}(v,\cdot))\|_{B^{\beta}_{p,q}}\Big)\|\big(\brho^{\varepsilon_k}_{\mu}-\brho^{\varepsilon_j}_{\mu}\big)(v,\cdot)\|_{B^{-\beta-\theta}_{p',1}}w^{-\gamma_2}_{-\gamma_1}(v) \\
			 & \quad+\Big(\|\big(b^{\varepsilon_k}-b^{\varepsilon_j}\big)(v,\cdot)\|_{B^{\beta-\delta'}_{p,q}}+\|\Div\big(b^{\varepsilon_k}-b^{\varepsilon_j}\big)(v,\cdot)\|_{B^{\beta-\delta'}_{p,q}}\Big)w^{-2\gamma_2}_{-2\gamma_1}(v)\Bigg)w^{-\frac{1}{\alpha}}_{\frac{\beta-\delta}{\alpha}}(s-v)dv.
		\end{align*}

		Applying the Hölder inequality with $\bar{r}$ defined in Proposition \ref{prop-b-conv} and $\delta'>0$ as in the proof of Lemma \ref{lemma-density-conv} yields
		\begin{align*}
			&\|(\brho^{\varepsilon_k}_{\mu}-\brho^{\varepsilon_j}_{\mu})(s,\cdot)\|_{B^{-\beta-\theta}_{p',1}} \\
			&\leq C\Big(\|b^{\varepsilon_k}\|_{L^r(B^{\beta}_{p,q})}+\|\Div(b^{\varepsilon_k})\|_{L^r(B^{\beta}_{p,q})}\Big)\Big(\int_0^s\|\big(\brho^{\varepsilon_k}_{\mu}-\brho^{\varepsilon_j}_{\mu}\big)(v,\cdot)\|_{B^{-\beta-\theta}_{p',1}}^{r'}w^{-r'\gamma_2}_{-r'\gamma_1}(v)w^{-\frac{r'}{\alpha}}_{r'\frac{\beta-\delta}{\alpha}}(s-v)dv\Big)^{\frac{1}{r'}}\\
			&\quad+C\Big(\|b^{\varepsilon_k}-b^{\varepsilon_j}\|_{L^{\bar{r}}(B^{\beta-\delta'}_{p,q})}+\|\Div(b^{\varepsilon_k}-b^{\varepsilon_j})\|_{L^{\bar{r}}(B^{\beta-\delta'}_{p,q})}\Big)\Big(\int_0^sw^{-2\bar{r}'\gamma_2}_{-2r'\gamma_1}(v)w^{-\frac{\bar{r}'}{\alpha}}_{\bar{r}'\frac{\beta-\delta}{\alpha}}(s-v)dv\Big)^{\frac{1}{\bar{r}'}}\\
			&\leq C\Big(\|b^{\varepsilon_k}\|_{L^r(B^{\beta}_{p,q})}+\|\Div(b^{\varepsilon_k})\|_{L^r(B^{\beta}_{p,q})}\Big)\Big(\int_0^s\|\big(\brho^{\varepsilon_k}_{\mu}-\brho^{\varepsilon_j}_{\mu}\big)(v,\cdot)\|_{B^{-\beta-\theta}_{p',1}}^{r'}w^{-r'\gamma_2}_{-r'\gamma_1}(v)w^{-\frac{r'}{\alpha}}_{r'\frac{\beta-\delta}{\alpha}}(s-v)dv\Big)^{\frac{1}{r'}}\\
			&\quad+C\Big(\|b^{\varepsilon_k}-b^{\varepsilon_j}\|_{L^{\bar{r}}(B^{\beta-\delta'}_{p,q})}+\|\Div(b^{\varepsilon_k}-b^{\varepsilon_j})\|_{L^{\bar{r}}(B^{\beta-\delta'}_{p,q})}\Big)w^{\frac{1}{\bar{r}'}-2\gamma_2-\frac{1}{\alpha}}_{\frac{1}{\bar{r}'}-2\gamma_1+\frac{\beta-\delta}{\alpha}}(s),
		\end{align*}
		where the last inequality follows from Lemma \ref{lemma-beta-function-lt}. Now, we multiply the r.h.s. by $w^{\gamma_2}_{\gamma_1}(s)$, take the supremum over $[0,S]$ and make the same norm appear in the integrand:
		\begin{align*}
			&\sup_{s\in(0,S]}w^{\gamma_2}_{\gamma_1}(s)\|(\brho^{\varepsilon_k}_{\mu}-\brho^{\varepsilon_j}_{\mu})(s,\cdot)\|_{B^{-\beta-\theta}_{p',1}}\\
			&\leq C\Big(\|b^{\varepsilon_k}\|_{L^r(B^{\beta}_{p,q})}+\|\Div(b^{\varepsilon_k})\|_{L^r(B^{\beta}_{p,q})}\Big)w^{\gamma_2}_{\gamma_1}(S)\\
			&\quad\cdot\Big(\int_0^S\big(\sup_{s\in(0,v]}w^{\gamma_2}_{\gamma_1}(s)\|\big(\brho^{\varepsilon_k}_{\mu}-\brho^{\varepsilon_j}_{\mu}\big)(s,\cdot)\|_{B^{-\beta-\theta}_{p',1}}\big)^{r'}w^{-2r'\gamma_2}_{-2r'\gamma_1}(v)w^{-\frac{r'}{\alpha}}_{r'\frac{\beta-\delta}{\alpha}}(S-v)dv\Big)^{\frac{1}{r'}}\\
			&\quad+C\Big(\|b^{\varepsilon_k}-b^{\varepsilon_j}\|_{L^{\bar{r}}(B^{\beta-\delta'}_{p,q})}+\|\Div(b^{\varepsilon_k}-b^{\varepsilon_j})\|_{L^{\bar{r}}(B^{\beta-\delta'}_{p,q})}\Big)w^{\frac{1}{\bar{r}'}-\gamma_2-\frac{1}{\alpha}}_{\frac{1}{\bar{r}'}-\gamma_1+\frac{\beta-\delta}{\alpha}}(S).
		\end{align*}

		We want to apply the Grönwall-Volterra lemma (\cite[Lemma 2.2]{ZHANGGronwallVolterra}) to the control given by
		\begin{align*}
			&\big(\sup_{s\in(0,S]}w^{\gamma_2}_{\gamma_1}(s)\|(\brho^{\varepsilon_k}_{\mu}-\brho^{\varepsilon_j}_{\mu})(s,\cdot)\|_{B^{-\beta-\theta}_{p',1}}\big)^{r'}\\
			&\leq C\Big(\|b^{\varepsilon_k}\|_{L^r(B^{\beta}_{p,q})}+\|\Div(b^{\varepsilon_k})\|_{L^r(B^{\beta}_{p,q})}\Big)^{r'}w^{r'\gamma_2}_{r'\gamma_1}(S)\\
			&\quad\cdot\int_0^S\big(\sup_{s\in(0,v]}w^{\gamma_2}_{\gamma_1}(s)\|\big(\brho^{\varepsilon_k}_{\mu}-\brho^{\varepsilon_j}_{\mu}\big)(s,\cdot)\|_{B^{-\beta-\theta}_{p',1}}\big)^{r'}w^{-2r'\gamma_2}_{-2r'\gamma_1}(v)w^{-\frac{r'}{\alpha}}_{r'\frac{\beta-\delta}{\alpha}}(S-v)dv\\
			&\quad+C\Big(\|b^{\varepsilon_k}-b^{\varepsilon_j}\|_{L^{\bar{r}}(B^{\beta-\delta'}_{p,q})}+\|\Div(b^{\varepsilon_k}-b^{\varepsilon_j})\|_{L^{\bar{r}}(B^{\beta-\delta'}_{p,q})}\Big)^{r'}w^{\frac{r'}{\bar{r}'}-r'\gamma_2-\frac{r'}{\alpha}}_{\frac{r'}{\bar{r}'}-r'\gamma_1+r'\frac{\beta-\delta}{\alpha}}(S).
		\end{align*}
		
		For this, we need to check is the coefficients of this inequality satisfy the required conditions dictated by Lemma 2.1 and Lemma 2.2 in the mentioned article. Note that our so-called kernel
		\begin{align*}
			\kappa(S,v):=C_{b^{\varepsilon_k}}w^{r'\gamma_2}_{r'\gamma_1}(S)w^{-2r'\gamma_2}_{-2r'\gamma_1}(v)w^{-\frac{r'}{\alpha}}_{r'\frac{\beta-\delta}{\alpha}}(S-v)
		\end{align*}
		with $C_{b^{\varepsilon_k}}=C\Big(\|b^{\varepsilon_k}\|_{L^r(B^{\beta}_{p,q})}+\|\Div(b^{\varepsilon_k})\|_{L^r(B^{\beta}_{p,q})}\Big)^{r'}$, is very close to that defined in \cite[Example 2.4]{ZHANGGronwallVolterra}. We check that this kernel satisfies the conditions of Lemma 2.1. Firstly, from Lemma \ref{lemma-beta-function-lt} we derive
		\begin{align*}
			w^{r'\gamma_2}_{r'\gamma_1}(S)\int_0^Sw^{-2r'\gamma_2}_{-2r'\gamma_1}(v)w^{-\frac{r'}{\alpha}}_{r'\frac{\beta-\delta}{\alpha}}(S-v)dv=w^{1-r'\gamma_2-\frac{r'}{\alpha}}_{1-r'\gamma_1+r'\frac{\beta-\delta}{\alpha}}(S).
		\end{align*}

		Recalling that $1-r'\gamma_2-\frac{r'}{\alpha}\leq0$ and $1-r'\gamma_1+r'\frac{\beta-\delta}{\alpha}>0$ which is guaranteed by \eqref{ineq-gamma2-condition} and \eqref{ineq-gamma1-condition} respectively, we obtain that
		\begin{align*}
			S\mapsto\int_0^S\kappa(S,v)dv\in L^\infty(0,T).
		\end{align*}

		Secondly, we verify that
		\begin{align*}
			\limsup_{\varepsilon\to0}\Big\|\int_{\cdot}^{\cdot+\varepsilon}\kappa(\cdot+\varepsilon,v)dv\Big\|_{L^\infty(0,T)}<1.
		\end{align*}

		We have
		\begin{align*}
			\int_{s}^{s+\varepsilon}\kappa(s+\varepsilon,v)dv&=C_{b^{\varepsilon_k}}w^{r'\gamma_2}_{r'\gamma_1}(s+\varepsilon)\int_{s}^{s+\varepsilon}w^{-2r'\gamma_2}_{-2r'\gamma_1}(v)w^{-\frac{r'}{\alpha}}_{r'\frac{\beta-\delta}{\alpha}}(s+\varepsilon-v)dv\\
			&=C_{b^{\varepsilon_k}}w^{1-r'\gamma_2-\frac{r'}{\alpha}}_{1-r'\gamma_1+r'\frac{\beta-\delta}{\alpha}}(s+\varepsilon)\int_{\frac{s}{s+\varepsilon}}^{1}w^{-2r'\gamma_2}_{-2r'\gamma_1}(v)w^{-\frac{r'}{\alpha}}_{r'\frac{\beta-\delta}{\alpha}}(1-v)dv.
		\end{align*}

		Then it is clear that
		\begin{align*}
			\limsup_{\varepsilon\to0}\int_{1-\varepsilon}^{1}w^{-2r'\gamma_2}_{-2r'\gamma_1}(v)w^{-\frac{r'}{\alpha}}_{r'\frac{\beta-\delta}{\alpha}}(1-v)dv\Big\|w^{1-r'\gamma_2-\frac{r'}{\alpha}}_{1-r'\gamma_1+r'\frac{\beta-\delta}{\alpha}}(\cdot+\varepsilon)\Big\|_{L^\infty(0,T)}=0,
		\end{align*}
		and the kernel $\kappa$ satisfies \cite[Lemma 2.1]{ZHANGGronwallVolterra}. It implies that defining
		\begin{align*}
			r_1(S,v):=\kappa(S,v),\quad r_n(S,v):=\int_v^S\kappa(S,u)r_{n+1}(u,v)du,\quad r(S,v):=\sum_{n=1}^{+\infty}r_n(S,v),
		\end{align*}
		the series $r(S,v)$ is convergent and is integrable in $v$ over $[0,S]$. Moreover,
		\begin{align*}
			C\Big(\|b^{\varepsilon_k}-b^{\varepsilon_j}\|_{L^{\bar{r}}(B^{\beta}_{p,q})}+\|\Div(b^{\varepsilon_k}-b^{\varepsilon_j})\|_{L^{\bar{r}}(B^{\beta}_{p,q})}\Big)^{r'}\int_0^Sr(S,v)w^{\frac{r'}{\bar{r}'}-r'\gamma_2-\frac{r'}{\alpha}}_{\frac{r'}{\bar{r}'}-r'\gamma_1+r'\frac{\beta-\delta}{\alpha}}(v)dv<+\infty
		\end{align*}
		and
		\begin{align*}
			S\mapsto\int_0^Sr_1(S,v)\big(\sup_{s\in(0,v]}w^{\gamma_2}_{\gamma_1}(s)\|(\brho^{\varepsilon_k}_{\mu}-\brho^{\varepsilon_j}_{\mu})(s,\cdot)\|_{B^{-\beta-\theta}_{p',1}}\big)^{r'}dv\in L^\infty(0,T).
		\end{align*}

		All these properties allow us to conclude that 
		\begin{align*}
			&\big(\sup_{s\in(0,S]}w^{\gamma_2}_{\gamma_1}(s)\|(\brho^{\varepsilon_k}_{\mu}-\brho^{\varepsilon_j}_{\mu})(s,\cdot)\|_{B^{-\beta-\theta}_{p',1}}\big)^{r'}\\
			&\leq C\Big(\|b^{\varepsilon_k}-b^{\varepsilon_j}\|_{L^{\bar{r}}(B^{\beta}_{p,q})}+\|\Div(b^{\varepsilon_k}-b^{\varepsilon_j})\|_{L^{\bar{r}}(B^{\beta}_{p,q})}\Big)^{r'}\\
			&\quad\cdot\Big(w^{\frac{r'}{\bar{r}'}-r'\gamma_2-\frac{r'}{\alpha}}_{\frac{r'}{\bar{r}'}-r'\gamma_1+r'\frac{\beta-\delta}{\alpha}}(S)+\int_0^Sr(S,v)w^{\frac{r'}{\bar{r}'}-r'\gamma_2-\frac{r'}{\alpha}}_{\frac{r'}{\bar{r}'}-r'\gamma_1+r'\frac{\beta-\delta}{\alpha}}(v)dv\Big)\\
			&\xrightarrow[j,k\to+\infty]{\text{Proposition }\ref{prop-b-conv}}0.
		\end{align*}

		In other words, $\big(\brho^{\varepsilon_k}_{\mu}\big)_{k\geq1}$ is Cauchy in $L^{\infty}_{w^{\gamma_2}_{\gamma_1}}((0,S],B^{-\beta-\theta}_{p',1})$.

		\vspace{2mm}

		\textbullet\ For the second part, i.e., for showing that $\big(\brho^{\varepsilon_k}_{\mu}\big)_{k\geq1}$ is Cauchy in $L^\infty((0,S],L^1)$, we again mimic the approach used in Lemma \ref{lemma-density-conv} taking into account the long time factor for the heat kernel \eqref{besov-prop-hk-lt}, Lemma \ref{lemma-density-estimate-lt}:
		\begin{align*}
			& \int_0^s\Big\|\Big(\mathcal{B}^{\varepsilon_k}_{\brho^{\varepsilon_k}_{\mu}}(v,\cdot)\brho^{\varepsilon_k}_{\mu}(v,\cdot)-\mathcal{B}^{\varepsilon_j}_{\brho^{\varepsilon_j}_{\mu}}(v,\cdot)\brho^{\varepsilon_j}_{\mu}(v,\cdot)\Big)\ast\nabla p^\alpha_{s-v}\Big\|_{B^{\theta+\delta}_{1,1}}dv                                                                                                            \\
			& \leq C\int_0^s\Big(\|b^{\varepsilon_k}(v,\cdot)\|_{B^{\beta}_{p,q}}w^{-\gamma_2}_{-\gamma_1}(v)\|\big(\brho^{\varepsilon_k}_{\mu}-\brho^{\varepsilon_j}_{\mu}\big)(v,\cdot)\|_{B^{\theta+\delta}_{1,1}}                                                                                                                 \\
			 &\quad+\|b^{\varepsilon_k}(v,\cdot)\|_{B^{\beta}_{p,q}}\|\big(\brho^{\varepsilon_k}_{\mu}-\brho^{\varepsilon_j}_{\mu}\big)(v,\cdot)\|_{B^{-\beta-\theta}_{p',1}}\|\brho^{\varepsilon_j}_{\mu}(v,\cdot)\|_{B^{\theta+\delta}_{1,1}}\\
			 & \quad+\|\big(b^{\varepsilon_k}-b^{\varepsilon_j}\big)(v,\cdot)\|_{B^{\beta-\delta'}_{p,q}}\|\brho^{\varepsilon_j}_{\mu}(v,\cdot)\|_{B^{\theta+\delta}_{1,1}}w^{-\gamma_2}_{-\gamma_1}(v)\Big)w^{-\frac{1}{\alpha}}_{-\frac{1+2\theta+2\delta}{\alpha}}(s-v)dv.
		\end{align*}

		Hölder's inequality in time yields
		\begin{align*}
			& \int_0^s\Big\|\Big(\mathcal{B}^{\varepsilon_k}_{\brho^{\varepsilon_k}_{\mu}}(v,\cdot)\brho^{\varepsilon_k}_{\mu}(v,\cdot)-\mathcal{B}^{\varepsilon_j}_{\brho^{\varepsilon_j}_{\mu}}(v,\cdot)\brho^{\varepsilon_j}_{\mu}(v,\cdot)\Big)\ast\nabla p^\alpha_{s-v}\Big\|_{B^{\theta+\delta}_{1,1}}dv                                                                                                            \\
			& \leq C\|b^{\varepsilon_k}\|_{L^{r}(B^{\beta}_{p,q})}\Big(\int_0^s\|\big(\brho^{\varepsilon_k}_{\mu}-\brho^{\varepsilon_j}_{\mu}\big)(v,\cdot)\|_{B^{\theta+\delta}_{1,1}}^{r'}w^{-r'\gamma_2}_{-r'\gamma_1}(v)w^{-\frac{r'}{\alpha}}_{-r'\frac{1+2\theta+2\delta}{\alpha}}(s-v)dv\Big)^{\frac{1}{r'}}                         \\
			 &\quad+C\|b^{\varepsilon_k}\|_{L^{r}(B^{\beta}_{p,q})}\Big(\int_0^s\|\brho^{\varepsilon_j}_{\mu}(v,\cdot)\|_{B^{\theta+\delta}_{1,1}}^{r'}\|\big(\brho^{\varepsilon_k}_{\mu}-\brho^{\varepsilon_j}_{\mu}\big)(v,\cdot)\|_{B^{-\beta-\theta}_{p',1}}^{r'}w^{-\frac{r'}{\alpha}}_{-r'\frac{1+2\theta+2\delta}{\alpha}}(s-v)dv\Big)^{\frac{1}{r'}}\\
			 & \quad+C\|b^{\varepsilon_k}-b^{\varepsilon_j}\|_{L^{\bar{r}}(B^{\beta-\delta'}_{p,q})}\Big(\int_0^s\|\brho^{\varepsilon_j}_{\mu}(v,\cdot)\|_{B^{\theta+\delta}_{1,1}}^{\bar{r}'}w^{-\bar{r}'\gamma_2}_{-\bar{r}'\gamma_1}(v)w^{-\frac{\bar{r}'}{\alpha}}_{-\bar{r}'\frac{1+2\theta+2\delta}{\alpha}}(s-v)dv\Big)^{\frac{1}{\bar{r}'}},
		\end{align*}
		and after computing the time integrals, we have
		\begin{align}
			\label{ineq-l1-cauchy-lt}
			\begin{split}
			& \|\big(\brho^{\varepsilon_k}_{\mu}-\brho^{\varepsilon_j}_{\mu}\big)(s,\cdot)\|_{B^{\theta+\delta}_{1,1}}\\
			&\leq C\|b^{\varepsilon_k}\|_{L^{r}(B^{\beta}_{p,q})}\Big(\int_0^s\|\big(\brho^{\varepsilon_k}_{\mu}-\brho^{\varepsilon_j}_{\mu}\big)(v,\cdot)\|_{B^{\theta+\delta}_{1,1}}^{r'}w^{-r'\gamma_2}_{-r'\gamma_1}(v)w^{-\frac{r'}{\alpha}}_{-r'\frac{1+2\theta+2\delta}{\alpha}}(s-v)dv\Big)^{\frac{1}{r'}}                         \\
			 &\quad+C\|b^{\varepsilon_k}\|_{L^{r}(B^{\beta}_{p,q})}\sup_{v\in(0,s]}\|\brho^{\varepsilon_j}_{\mu}(v,\cdot)\|_{B^{\theta+\delta}_{1,1}}\|\brho^{\varepsilon_k}_{\mu}-\brho^{\varepsilon_j}_{\mu}\|_{L^{\infty}_{w^{\gamma_2}_{\gamma_1}}((0,s],B^{-\beta-\theta}_{p',1})}w^{\frac{1}{r'}-\gamma_2-\frac{r'}{\alpha}}_{\frac{1}{r'}-\gamma_1-\frac{1+2\theta+2\delta}{\alpha}}(s)\\
			 & \quad+C\|b^{\varepsilon_k}-b^{\varepsilon_j}\|_{L^{\bar{r}}(B^{\beta-\delta'}_{p,q})}\sup_{v\in(0,s]}\|\brho^{\varepsilon_j}_{\mu}(v,\cdot)\|_{B^{\theta+\delta}_{1,1}}w^{\frac{1}{\bar{r}'}-\gamma_2-\frac{1}{\alpha}}_{\frac{1}{\bar{r}'}-\gamma_1-\frac{1+2\theta+2\delta}{\alpha}}(s).
			\end{split}
		\end{align}

		We bound the norm $\|\brho^{\varepsilon_j}_{\mu}(v,\cdot)\|_{B^{\theta+\delta}_{1,1}}$ similarly as in the proof of Lemma \ref{lemma-density-conv} using \eqref{besov-prop-hk-lt} as the heat kernel estimate on the long time interval:
		\begin{align*}
			&\|\brho^{\varepsilon_j}_{\mu}(v,\cdot)\|_{B^{\theta+\delta}_{1,1}}\\
			&\leq\|\mu\ast p^\alpha_{ v}\|_{B^{\theta+\delta}_{1,1}}+\int_0^v\|(\mathcal{B}^{\varepsilon}_{\brho^{\varepsilon}_{\mu}}(t,\cdot)\brho^{\varepsilon}_{\mu}(t,\cdot))\ast\nabla p^\alpha_{v}\|_{B^{\theta+\delta}_{1,1}}dt\\
			&\leq C\|\mu\|_{B^{\theta+\delta}_{1,\infty}}+\int_0^{v-1}\|(\mathcal{B}^{\varepsilon}_{\brho^{\varepsilon}_{\mu}}(t,\cdot)\brho^{\varepsilon}_{\mu}(t,\cdot))\ast\nabla p^\alpha_{v}\|_{B^{\theta+\delta}_{1,1}}dt\\
			 &\quad+\int_{v-1}^{v}\|(\mathcal{B}^{\varepsilon}_{\brho^{\varepsilon}_{\mu}}(t,\cdot)\brho^{\varepsilon}_{\mu}(t,\cdot))\ast\nabla p^\alpha_{v}\|_{B^{\theta+\delta}_{1,1}}dt\\
			&\leq C\|\mu\|_{B^{\theta+\delta}_{1,\infty}}+C\int_0^{v-1}\|b^\varepsilon(t,\cdot)\|_{B^\beta_{p,q}}\|\brho^{\varepsilon}_{\mu}(t,\cdot)\|_{B^{-\beta-\theta}_{p',1}}\|\brho^{\varepsilon}_{\mu}(t,\cdot)\|_{B^{\theta+\delta}_{1,1}}(v-t)^{-\frac{1}{\alpha}}dt\\
			&\quad+C\int_{v-1}^v\|b^\varepsilon(t,\cdot)\|_{B^\beta_{p,q}}\|\brho^{\varepsilon}_{\mu}(t,\cdot)\|_{B^{-\beta-\theta}_{p',1}}(v-t)^{-\frac{1+2\theta+\delta}{\alpha}}dt\\
			&\leq C\|\mu\|_{B^{\theta+\delta}_{1,\infty}}+C\|b^\varepsilon\|_{L^r(B^\beta_{p,q})}\|\brho^{\varepsilon}_{\mu}\|_{L^\infty_{w^{\gamma_2}_{\gamma_1}}((0,v],B^{-\beta-\theta}_{p',1})}\\
			&\quad \cdot\Big(\int_0^v\|\brho^{\varepsilon}_{\mu}(t,\cdot)\|_{B^{\theta+\delta}_{1,1}}^{r'}w^{-r'\gamma_2}_{-r'\gamma_1}(t)w^{-\frac{r'}{\alpha}}_{-r'\frac{1+2\theta+\delta}{\alpha}}(v-t)dt\Big)^{\frac{1}{r'}}.
		\end{align*}

		Then we apply the Grönwall-Volterra lemma from \cite{ZHANGGronwallVolterra}, and we obtain that there exists $\bar{C}>0$ such that
		\begin{align*}
			\|\brho^{\varepsilon_j}_{\mu}(v,\cdot)\|_{B^{\theta+\delta}_{1,1}}\leq \bar{C}.
		\end{align*}

		Substituting this control back into \eqref{ineq-l1-cauchy-lt}, we end up with
		\begin{align*}
			& \|\big(\brho^{\varepsilon_k}_{\mu}-\brho^{\varepsilon_j}_{\mu}\big)(s,\cdot)\|_{B^{\theta+\delta}_{1,1}}\\
			&\leq C\|b^{\varepsilon_k}\|_{L^{r}(B^{\beta}_{p,q})}\Big(\int_0^s\|\big(\brho^{\varepsilon_k}_{\mu}-\brho^{\varepsilon_j}_{\mu}\big)(v,\cdot)\|_{B^{\theta+\delta}_{1,1}}^{r'}w^{-r'\gamma_2}_{-r'\gamma_1}(v)w^{-\frac{r'}{\alpha}}_{-r'\frac{1+2\theta+2\delta}{\alpha}}(s-v)dv\Big)^{\frac{1}{r'}}                         \\
			 &\quad+C\bar{C}\|b^{\varepsilon_k}\|_{L^{r}(B^{\beta}_{p,q})}\|\brho^{\varepsilon_k}_{\mu}-\brho^{\varepsilon_j}_{\mu}\|_{L^{\infty}_{w^{\gamma_2}_{\gamma_1}}((0,s],B^{-\beta-\theta}_{p',1})}w^{\frac{1}{r'}-\gamma_2-\frac{r'}{\alpha}}_{\frac{1}{r'}-\gamma_1-\frac{1+2\theta+2\delta}{\alpha}}(s)\\
			 & \quad+C\bar{C}\|b^{\varepsilon_k}-b^{\varepsilon_j}\|_{L^{\bar{r}}(B^{\beta-\delta'}_{p,q})}w^{\frac{1}{\bar{r}'}-\gamma_2-\frac{1}{\alpha}}_{\frac{1}{\bar{r}'}-\gamma_1-\frac{1+2\theta+2\delta}{\alpha}}(s).
		\end{align*}

		One can check that we are again in the setting of \cite{ZHANGGronwallVolterra} to apply the Grönwall-Volterra lemma, and therefore, there exists a constant $C(k,j,S)>0$ such that
		\begin{align*}
			 & \sup_{s\in(0,S]}\|(\brho^{\varepsilon_k}_{\mu}-\brho^{\varepsilon_j}_{\mu})(s,\cdot)\|_{B^{\theta+\delta}_{1,1}}\leq C(k,j,S)\xrightarrow[k,j\to+\infty]{}0.
		\end{align*}

		The proof is thus complete.

	\end{proof}

\end{lemma}

\begin{lemma}[Existence of a solution for the limit Fokker-Planck equation in long time]
	\label{lemma-density-existence-lt}
	Let conditions \eqref{mu-range}, \eqref{C3} and \eqref{C3LT} be satisfied. Let $(\varepsilon_k)_{k\geq1}$ be a decreasing sequence and $\brho_{\mu}$ be the limit point of $(\brho^{\varepsilon_k}_{\mu})_{k\geq1}$, which exists by Lemma \ref{lemma-density-conv}. Then $\brho_{\mu}$ satisfies the Fokker-Planck equation \eqref{pde-main} in a distributional sense and it also satisfies the Duhamel representation \eqref{duhamel-main}.

	\begin{proof}
		The modifications we need to make in Lemma \ref{lemma-density-existence} concerns only the functional $\Delta^2_{\brho_{\mu},\brho_{\mu}^{\varepsilon_k}}$. For $\rho>\theta$, $\varepsilon>0$, we perform the same splitting as in the short time regime, and we apply the Besov duality \eqref{besov-prop-d}:
		\begin{align*}
			 & |\Delta^2_{\brho_{\mu},\brho_{\mu}^{\varepsilon_k}}(\phi)|                                                                                                                                                                                                                                                                                                                                                                           \\
			 & \leq \Big|\int_0^S\int_{\R^d}\big(\brho_{\mu}-\brho_{\mu}^{\varepsilon_k}\big)(s,x)\big(\mathcal{B}_{\brho_{\mu}}\cdot\nabla\phi\big)(s,x)dsdx\Big|\\
			 &\quad+\Big|\int_0^S\int_{\R^d}\brho_{\mu}^{\varepsilon_k}(s,x)\big((\mathcal{B}_{\brho_{\mu}}-\mathcal{B}^{\varepsilon_k}_{\brho^{\varepsilon_k}_{\mu}})\cdot\nabla\phi\big)(s,x)dsdx\Big|                                                                                                              \\
			 & \leq\int_0^S\|(\brho_{\mu}-\brho_{\mu}^{\varepsilon_k})(s,\cdot)\|_{B^{\theta}_{1,1}}\|\big(\mathcal{B}_{\brho_{\mu}}\cdot\nabla\phi\big)(s,\cdot)\|_{B^{-\theta}_{\infty,\infty}}ds\\
			 &\quad+\int_0^S\|\brho_{\mu}^{\varepsilon_k}(s,\cdot)\|_{B^{\theta+\varepsilon}_{1,1}}\|\big((\mathcal{B}_{\brho_{\mu}}-\mathcal{B}^{\varepsilon_k}_{\brho^{\varepsilon_k}_{\mu}})\cdot\nabla\phi\big)(s,\cdot)\|_{B^{-\theta-\varepsilon}_{\infty,\infty}}ds.
		\end{align*}

		Then, using Lemma \ref{lemma-density-estimate-lt}, the product rule \ref{besov-prop-pr2} and the Young inequality \eqref{besov-prop-y}, we obtain
		\begin{align*}
			& |\Delta^2_{\brho_{\mu},\brho_{\mu}^{\varepsilon_k}}(\phi)|                                                                                                                                                                                                                                                                                                                                                                           \\
			&\leq C\|\nabla\phi\|_{B^\rho_{\infty,\infty}}\Bigg(\int_0^Sw^{-\gamma_2}_{-\gamma_1}(s)\|b(s,\cdot)\|_{B^{\beta}_{p,q}}\|\big(\brho_{\mu}-\brho_{\mu}^{\varepsilon_k}\big)(s,\cdot)\|_{B^{-\beta-\theta}_{p',1}}ds                                                                                                                                                                                                                               \\
			 & \quad+\int_0^Sw^{-\gamma_2}_{-\gamma_1}(s) \Big(\|(b-b^{\varepsilon_k})(s,\cdot)\|_{B^{\beta-\varepsilon}_{p,q}}\|\brho^{\varepsilon_k}_{\mu}(s,\cdot)\|_{B^{-\beta-\theta}_{p',1}}+\|b(s,\cdot)\|_{B^{\beta}_{p,q}}\|\big(\brho_{\mu}-\brho_{\mu}^{\varepsilon_k}\big)(s,\cdot)\|_{B^{-\beta-\theta}_{p',1}}\Big)ds\Bigg).
		\end{align*}

		Finally, by the Hölder inequality and Lemma \ref{lemma-beta-function-lt}, we have
		\begin{align*}
			& |\Delta^2_{\brho_{\mu},\brho_{\mu}^{\varepsilon_k}}(\phi)|                                                                                                                                                                                                                                                                                                                                                                           \\
			&\leq C\|\nabla\phi\|_{B^\rho_{\infty,\infty}}\Bigg(\|b\|_{L^r(B^{\beta}_{p,q})}\|\brho_{\mu}-\brho_{\mu}^{\varepsilon_k}\|_{L^\infty_{w^{\gamma_2}_{\gamma_1}}((0,S],B^{-\beta-\theta}_{p',1})}w^{\frac{1}{r'}-2\gamma_2}_{\frac{1}{r'}-2\gamma_1}(S)+\|b-b^{\varepsilon_k}\|_{L^{\bar{r}}(B^{\beta-\varepsilon}_{p,q})}w^{\frac{1}{\bar{r}'}-2\gamma_2}_{\frac{1}{\bar{r}'}-2\gamma_1}(S)\Bigg).
		\end{align*}

		From Proposition \ref{prop-b-conv} and Lemmas \ref{lemma-density-conv}, \ref{lemma-density-conv-lt} we see that $|\Delta^2_{\brho_{\mu},\brho_{\mu}^{\varepsilon_k}}(\phi)|\xrightarrow[k\to\infty]{}0$.

	\end{proof}

\end{lemma}

\begin{lemma}[Uniqueness of the limit Fokker-Planck solution in long time]
	\label{lemma-density-uniqueness-lt}
	Under conditions \eqref{mu-range}, \eqref{C3} and \eqref{C3LT}, the Fokker-Planck equation \eqref{pde-main} admits a unique solution in $L^\infty_{w^{\gamma_2}_{\gamma_1}}((0,S],B^{-\beta-\theta}_{p',1})$ for any $S<\mathcal{T}_1$, where $\mathcal{T}_1$ is defined in Lemma \ref{lemma-gronwall-control-lt}.

	\begin{proof}
		Let $\brho^1_{\mu}$ and $\brho^2_{\mu}$ be two possible solutions to \eqref{pde-main}. Repeating the argument from Lemma \ref{lemma-density-uniqueness} and the argument typical for long time, we get for $s\in(0,S]$
		\begin{align*}
			 & \|(\brho^1_{\mu}-\brho^2_{\mu})(s,\cdot)\|_{B^{-\beta-\theta}_{p',1}}\\
			 &\leq\int_0^{s-1}\|\big(\mathcal{B}_{\brho^1_{\mu}}(v,\cdot)\brho^1_{\mu}(v,\cdot)-\mathcal{B}_{\brho^2_{\mu}}(v,\cdot)\brho^2_{\mu}(v,\cdot)\big)\ast\nabla p^\alpha_{s-v}\|_{B^{-\beta-\theta}_{p',1}}dv                                                          \\
			 & \quad+\int_{s-1}^s\|\big(\mathcal{B}_{\brho^1_{\mu}}(v,\cdot)\brho^1_{\mu}(v,\cdot)-\mathcal{B}_{\brho^2_{\mu}}(v,\cdot)\brho^2_{\mu}(v,\cdot)\big)\ast\nabla p^\alpha_{s-v}\|_{B^{-\beta-\theta}_{p',1}}dv                                                                                                                                 \\
			 & \leq C\int_0^{s-1}\Big(\|\brho^1_{\mu}(v,\cdot)\|_{B^{-\beta-\theta}_{p',1}}+\|\brho^2_{\mu}(v,\cdot)\|_{B^{-\beta-\theta}_{p',1}}\Big)\|b(v,\cdot)\|_{B^{\beta}_{p,q}}\|(\brho^1_{\mu}-\brho^2_{\mu})(v,\cdot)\|_{B^{-\beta-\theta}_{p',1}}\|\nabla p^\alpha_{s-v}\|_{B^{-\beta}_{1,1}}dv                                           \\
			 & \quad+\int_{s-1}^s\Big(\|\brho^1_{\mu}(v,\cdot)\|_{B^{-\beta-\theta}_{p',1}}+\|\brho^2_{\mu}(v,\cdot)\|_{B^{-\beta-\theta}_{p',1}}\Big)\Big(\|b(v,\cdot)\|_{B^{\beta}_{p,q}}+\|\Div(b(v,\cdot))\|_{B^{\beta}_{p,q}}\Big)\\
			 &\quad\cdot\|(\brho^1_{\mu}-\brho^2_{\mu})(v,\cdot)\|_{B^{-\beta-\theta}_{p',1}}\|p^\alpha_{s-v}\|_{B^{-\beta}_{1,1}}dv.
		\end{align*}

		Using Lemma \ref{lemma-density-estimate-lt}, \eqref{besov-prop-hk-lt} and the Hölder inequality, we get
		\begin{align*}
			& \|(\brho^1_{\mu}-\brho^2_{\mu})(s,\cdot)\|_{B^{-\beta-\theta}_{p',1}}\\
			&\leq C\Big(\|b\|_{L^r(B^{\beta}_{p,q})}+\|\Div(b)\|_{L^r(B^{\beta}_{p,q})}\Big)\Big(\int_0^s\|(\brho^1_{\mu}-\brho^2_{\mu})(v,\cdot)\|_{B^{-\beta-\theta}_{p',1}}^{r'}w^{-r'\gamma_2}_{-r'\gamma_1}(v)w^{-\frac{r'}{\alpha}}_{r'\frac{\beta}{\alpha}}(s-v)dv\Big)^{\frac{1}{r'}}\\
			&\leq C\Big(\|b\|_{L^r(B^{\beta}_{p,q})}+\|\Div(b)\|_{L^r(B^{\beta}_{p,q})}\Big)\\
			&\quad\cdot\Big(\int_0^s\big(\sup_{t\in(0,v]}w^{\gamma_2}_{\gamma_1}(t)\|(\brho^1_{\mu}-\brho^2_{\mu})(t,\cdot)\|_{B^{-\beta-\theta}_{p',1}}\big)^{r'}w^{-2r'\gamma_2}_{-2r'\gamma_1}(v)w^{-\frac{r'}{\alpha}}_{r'\frac{\beta}{\alpha}}(s-v)dv\Big)^{\frac{1}{r'}},
		\end{align*}
		and
		\begin{align*}
			&\sup_{s\in(0,S]}w^{\gamma_2}_{\gamma_1}(s)\|(\brho^1_{\mu}-\brho^2_{\mu})(s,\cdot)\|_{B^{-\beta-\theta}_{p',1}}\\
			&\leq C\Big(\|b\|_{L^r(B^{\beta}_{p,q})}+\|\Div(b)\|_{L^r(B^{\beta}_{p,q})}\Big)\\
			&\quad\cdot\Big(\int_0^s\big(\sup_{t\in(0,v]}w^{\gamma_2}_{\gamma_1}(t)\|(\brho^1_{\mu}-\brho^2_{\mu})(t,\cdot)\|_{B^{-\beta-\theta}_{p',1}}\big)^{r'}w^{-2r'\gamma_2}_{-2r'\gamma_1}(v)w^{-\frac{r'}{\alpha}}_{r'\frac{\beta}{\alpha}}(s-v)dv\Big)^{\frac{1}{r'}}.
		\end{align*}

		Then applying the Grönwall-Volterra lemma \cite[Lemma 2.2]{ZHANGGronwallVolterra} similarly as in Lemma \ref{lemma-density-conv-lt}, we deduce that $\brho^1_{\mu}=\brho^2_{\mu}$ in the given norm.

	\end{proof}
\end{lemma}

\subsection{Global SDE well-posedness}

\begin{lemma}[Integrability property of the mollified drift in long time]
	Assume that \eqref{mu-range}, \eqref{C3} and \eqref{C3LT} hold. Assume also that \eqref{good-relation} or \eqref{good-relation-enhaced} hold. Then for any $S\in(0,\mathcal{T}_1)$,
	\begin{align*}
		\mathcal{B}^\varepsilon_{\brho^\varepsilon_{\mu}}\in L^{r_\theta}((0,S],B^{-\theta}_{\infty,\infty})
	\end{align*}
	and
	\begin{align*}
		\|\mathcal{B}^\varepsilon_{\brho^\varepsilon_{\mu}}(s,\cdot)\|_{B^{-\theta}_{\infty,\infty}}\leq C\|b\|_{L^r(B^\beta_{p,q})}w^{\frac{1}{r_\theta}-\frac{1}{r}-\gamma_2}_{\frac{1}{r_\theta}-\frac{1}{r}-\gamma_1}(S).
	\end{align*}

	\begin{proof}
		We apply the Hölder inequality $L^{r_{\theta}}:L^r-L^{(r_\theta^{-1}-r^{-1})^{-1}}$ and use Lemma \ref{lemma-density-estimate-lt}:
		\begin{align*}
			\|\mathcal{B}^\varepsilon_{\brho^\varepsilon_{\mu}}\|_{L^{r_\theta}(B^{-\theta}_{\infty,\infty})} & \leq C\|b^\varepsilon\|_{L^r(B^\beta_{p,q})}\Big(\int_0^S\|\brho^\varepsilon_{\mu}(s,\cdot)\|_{B^{-\beta-\theta}_{p',1}}^{\frac{r_\theta r}{r-r_\theta}}ds\Big)^{\frac{r-r_\theta}{r_\theta r}}                                                                                            \\
			                                                                                                    & \leq C\|b^\varepsilon\|_{L^r(B^\beta_{p,q})}\|\brho^\varepsilon_{\mu}\|_{L^{\infty}_{w^{\gamma_2}_{\gamma_1}}((0,S],B^{-\beta-\theta}_{p',1})}\Big(\int_0^Sw^{-\gamma_2\frac{r_\theta r}{r-r_\theta}}_{-\gamma_1\frac{r_\theta r}{r-r_\theta}}( s )ds\Big)^{\frac{r-r_\theta}{r_\theta r}} \\
			                                                                                                    & \leq C\|b\|_{L^r(B^\beta_{p,q})}\|\brho^\varepsilon_{\mu}\|_{L^{\infty}_{w^{\gamma_2}_{\gamma_1}}((0,S],B^{-\beta-\theta}_{p',1})}w^{\frac{1}{r_\theta}-\frac{1}{r}-\gamma_2}_{\frac{1}{r_\theta}-\frac{1}{r}-\gamma_1}(S).
		\end{align*}

		Above, we have integrable singularities when $\gamma_1\frac{r_\theta r}{r-r_\theta}<1$, $\gamma_2\frac{r_\theta r}{r-r_\theta}<1$, meaning that $r_\theta$ should be such that $r_\theta<\min\Big(\frac{r}{1+r\gamma_1},\frac{r}{1+r\gamma_2}\Big)<r$.
	\end{proof}

	\begin{proposition}
		Let conditions \eqref{mu-range}, \eqref{C3} and \eqref{C3LT} be satisfied. Let $\theta,r$ satisfy \eqref{good-relation}. Denote by $(\Prob^\alpha_\varepsilon)_{\varepsilon>0}$ the solution to the non-linear martingale problem associated with \eqref{main-sde-mollified}. Then there exists a unique limit point of a converging subsequence $(\Prob^\alpha_{\varepsilon_k})_{k\geq1}$, $\varepsilon_k\downarrow_{k\to\infty}0$, in $\mathcal{P}(\Omega_\alpha)$ equipped with its weak topology, such that it solves the non-linear martingale problem associated with \eqref{main-sde} in the sense of Def. \ref{def-martingale-problem}.

		\begin{proof}
			We follow the strategy of the proof of Proposition \ref{prop-martingale-problem-wp}. Tightness of $(\Prob_\varepsilon^\alpha)$ follows from the proof of \ref{prop-martingale-problem-wp}. To identify the limit points, it is enough to observe that
			\begin{align*}
				\|\mathcal{B}^{\varepsilon_k}_{\brho^{\varepsilon_k}_{\mu}}-\mathcal{B}_{\brho_{\mu}}\|_{L^{r_\theta}((0,S],B^{-\theta-\varepsilon}_{\infty,\infty})} & \leq C\Big(\|\brho^{\varepsilon_k}_{\mu}\|_{L^{\infty}_{w^{\gamma_2}_{\gamma_1}}((0,S],B^{-\beta-\theta}_{p',1})}\|b^{\varepsilon_k}-b\|_{L^{\bar{r}}((0,S],B^{\beta-\varepsilon}_{p,q})}w^{\frac{1}{r_\theta}-\frac{1}{\bar{r}}-\gamma_2}_{\frac{1}{r_\theta}-\frac{1}{\bar{r}}-\gamma_1}(S) \\
				                                                                                                                                                       & \quad+\|b\|_{L^{r}((0,S],B^\beta_{p,q})}\|\brho_{\mu}^{\varepsilon_k}-\brho_{\mu}\|_{L^\infty_{w^{\gamma_2}_{\gamma_1}}((0,S],B^{-\beta-\theta}_{p',1})}w^{\frac{1}{r_\theta}-\frac{1}{r}-\gamma_2}_{\frac{1}{r_\theta}-\frac{1}{r}-\gamma_1}(S) \Big)                                                                            \\
				                                                                                                                                                       & \xrightarrow[k\to+\infty]{}0,
			\end{align*}
			and $\frac{1}{r_\theta}-\frac{1}{r}-\gamma_2>0$, $\frac{1}{r_\theta}-\frac{1}{r}-\gamma_1>0$.

			\vspace{2mm}

			Uniqueness of the limit is proved similarly, noting that
			\begin{align*}
				\|\mathcal{B}_{\brho_{\mu}}\|_{L^{r_\theta}((0,S],B^{-\theta}_{\infty,\infty})}\leq C\|b\|_{L^r(B^\beta_{p,q})}\|\brho_{\mu}\|_{L^\infty_{w^{\gamma_2}_{\gamma_1}}((0,S],B^{-\beta-\theta}_{p',1})}<+\infty,
			\end{align*}
			where $\brho_{\mu}\in L^\infty_{w^{\gamma_2}_{\gamma_1}}((0,S],B^{-\beta-\theta}_{p',1})$ is a unique solution to \eqref{pde-main}.
		\end{proof}
	\end{proposition}

\end{lemma}

The rest of Theorem \ref{thm-main-sde-lt}, i.e., weak well-posedness, stays the same as in Theorem \ref{thm-main-sde} once the martingale problem well-posedness is provided.

\section{Classical well-posedness in short time}
\label{section-sde-lplq}

This section aims at showing that the results of work \cite{deraynal2023multidimensionalstabledrivenmckeanvlasov} hold not only for $\beta=-1$ but for $\beta\in(-2,-1]$ provided the initial condition is sufficiently smooth. We insist that the result proven in this section allows us to understand a weak solution (if it exists) in a classical sense; however, this requires more regularity on the initial measure compared to Theorem \ref{thm-main-sde}.

\vspace{1mm}

Recall that
\begin{align*}
	\gamma^*=\frac{1}{\alpha}\Big(-\beta+\frac{d}{p}-\zeta_0+\Big(\frac{1+\eta}{2\eta}\Big)\Gamma\Big)>0,
\end{align*}
and $\Gamma\in(0,1)$ is given by
\begin{align}
	\label{def-lplq-Gamma}
	\Gamma=\eta\Big(\alpha-1+\beta-\frac{\alpha}{r}-\frac{d}{p}+\zeta_0\Big),\quad \eta\in(0,1).
\end{align}

As before, we prove an a priori result for the non-linear Fokker-Planck equation.
\begin{proposition}[Well-posedness of the Fokker-Planck PDE]
	\label{prop-lplq-main-pde}
	Assume that condition \eqref{C4} holds for $\beta\in(-2,-1]$. Then, for any $\mu\in\mathcal{P}(\R^d)\cap B^{\beta_0}_{p_0,q_0}$, there exists a time horizon $\mathcal{T}_0\in(0,T]$ small enough such that the non-linear Fokker-Planck equation \eqref{pde-main} admits a unique mild solution in $L^\infty_{\gamma^*}((0,S],B^{-\beta+\Gamma}_{p',1})$ for any $S\in(0,\mathcal{T}_0]$ with $\Gamma,\theta$ defined in Theorem \ref{thm-main-sde-lplq}.

	\vspace{1mm}
	\noindent
	Moreover, for all $s\in[0,S]$, $\brho_{\mu}(s,\cdot)\in\mathcal{P}(\R^d)$, and for a.e. $s\in(0,S]$, $\brho_{\mu}(s,\cdot)$ is absolutely continuous w.r.t. the Lebesgue measure and it satisfies the Duhamel representation \eqref{pde-main}.
\end{proposition}

The main difference between Theorem \ref{thm-main-sde-lplq} and the results of \cite{deraynal2023multidimensionalstabledrivenmckeanvlasov} is the a priori estimates on the mollified density, which is described in the following Lemma.

\begin{lemma}[A priori estimates on the mollified density]
	\label{lemma-lplq-density-estimate}
	Let condition \eqref{C4} be satisfied. Then for any small $S\leq T$,
	\begin{align*}
		 & \sup_{s\in(0,S]} s ^{\gamma^*}\|\brho^\varepsilon_{\mu}(s,\cdot)\|_{B^{-\beta+\Gamma}_{p',1}}\leq C\|\mu\|_{B^{\bar{\beta}_0}_{\bar{p}_0,\bar{q}_0}}S^{\frac{1-\eta}{2\eta}\frac{\Gamma}{\alpha}}                                   \\
		 & \quad+C\Big(\|b\|_{L^r(B^\beta_{p,q})}+\|\Div(b)\|_{L^r(B^\beta_{p,q})}\Big)\big(\sup_{s\in(0,S]}s^{\gamma^*}\|\brho^\varepsilon_{\mu}(s,\cdot)\|_{B^{-\beta+\Gamma}_{p',1}}\big)^2S^{\frac{1-\eta}{2\eta}\frac{\Gamma}{\alpha}}.
	\end{align*}

	\begin{proof}
		Let $\beta\in(-2,-1]$. We start with the Duhamel representation \eqref{duhamel-main-mollified} under \eqref{C4} for fixed $s\in(0,S]$ by applying integration by parts in the bilinear term.
		\begin{align}
			\label{prf-ineq-lplq-duh}
			\|\brho^\varepsilon_{\mu}(s,\cdot)\|_{B^{-\beta+\Gamma}_{p',1}} & \leq\|\mu\ast p^\alpha_{s}\|_{B^{-\beta+\Gamma}_{p',1}}+\int_0^s\|\big(\Div(\mathcal{B}^\varepsilon_{\brho^\varepsilon_{\mu}})(v,\cdot)\brho^\varepsilon_{\mu}(v,\cdot)\big)\ast p^\alpha_{s-v}\|_{B^{-\beta+\Gamma}_{p',1}}dv \\
			                                                                  & \quad+\int_0^s\|\big(\mathcal{B}^\varepsilon_{\brho^\varepsilon_{\mu}}(v,\cdot)\cdot\nabla\brho^\varepsilon_{\mu}(v,\cdot)\big)\ast p^\alpha_{s-v}\|_{B^{-\beta+\Gamma}_{p',1}}dv.
		\end{align}

		The initial data term $\|\mu\ast p^\alpha_{s}\|$ is handled in the same way as in Lemma \ref{lemma-density-estimate} by replacing $-\theta$ with $\Gamma$ and noting that the estimate does not depend on the range of $\beta$:
		\begin{align*}
			\|\mu\ast p^\alpha_{s}\|_{B^{-\beta+\Gamma}_{p',1}}\leq C\|\mu\|_{B^{\bar{\beta}_0}_{\bar{p}_0,\bar{q}_0}} v^{-\gamma_0^*},
		\end{align*}
		where
		\begin{align*}
			\gamma_0^*=\frac{1}{\alpha}\Big(\Gamma-\beta+\frac{d}{p}-\zeta_0\Big)_+.
		\end{align*}

		\vspace{1mm}

		Now, note that in Lemma 6 of \cite{deraynal2023multidimensionalstabledrivenmckeanvlasov}, the computations are made with a fixed $\beta=-1$. Here, we do not fix $\beta$, and instead we show that it can vary over the interval $(-2,-1]$. Let us handle the integrand of the first integral:
		\begin{align*}
			\|\big(\Div(\mathcal{B}^\varepsilon_{\brho^\varepsilon_{\mu}})(v,\cdot)\brho^\varepsilon_{\mu}(v,\cdot)\big)\ast p^\alpha_{s-v}\|_{B^{-\beta+\Gamma}_{p',1}} & \stackrel{\eqref{besov-prop-y}}{\leq} C \|\Div(\mathcal{B}^\varepsilon_{\brho^\varepsilon_{\mu}}(v,\cdot))\brho^\varepsilon_{\mu}(v,\cdot)\|_{B^{\Gamma}_{p',\infty}}\|p^\alpha_{s-v}\|_{B^{-\beta}_{1,1}}                          \\
			                                                                                                                                                                 & \stackrel{\eqref{besov-prop-pr1}}{\leq} C\|\Div(\mathcal{B}^\varepsilon_{\brho^\varepsilon_{\mu}}(v,\cdot))\|_{B^\Gamma_{\infty,\infty}}\|\brho^\varepsilon_{\mu}(v,\cdot)\|_{B^{\Gamma}_{p',1}}\|p^\alpha_{s-v}\|_{B^{-\beta}_{1,1}} \\
			                                                                                                                                                                 & \stackrel{\eqref{besov-prop-y}, \eqref{besov-prop-e2}, \eqref{besov-prop-hk}}{\leq} C\|\Div(b^\varepsilon(v,\cdot))\|_{B^\beta_{p,q}}\|\brho^\varepsilon_{\mu}(v,\cdot)\|^2_{B^{-\beta-\Gamma}_{p',1}}(s-v)^{\frac{\beta}{\alpha}}.
		\end{align*}

		Similarly, for the second integral in \eqref{prf-ineq-lplq-duh} using additionally \eqref{besov-prop-l} for the third inequality,
		\begin{align*}
			\|\big(\mathcal{B}^\varepsilon_{\brho^\varepsilon_{\mu}}(v,\cdot)\cdot\nabla\brho^\varepsilon_{\mu}(v,\cdot)\big)\ast p^\alpha_{s-v}\|_{B^{-\beta+\Gamma}_{p',1}} & \stackrel{\eqref{besov-prop-y}}{\leq} C\|\mathcal{B}^\varepsilon_{\brho^\varepsilon_{\mu}}(v,\cdot)\cdot\nabla\brho^\varepsilon_{\mu}(v,\cdot)\|_{B^{\Gamma}_{p',\infty}}\|p^\alpha_{s-v}\|_{B^{-\beta}_{1,1}}                                  \\
			                                                                                                                                                                      & \stackrel{\eqref{besov-prop-pr1}}{\leq} C\|\mathcal{B}^\varepsilon_{\brho^\varepsilon_{\mu}}(v,\cdot)\|_{B^\Gamma_{\infty,\infty}}\|\nabla\brho^\varepsilon_{\mu}(v,\cdot)\|_{B^\Gamma_{p',1}}\|p^\alpha_{s-v}\|_{B^{-\beta}_{1,1}}               \\
			                                                                                                                                                                      & \stackrel{\eqref{besov-prop-l}}{\leq} C\|b^\varepsilon(v,\cdot)\|_{B^\beta_{p,q}}\|\brho^\varepsilon_{\mu}(v,\cdot)\|_{B^{-\beta+\Gamma}_{p',1}}\|\brho^\varepsilon_{\mu}(v,\cdot)\|_{B^{1+\Gamma}_{p',1}}\|p^\alpha_{s-v}\|_{B^{-\beta}_{1,1}} \\
			                                                                                                                                                                      & \stackrel{\eqref{besov-prop-e2}, \eqref{besov-prop-hk}}{\leq} C\|b^\varepsilon(v,\cdot)\|_{B^\beta_{p,q}}\|\brho^\varepsilon_{\mu}(v,\cdot)\|^2_{B^{-\beta+\Gamma}_{p',1}}(s-v)^{\frac{\beta}{\alpha}},
		\end{align*}
		where the last inequality is possible from \eqref{besov-prop-e2} because $\beta\leq-1$.

		\vspace{1mm}

		Note that the main difference in handling this quadratic term compared to Lemma \ref{lemma-density-estimate} is the product rule that is used. Indeed, here we apply \eqref{besov-prop-pr1} for positive regularity, which is possible thanks to the parameter $\Gamma>0$, whereas in Lemma \ref{lemma-density-estimate}, $\Gamma$ was replaced by $-\theta$, and thus the product rule \eqref{ineq-product-rule} was used.

		\vspace{1mm}

		Now, substituting the obtained bounds in the Duhamel representation,
		\begin{align*}
			\|\brho^\varepsilon_{\mu}(s,\cdot)\|_{B^{-\beta+\Gamma}_{p',1}} & \leq C\|\mu\|_{B^{\bar{\beta}_0}_{\bar{p}_0,\bar{q}_0}}s^{-\gamma_0^*}                                                                                                                                                 \\
			                                                                  & \quad+C\int_0^s\big(\|b^\varepsilon(v,\cdot)\|_{B^{\beta}_{p,q}}+\|\Div(b^\varepsilon(v,\cdot))\|_{B^{\beta}_{p,q}}\big)\|\brho^\varepsilon_{\mu}(v,\cdot)\|^2_{B^{-\beta+\Gamma}_{p',1}}(s-v)^{\frac{\beta}{\alpha}}dv.
		\end{align*}

		Repeating the same steps as in Lemma \ref{lemma-density-estimate}, we obtain
		\begin{align*}
			 & \sup_{s\in(0,S]}s^{\gamma^*}\|\brho^\varepsilon_{\mu}(s,\cdot)\|_{B^{-\beta+\Gamma}_{p',1}}\leq C\|\mu\|_{B^{\bar{\beta}_0}_{\bar{p}_0,\bar{q}_0}}S^{\gamma^*-\gamma_0^*}                                                                                              \\
			 & \quad+\Big(\|b\|_{L^r(B^{\beta}_{p,q})}+\|\Div(b)\|_{L^r(B^{\beta}_{p,q})}\Big)\big(\sup_{s\in(0,S]}s^{\gamma^*}\|\brho^\varepsilon_{\mu}(s,\cdot)\|_{B^{-\beta+\Gamma}_{p',1}}\big)^2\Big(\int_0^S(s-v)^{\frac{r'\beta}{\alpha}} v^{-2r'\gamma^*}dv\Big)^{\frac{1}{r'}},
		\end{align*}
		where the integral is finite if $\gamma^*<\frac{1}{2r'}$.
		\begin{align}
			\label{prf-lplq-beta-cond}
			\frac{1}{r}<1+\frac{\beta}{\alpha}\quad\text{ and }\quad\beta>-\alpha.
		\end{align}
		This leads to
		\begin{align*}
			 & \sup_{s\in(0,S]}s^{\gamma^*}\|\brho^\varepsilon_{\mu}(s,\cdot)\|_{B^{-\beta+\Gamma}_{p',1}}\leq C\|\mu\|_{B^{\bar{\beta}_0}_{\bar{p}_0,\bar{q}_0}}S^{\gamma^*-\gamma_0^*}                                                             \\
			 & \quad+\Big(\|b\|_{L^r(B^{\beta}_{p,q})}+\|\Div(b)\|_{L^r(B^{\beta}_{p,q})}\Big)\big(\sup_{s\in(0,S]}s^{\gamma^*}\|\brho^\varepsilon_{\mu}(s,\cdot)\|_{B^{-\beta+\Gamma}_{p',1}}\big)^2S^{\frac{1}{r'}+\frac{\beta}{\alpha}-\gamma^*}.
		\end{align*}

		In \cite{deraynal2023multidimensionalstabledrivenmckeanvlasov}, condition \eqref{prf-lplq-beta-cond} was obtained with $\beta=-1$, which led to a more relaxed condition on $r$. Here, however, we can take $\beta$ close to $-\alpha$ the price of which requires time integrability $r$ of the kernel to be close to infinity. The rest of the proof is the same up to the constraint \eqref{prf-lplq-beta-cond}. More precisely, we find $\Gamma,\theta$ from the conditions
		\begin{align*}
			\gamma_0^*=\frac{1}{\alpha}\Big(\Gamma-\beta+\frac{d}{p}-\zeta_0\Big)_+<\gamma^*<\frac{1}{r'}+\frac{\beta}{\alpha},
		\end{align*}
		which leads to
		\begin{align*}
			\beta>-\alpha+\frac{\alpha}{r}+\Big(-\beta+\frac{d}{p}-\zeta_0\Big)_+,
		\end{align*}
		and this is exactly \eqref{C4}.
	\end{proof}
\end{lemma}

With Lemma \ref{lemma-lplq-density-estimate}, the rest of the proof of Proposition \ref{prop-lplq-main-pde} stays the same as in \cite{deraynal2023multidimensionalstabledrivenmckeanvlasov} as well as the weak well-posedness of \eqref{main-sde} under \eqref{C4}.

\begin{remark}[Strong well-posedness]
	The strong well-posedness of \eqref{main-sde} under \eqref{C4S} follows from Proposition 15 in \cite{deraynal2023multidimensionalstabledrivenmckeanvlasov} by replacing $\beta=-1$ by $\beta\in(-2,-1]$. The only thing that appears at the end of the proof of Proposition 15 in \cite{deraynal2023multidimensionalstabledrivenmckeanvlasov} to be checked is
	\begin{align*}
		\frac{\Gamma}{\bar{\eta}}>1-\frac{\alpha}{2}.
	\end{align*}
	From \eqref{def-lplq-Gamma}, we obtain condition \eqref{C4S}. More precisely,
	\begin{align*}
		\frac{\Gamma}{\bar{\eta}}=\alpha+\beta-\frac{\alpha}{r}-\Big(-\beta+\frac{d}{p}-\zeta_0\Big)-\eta>1-\frac{\alpha}{2} \\
		\iff\beta>\frac{1}{2}\Big(1-\frac{3}{2}\alpha+\frac{\alpha}{r}+\frac{d}{p}-\zeta_0\Big).
	\end{align*}
\end{remark}

\printbibliography

@article{zvonkin,
author = {Zvonkin, Alexander},
year = {2007},
month = {10},
pages = {129},
title = {A transformation of the phase space of a diffusion process that removes the drift},
volume = {22},
journal = {Mathematics of the USSR-Sbornik},
doi = {10.1070/SM1974v022n01ABEH001689}}

@misc{deraynal2022multidimensional,
author = {Chaudru de Raynal, Paul-Eric and Jabir, Jean-François and Menozzi, Stéphane},
year = {2022},
month = {05},
title = {Multidimensional Stable driven McKean-Vlasov SDEs with distributional interaction kernel -- a regularization by noise perspective},
doi = {10.48550/arXiv.2205.11866}}

@article{deraynal2022-linear-besov,
author = {Chaudru de Raynal, Paul-Eric and Menozzi, Stéphane},
title = {{On multidimensional stable-driven stochastic differential equations with Besov drift}},
volume = {27},
journal = {Electronic Journal of Probability},
publisher = {Institute of Mathematical Statistics and Bernoulli Society},
pages = {1 -- 52},
year = {2022},
doi = {10.1214/22-EJP864},
URL = {https://doi.org/10.1214/22-EJP864}}

@article{Veretennikov1981,
title={On strong solutions and explicit formulas for solutions of stochastic integral equations},
author={Veretennikov, Alexander},
journal={Mathematics of The Ussr-sbornik},
year={1981},
volume={39},
pages={387-403},
url={10.1070/SM1981v039n03ABEH001522}}

@book{Sawano2018BesovSpaces,
title={Theory of Besov Spaces},
author={Sawano, Yoshihiro},
year={2018},
publisher={Springer Singapore},
url={https://doi.org/10.1007/978-981-13-0836-9}}

@article{Burenkov1990BesovYoungIneq,
Author = {Burenkov, Victor},
Title = {On estimates of Fourier transforms and convolutions in {Nikol}'skij-Besov spaces},
FJournal = {Proceedings of the Steklov Institute of Mathematics},
Journal = {Proc. Steklov Inst. Math.},
ISSN = {0081-5438},
Volume = {187},
Pages = {35--44},
Year = {1990},
zbMATH = {4186334},
Zbl = {0719.42021}}

@book{LemariRieusset2002RecentDI,
title={Recent Developments in the Navier-Stokes Problem},
author={Lemari{\'e}–Rieusset, Pierre Gilles},
year={2002},
url={https://doi.org/10.1201/9780367801656},
publisher={CRC Press}}

@book{Sato1999LvyPA,
title={L{\'e}vy Processes and Infinitely Divisible Distributions},
author={Sato, Ken{-}iti},
year={1999},
publisher={New York: Cambridge University Press},
url={https://api.semanticscholar.org/CorpusID:209888537}}

@article{Kolokoltsov2000SymmetricSL,
title={Symmetric Stable Laws and Stable‐Like Jump‐Diffusions},
author={Kolokoltsov, Vassili},
journal={Proceedings of the London Mathematical Society},
year={2000},
volume={80},
url={https://doi.org/10.1112/S0024611500012314}}

@book{RunstSickel1996,
address = {Berlin, New York},
author = {Runst, Thomas and Sickel, Winfried},
doi = {doi:10.1515/9783110812411},
isbn = {9783110812411},
lastchecked = {2024-08-13},
publisher = {De Gruyter},
title = {Sobolev Spaces of Fractional Order, Nemytskij Operators, and Nonlinear Partial Differential Equations},
url = {https://doi.org/10.1515/9783110812411},
year = {1996}}

@article{flandoli2015multidimensionalsde,
author = {Flandoli, Franco and Issoglio, Elena and Russo, Francesco},
year = {2017},
month = {04},
pages = {},
title = {Multidimensional stochastic differential equations with distributional drift},
journal = {Transactions of the American Mathematical Society},
doi = {10.1090/tran/6729}}

@article{Priola2012PathwiseUniq,
author = {Priola, Enrico},
journal = {Osaka Journal of Mathematics},
number = {2},
pages = {421 -- 447},
publisher = {Osaka University and Osaka Metropolitan University, Departments of Mathematics},
title = {{Pathwise uniqueness for singular SDEs driven by stable processes}},
volume = {49},
year = {2012}}

@article{ABM2020StableSDE,
author = {Athreya, Siva and Butkovsky, Oleg and Mytnik, Leonid},
doi = {10.1214/19-AOP1358},
journal = {The Annals of Probability},
number = {1},
pages = {178 -- 210},
publisher = {Institute of Mathematical Statistics},
title = {{Strong existence and uniqueness for stable stochastic differential equations with distributional drift}},
url = {https://doi.org/10.1214/19-AOP1358},
volume = {48},
year = {2020}}

@book{stroock1997multidimensional,
author = {Stroock, Daniel and Varadhan, Srinivasa},
isbn = {9783540903536},
lccn = {96027093},
publisher = {Springer Berlin Heidelberg},
series = {Grundlehren der mathematischen Wissenschaften},
title = {Multidimensional Diffusion Processes},
url = {https://books.google.de/books?id=DuDsmoyqCy4C},
year = 1997}

@article{DelarueDiel2016RoughPaths,
author = {Delarue, François and Diel, Roland},
year = {2016},
month = {02},
pages = {},
title = {Rough paths and 1d SDE with a time dependent distributional drift: application to polymers},
volume = {165},
journal = {Probability Theory and Related Fields},
doi = {10.1007/s00440-015-0626-8}}

@article{KrylovRockner2005LpLq,
author = {Krylov, Nicolai and Röckner, Michael},
year = {2005},
month = {02},
pages = {154-196},
title = {Strong solutions of stochastic equations with singular time dependent drift},
volume = {131},
journal = {Probab. Theory Related Fields},
doi = {10.1007/s00440-004-0361-z}}

@article{XieZhang2020LpLq,
author = {Xie, Longjie and Zhang, Xicheng},
doi = {10.1214/19-AIHP959},
journal = {Annales de l'Institut Henri Poincar{\'e}, Probabilit{\'e}s et Statistiques},
number = {1},
pages = {175 -- 229},
publisher = {Institut Henri Poincar{\'e}},
title = {{Ergodicity of stochastic differential equations with jumps and singular coefficients}},
url = {https://doi.org/10.1214/19-AIHP959},
volume = {56},
year = {2020}}

@article{KrempPerkowski2022,
author = {Kremp, Helena and Perkowski, Nicolas},
year = {2022},
month = {08},
pages = {},
title = {Multidimensional SDE with distributional drift and Lévy noise},
volume = {28},
journal = {Bernoulli},
doi = {10.3150/21-BEJ1394}}

@article{KurtzWeakStrongSolutions,
author = {Thomas Kurtz},
title = {{Weak and strong solutions of general stochastic models}},
volume = {19},
journal = {Electronic Communications in Probability},
number = {none},
publisher = {Institute of Mathematical Statistics and Bernoulli Society},
pages = {1 -- 16},
year = {2014},
doi = {10.1214/ECP.v19-2833},
URL = {https://doi.org/10.1214/ECP.v19-2833}}

@misc{deraynal2023multidimensionalstabledrivenmckeanvlasov,
title={Multidimensional Stable Driven McKean-Vlasov SDEs with Distributional Interaction Kernel: Critical Thresholds and Related Models}, 
author={Chaudru de Raynal, Paul-Eric and Jabir, Jean-François and Menozzi, Stéphane},
year={2023},
eprint={2302.09900},
archivePrefix={arXiv},
primaryClass={math.AP},
url={https://arxiv.org/abs/2302.09900}}

@book{billing_prob_convergenc,
address = {New York},
author = {Billingsley, Patrick},
description = {q-paper},
edition = {Second},
isbn = {0-471-19745-9},
mrclass = {60B10 (28A33 60F17)},
mrnumber = {MR1700749 (2000e:60008)},
note = {A Wiley-Interscience Publication},
pages = {x+277},
publisher = {John Wiley \& Sons Inc.},
series = {Wiley Series in Probability and Statistics: Probability and Statistics},
title = {Convergence of probability measures},
year = 1999}

@book{markov-processes-conv,
address = {New York},
author = {Ethier, Stewart N. and Kurtz, Thomas G.},
description = {MR: Publications results for "MR Number=(838085)"},
isbn = {0-471-08186-8},
mrclass = {60J25 (60B10 60F05 60F17 60G44 60J80)},
mrnumber = {MR838085 (88a:60130)},
mrreviewer = {S. R. S. Varadhan},
pages = {x+534},
publisher = {John Wiley \& Sons Inc.},
series = {Wiley Series in Probability and Mathematical Statistics: Probability and Mathematical Statistics},
title = {Markov processes -- characterization and convergence},
year = 1986}

@book{Bass_2011-stoch-pr, 
place={Cambridge}, 
series={Cambridge Series in Statistical and Probabilistic Mathematics}, 
title={Stochastic Processes}, 
publisher={Cambridge University Press}, 
author={Bass, Richard F.}, 
year={2011}, 
collection={Cambridge Series in Statistical and Probabilistic Mathematics}}

@article{WatanabeSingularMeasures,
author = {Watanabe, Toshiro},
year = {2007},
month = {01},
pages = {2851-2879},
title = {Asymptotic estimates of multi-dimensional stable densities and their applications},
volume = {359},
journal = {Transactions of the American Mathematical Society},
doi = {10.1090/S0002-9947-07-04152-9}}

@book{BanachSpaces,
title = "Analysis in Banach Spaces: Volume I: Martingales and Littlewood-Paley Theory",
author = "Hyt{\"o}nen, Tuomas and {van Neerven}, Jan and Veraar, Mark and Weis, Lutz",
year = "2016",
doi = "10.1007/978-3-319-48520-1",
language = "English",
isbn = "978-3-319-48519-5",
series = "Ergebnisse der Mathematik und ihrer Grenzgebiete. 3. Folge",
publisher = "Springer"}

@misc{rockner2025sdescriticaltimedependent,
title={SDEs with critical time dependent drifts: strong solutions}, 
author={Röckner, Michael and Zhao, Guohuan},
year={2025},
eprint={2103.05803},
archivePrefix={arXiv},
primaryClass={math.PR},
url={https://arxiv.org/abs/2103.05803}}

@misc{krylov21pd,
author={Krylov, Nicolai},
year={2020},
month={07},
title={On strong solutions of It\^o's equations with a$\,\in W^{1}_{d}$ and b$\,\in L_{d}$},
doi={10.48550/arXiv.2007.06040}}

@misc{fitoussi2024weakwellposednessweakdiscretization,
title={Weak well-posedness and weak discretization error for stable-driven SDEs with Lebesgue drift}, 
author={Fitoussi, Mathis and Jourdain, Benjamin and Menozzi, Stéphane},
year={2024},
eprint={2405.08378},
archivePrefix={arXiv},
primaryClass={math.PR},
url={https://arxiv.org/abs/2405.08378}}

@article{ZHANGGronwallVolterra,
title = {Stochastic Volterra equations in Banach spaces and stochastic partial differential equation},
journal = {Journal of Functional Analysis},
volume = {258},
number = {4},
pages = {1361-1425},
year = {2010},
issn = {0022-1236},
doi = {https://doi.org/10.1016/j.jfa.2009.11.006},
url = {https://www.sciencedirect.com/science/article/pii/S0022123609004650},
author = {Xicheng Zhang}}

\end{document}